\documentclass[11pt]{amsart}
\usepackage{amsmath,amssymb,amsthm,mathtools,mathrsfs}
\usepackage{enumitem}
\usepackage[margin=1.08in]{geometry}
\usepackage{array}
\allowdisplaybreaks[2]
\usepackage{booktabs}
\usepackage[protrusion=true,expansion=false]{microtype}
\usepackage{xcolor}
\usepackage{hyperref}
\hypersetup{colorlinks=true,citecolor=blue!50!black,linkcolor=blue!50!black,urlcolor=blue!50!black,pdfauthor={Hamed Ebadi},pdftitle={Stationary Phase along Rays and the Determinant of Linearly Constrained Kloosterman Families}}
\numberwithin{equation}{section}
\theoremstyle{plain}
\newtheorem{theorem}[equation]{Theorem}
\newtheorem{proposition}[equation]{Proposition}
\newtheorem{lemma}[equation]{Lemma}
\newtheorem{corollary}[equation]{Corollary}
\newtheorem*{theoremA}{Theorem A}
\newtheorem*{theoremB}{Theorem B}
\newtheorem*{theoremC}{Theorem C}
\newtheorem*{theoremD}{Theorem D}
\newtheorem*{corollaryE}{Corollary E}
\newtheorem*{theorem0}{Theorem 0}
\theoremstyle{definition}
\newtheorem{definition}[equation]{Definition}
\newtheorem{example}[equation]{Example}
\theoremstyle{remark}
\newtheorem{remark}[equation]{Remark}

\newcommand{\Qlb}{\overline{\mathbb{Q}}_\ell}

\newcommand{\Fq}{\mathbb{F}_q}
\newcommand{\Fp}{\mathbb{F}_p}
\newcommand{\Fpbar}{\overline{\mathbb{F}}_p}
\newcommand{\Fqbar}{\overline{\mathbb{F}}_q}
\newcommand{\Gm}{\mathbb{G}_m}
\newcommand{\Gmn}[1]{\mathbb{G}_{m,#1}}
\newcommand{\A}{\mathbb{A}}
\newcommand{\PP}{\mathbb{P}}
\newcommand{\Z}{\mathbb{Z}}
\newcommand{\Q}{\mathbb{Q}}
\newcommand{\C}{\mathbb{C}}
\newcommand{\Frob}{\mathrm{Frob}}
\newcommand{\Tr}{\mathrm{Tr}}
\newcommand{\Swan}{\mathrm{Swan}}
\newcommand{\ggeom}{G_{\mathrm{geom}}}
\newcommand{\garith}{G_{\mathrm{arith}}}
\newcommand{\pgeom}{\pi_1^{\mathrm{geom}}}
\newcommand{\Lpsi}[1]{\mathcal{L}_\psi(#1)}
\newcommand{\Lchi}[1]{\mathcal{L}_{\chi_2}(#1)}
\newcommand{\Kl}{\mathcal{K}\ell}
\newcommand{\cW}{\mathcal{W}}
\newcommand{\cH}{\mathcal{H}}
\newcommand{\cE}{\mathcal{E}}
\newcommand{\cN}{\mathcal{N}}
\newcommand{\cK}{\mathcal{K}}
\newcommand{\cF}{\mathcal{F}}
\newcommand{\cG}{\mathcal{G}}
\newcommand{\cL}{\mathcal{L}}
\newcommand{\Msf}{\mathsf{M}}
\newcommand{\cA}{\mathcal{A}}
\newcommand{\cD}{\mathcal{D}}
\newcommand{\FT}{\mathrm{FT}_\psi}
\newcommand{\Floc}{F^{(\infty,0')}}
\newcommand{\Flocz}{F^{(0,\infty')}}
\newcommand{\Flocii}{F^{(\infty,\infty')}}
\newcommand{\Up}[1]{\mathrm{U}(#1)}
\newcommand{\pr}{\mathrm{pr}}
\newcommand{\sw}{\mathrm{sw}}
\DeclareMathOperator{\SL}{SL}
\DeclareMathOperator{\GL}{GL}
\DeclareMathOperator{\Sp}{Sp}

\DeclareMathOperator{\Ind}{Ind}
\DeclareMathOperator{\Res}{Res}
\DeclareMathOperator{\rank}{rank}
\DeclareMathOperator{\Gal}{Gal}
\DeclareMathOperator{\Pf}{Pf}
\DeclareMathOperator{\gr}{gr}

\title[Stationary phase along rays]{Stationary Phase along Rays and the Determinant of Linearly Constrained Kloosterman Families}
\author{Hamed Ebadi}
\address{K. N. Toosi University of Technology, Tehran, Iran}
\email{h.ebadi1@email.kntu.ac.ir}
\date{12 November 2026}

\begin{document}
\begin{abstract}
We study the top-weight sheaf $\mathcal W_n$ of linearly constrained Kloosterman sums over finite fields of odd characteristic, of rank $r(n)=2^n-\binom n{\lfloor n/2\rfloor}$. A local Fourier-transform analysis describes wild monodromy at infinity along rays through nonescaping critical points and their quadratic characters, giving explicit Swan conductors and the escape locus. The resulting stationary-phase description yields a ray-theoretic proof of the geometric determinant. A reflection pairing and a Frobenius calculation at the tame endpoint recover the arithmetic determinant $\det\mathcal W_n\simeq\mathcal L_\psi(-2^{n-1}\lambda^{-1}\sum_i h_i)\otimes\overline{\mathbb Q}_\ell(-(n+1)r(n)/2)$. We obtain reciprocity for exterior-power Frobenius traces, the exact distribution of normalized determinants, and spectral compression on symmetric loci. Assuming the ambient connected-monodromy theorem from companion work, we determine the monodromy of the geometric generic ray. Further results include bounds for complete exponential sums, conductor identities for selected parameter curves, and the Euler characteristic of the parameter torus. Rank-two boundary conormal calculations are separated from higher-codimension characteristic-cycle recursions, which remain conditional on explicit microlocal hypotheses.
\end{abstract}

\maketitle

\section{Introduction}\label{sec:intro}

\subsection{The family}\label{subsec:family}
Throughout, $p$ is an odd prime, $\ell\ne p$ is a prime, $\psi$ is a nontrivial additive character of $\Fp$ with values in $\Qlb^\times$, and $\lambda\in\Fp^\times$. For a finite extension $\Fq/\Fp$ we write $\psi_q=\psi\circ\Tr_{\Fq/\Fp}$. For $n\ge1$ and $h\in\Fq^n$ put
\begin{equation}\label{eq:Sigma}
\Sigma_n(h)=\Sigma_{n,q}(h):=\sum_{\substack{u\in(\Fq^\times)^n\\ \langle h,u\rangle=-\lambda}}\psi_q\bigl(u_1^{-1}+\cdots+u_n^{-1}\bigr),\qquad\langle h,u\rangle=\sum_{i=1}^nh_iu_i,
\end{equation}
and
\[
r(n)=2^n-\binom n{\lfloor n/2\rfloor},\qquad b_n=\binom n{\lfloor n/2\rfloor},\qquad N(n)=\frac{(n+1)\,r(n)}2,
\]
\[
\nu_{n,k}=\binom nk-\binom n{k-1}\qquad\Bigl(0\le k\le\Bigl\lfloor\frac n2\Bigr\rfloor\Bigr),
\]
with $\binom n{-1}=0$. The number $N(n)$ is an integer. Indeed, if $n$ is even, $b_n=2\binom{n-1}{n/2-1}$ is even, hence so is $r(n)$. Inserting $\mathbf 1_{\langle h,u\rangle=-\lambda}=q^{-1}\sum_{t\in\Fq}\psi_q(t(\lambda+\langle h,u\rangle))$ into \eqref{eq:Sigma} gives
\begin{equation}\label{eq:t-sum}
q\,\Sigma_n(h)=\sum_{t\in\Fq}\psi_q(\lambda t)\prod_{i=1}^nK(th_i),\qquad K(a)=K_q(a)=\sum_{x\in\Fq^\times}\psi_q(x^{-1}+ax),
\end{equation}
Thus $\Sigma_n(h)$ is an additively twisted average of products of Kloosterman sums at proportional arguments. Let $S_n=\Gm^n$ over $\Fp$, with coordinates $h=(h_1,\dots,h_n)$. The following theorem is established in the companion preprints \cite{EbMom,EbCoh}; we recall the construction in Section~\ref{sec:prelim}.

\begin{theorem0}[{\cite[Theorem~2.10, Corollary~2.16]{EbMom}, \cite[Theorem~5.3]{EbCoh}}]
For every $n\ge1$ there is a lisse $\Qlb$-sheaf $\cW_n$ on $S_n$, pure of weight $n+1$ and of rank $r(n)$, such that for every finite extension $\Fq/\Fp$ and every $h\in S_n(\Fq)$
\begin{equation}\label{eq:trace}
(-1)^{n-1}q\,\Sigma_n(h)=\Tr\bigl(\Frob_h\mid\cW_n\bigr)+\sum_{k=1}^{\lfloor n/2\rfloor}\nu_{n,k}\,q^k .
\end{equation}
The fibre of $\cW_n$ at $h$ is the weight-$(n-1)$ part of $H^{n-1}_c\bigl(Y_h\otimes\Fqbar,\Lpsi{\sum_iu_i^{-1}}\bigr)$, twisted by $(-1)$, where $Y_h=\{u\in\Gm^n:\langle h,u\rangle=-\lambda\}$.
\end{theorem0}

For $n\ge2$ and $p>\max(2,n)$, the range used in the 169-page parent manuscript, the connected-monodromy theorem from the companion manuscripts gives $\ggeom(\cW_n)^0=\SL_{r(n)}$ \cite[Theorem~8.8]{EbCoh}, \cite[Theorem~4.19]{EbMom}. No extension of this quoted group-theoretic result to the remaining small-characteristic cases is asserted here. In \cite[Theorem~E]{EbCoh} the determinant of $\cW_n$ is computed,
\begin{equation}\label{eq:det-intro}
\det\cW_n\cong\Lpsi{-2^{n-1}\lambda^{-1}(h_1+\cdots+h_n)}\otimes\Qlb\bigl(-N(n)\bigr),
\end{equation}
and the geometric and arithmetic monodromy groups are deduced. There, the geometric part of \eqref{eq:det-intro} is obtained from the local monodromy at infinity along the coordinate lines $h_1\mapsto(h_1,c_2,\dots,c_n)$, whose part of slope one has the form $\Lpsi{-h_1/\lambda}\otimes X$ with $X$ of slopes $<1$ but otherwise undetermined, and the Tate twist is obtained by an induction on $n$ through vanishing cycles at the coordinate boundary $h_1=0$.

In this paper, we study $\cW_n$ along \emph{rays}, and we give a ray-theoretic derivation of the determinant formula. The function $-2^{n-1}\lambda^{-1}(h_1+\cdots+h_n)$ in \eqref{eq:det-intro} is the sum of the critical values of the phase of $\Sigma_n$ \cite[Remark~3.7]{EbCoh}, and \cite[Remark~3.7]{EbCoh} emphasizes that the proof given there does not pass through the critical points. Our first result (Theorem~A) shows that the critical points and their values are visible in the local monodromy of $\cW_n$ at the end at infinity of every ray, with their Hessians; the determinant character (Theorem~C) is then obtained by summing the critical values. For the arithmetic normalization (Theorem~D), we establish an integrality and Galois-theoretic constraint and present two sign calculations, one of which uses a perfect pairing between $\cW_n$ and its pull-back by $h\mapsto-h$, which for even $n$ becomes symplectic on the loci of symmetric parameters.

The algebraic critical-value polynomial, its first two coefficients and
its escape discriminant are given in Section~\ref{sec:critical-polynomial}.
Section~\ref{sec:reciprocity} records the resulting Frobenius reciprocity
and identifies the rank-one normalization factor in the
Fourier trace comparison. Dependencies on the companion manuscripts are identified where they enter the arguments.
The consequences developed below include generic-ray large monodromy (conditional on the ambient connected-monodromy theorem), exact characteristic-cycle pull-back for non-characteristic curves, explicit conductors for scalar-profile curves, explicit obstructions to coordinate-wise extension of those formulas, and the exterior-power reciprocity (Theorem~\ref{thm:exterior-reciprocity}), exact determinant distribution (Theorem~\ref{thm:det-distribution}), the complete-ray bound (Theorem~\ref{thm:complete-ray-bound}), and the four-variable escape equation (Theorem~\ref{thm:rank-four-divisor}).
\subsection{Rays and critical points}\label{subsec:rays-intro}
For $a=(a_1,\dots,a_n)\in(\Fpbar^\times)^n$ let
\[
\Delta_a\colon\Gmn{x}\longrightarrow S_n\otimes\Fp(a),\qquad x\longmapsto xa=(xa_1,\dots,xa_n),
\]
be the ray through $a$. Fix $s_i\in\Fpbar$ with $s_i^2=a_i$, put
\[
L_\varepsilon(a)=\sum_{i=1}^n\varepsilon_is_i\quad\bigl(\varepsilon\in\{\pm1\}^n\bigr),\qquad\mathsf E=\{\varepsilon\in\{\pm1\}^n:\varepsilon_1=1\},
\]
\[
P(a)=\#\{\varepsilon\in\mathsf E:L_\varepsilon(a)\ne0\}.
\]
Since $L_{-\varepsilon}=-L_\varepsilon$, the set $\mathsf E$ is a set of representatives of the pairs $\{\pm\varepsilon\}$, and the multiset $\{L_\varepsilon(a)^2:\varepsilon\in\mathsf E\}$ and the number $P(a)$ do not depend on the choice of the $s_i$. The phase of $\Sigma_n$ is
\[
F_n(u,t;h)=\sum_{i=1}^nu_i^{-1}+t\bigl(\lambda+\langle h,u\rangle\bigr)\qquad\text{on }\Gm^n\times\Gmn{t}.
\]
By \cite[Proposition~3.6 and Remark~3.7]{EbCoh}, for $h=a$ its critical points are in bijection with the $\varepsilon\in\mathsf E$ such that $L_\varepsilon(a)\ne0$; the critical point attached to $\varepsilon$ is $u_i=-\lambda\varepsilon_i/(s_iL_\varepsilon(a))$, $t=L_\varepsilon(a)^2/\lambda^2$, it is nondegenerate, with Hessian determinant $2^{n-1}\lambda/(t\,u_1^3\cdots u_n^3)$, and the critical value is $-L_\varepsilon(a)^2/\lambda$. Thus $P(a)$ is the number of critical points; it equals $2^{n-1}$ off the divisor $\prod_\varepsilon L_\varepsilon=0$, along which critical points escape to infinity \cite[Remark~5.4]{EbCoh}.

At the point $x=0$ of a ray, the local monodromy of $\Delta_a^*\cW_n$ is tame and unipotent \cite[Proposition~4.14]{EbMom} (Proposition~\ref{prop:ray-model} below). Our first theorem describes the other end. For a lisse sheaf $\cF$ on a punctured neighbourhood of a point $x$ of a curve over $\Fpbar$ we write $\cF(x)$ for the representation of the inertia group $I_x$ on its geometric generic stalk, and for a character $\chi$ of finite order of the tame inertia group we write $\mathcal L_\chi$ for the corresponding Kummer sheaf; $\chi_2$ is the quadratic character.

\subsection{Statement of the results}\label{subsec:results}

\begin{theoremA}[Stationary phase along rays]
Let $n\ge1$ and $a\in(\Fpbar^\times)^n$, and let $V=(\Delta_a^*\cW_n)(\infty)$.
\begin{enumerate}[label=\textup{(\roman*)}]
\item All breaks of $V$ are $0$ or $1$, and the part of $V$ of break $1$ is
\[
V^{(1)}\cong\bigoplus_{\substack{\varepsilon\in\mathsf E\\ L_\varepsilon(a)\ne0}}\Lpsi{-\frac{L_\varepsilon(a)^2}{\lambda}\,x}\otimes\Lchi{x}^{\otimes(n+1)} .
\]
\item $\Swan_\infty(\Delta_a^*\cW_n)=P(a)$, and the tame part $V^{(0)}$ has dimension $r(n)-P(a)$.
\item If $P(a)=2^{n-1}$, the inertia group acts trivially on $V^{(0)}$.
\item For every $c\in\Fpbar$,
\[
\chi_c\bigl(\Gm\otimes\Fpbar,\Lpsi{cx}\otimes\Delta_a^*\cW_n\bigr)=\begin{cases}-P(a),&c=0,\\-\bigl(r(n)-m_a(c)\bigr),&c\ne0,\end{cases}
\]
where $m_a(c)=\#\{\varepsilon\in\mathsf E:L_\varepsilon(a)^2=\lambda c\}$.
\end{enumerate}
\end{theoremA}

On the ray through $a$ the critical values of the phase at $h=xa$ are $-x\,L_\varepsilon(a)^2/\lambda$, and the Hessian determinant at the corresponding critical point is a nonzero constant times $x^{3n-1}$, whose quadratic character is $\chi_2(x)^{n+1}$. So (i) says that the wild part of the local monodromy at $x=\infty$ is exactly what the principle of stationary phase predicts: one Artin--Schreier sheaf of the critical value for each critical point, twisted by the quadratic character of the Hessian (Remark~\ref{rem:stationary}). When critical points escape to infinity, the corresponding summands disappear from $V^{(1)}$ and are replaced by tame ones; $\cW_n$ itself has no ramification along the critical-point escape divisor \cite[Theorem~5.3]{EbCoh}, and what changes is the behaviour at the end of the ray. Statement (iv) is a global consequence of (i)--(ii); we also prove it directly, by a global argument which does not use local Fourier transforms (Proposition~\ref{prop:global-euler}). It locates the Artin--Schreier exponents $L_\varepsilon(a)^2/\lambda$ as the values of $c$ at which the Euler characteristic of the twist jumps, and it is compatible with the finite algebraic checks recorded in Section~\ref{sec:numerics}.

The key local input is the following explicit local Fourier transform. Let $I_\infty$ be the inertia group at $\infty$ of the $y$-line, $J\subset I_\infty$ its open subgroup of index two, the inertia group of the double cover $w\mapsto w^2=y$, and for $c\in\Fpbar$ and $e\in\{0,1\}$ put
\[
\Msf_e(c)=\Ind_J^{I_\infty}\bigl(\Lpsi{cw}\otimes\Lchi{w}^{\otimes e}\bigr),
\]
a representation of dimension two, irreducible of slope $\frac12$ if $c\ne0$.

\begin{theoremB}[An explicit local Fourier transform]
For $c\in\Fpbar^\times$ and $e\in\{0,1\}$, Laumon's local Fourier transform satisfies
\[
\Floc\bigl(\Msf_e(c)\bigr)\cong\Lpsi{-\frac{c^2}{4z}}\otimes\Lchi{z}^{\otimes(1-e)} .
\]
\end{theoremB}

Theorem~B is the case $r=2$, $s=1$ of the explicit formulas for local Fourier transforms of Fu \cite[Theorem~0.1(iii)]{Fu} (see also Abbes--Saito \cite{AS} and, for $\mathcal D$-modules, Sabbah \cite{Sabbah}). We give an independent proof: we show that $\Floc(\Msf_e(c))$ governs the end of the rays of the family for every $n$, and we \emph{calibrate} it with the cases $n=1$ and $n=2$, where $\cW_1$ and $\cW_2$ are explicit; this recovers the constant $-\frac14$ and the quadratic twist without any computation of vanishing cycles (Section~\ref{subsec:calibration}).

\begin{theoremC}[The geometric determinant]
For every $n\ge1$ the lisse sheaf $\det\cW_n\otimes\Lpsi{2^{n-1}\lambda^{-1}(h_1+\cdots+h_n)}$ on $S_n$ is geometrically constant. Equivalently, there is $\alpha\in\Qlb^\times$ such that
\[
\det\cW_n\cong\Lpsi{-2^{n-1}\lambda^{-1}(h_1+\cdots+h_n)}\otimes\alpha^{\deg},
\]
where $\alpha^{\deg}$ is the geometrically constant sheaf on which $\Frob_p$ acts by $\alpha$. In particular the image of $\pgeom(S_n)$ under $\det\cW_n$ is $\mu_p$.
\end{theoremC}

The statement of Theorem~C is \cite[Theorem~8.11]{EbCoh}. Our proof is different: Theorem~A computes the determinant along every ray as the Artin--Schreier sheaf of the sum of the critical values, and a lemma on characters of $\pi_1(\Gm^n\otimes\Fpbar)$ (Lemma~\ref{lem:torus}) propagates this from rays to $S_n$, the only other input being the unipotence of the local monodromy at the coordinate hyperplanes \cite[Theorem~6.8(ii)]{EbCoh}.

\begin{theoremD}[The arithmetic normalization]
Let $n\ge1$ and let $\alpha$ be as in Theorem~C. Then $\alpha=p^{N(n)}$, that is,
\[
\det\cW_n\cong\Lpsi{-2^{n-1}\lambda^{-1}(h_1+\cdots+h_n)}\otimes\Qlb\bigl(-N(n)\bigr).
\]
More precisely:
\begin{enumerate}[label=\textup{(\roman*)}]
\item \textup{(Product formula and Galois symmetry.)} $\alpha\in\Z[\zeta_p]$, the number $\alpha\,p^{-N(n)}$ is a root of unity, and $\alpha\in\Q$; hence $\alpha=\pm p^{N(n)}$.
\item \textup{(Reflection.)} There is a perfect pairing of lisse sheaves on $S_n$ over $\Fp$
\[
\mathcal R\colon\cW_n\otimes[-1]^*\cW_n\longrightarrow\Qlb(-n-1),\qquad[-1](h)=-h,
\]
such that $[-1]^*\mathcal R=(-1)^{n+1}\,\mathcal R\circ\sw$. Consequently $[-1]^*\cW_n\cong\cW_n^\vee(-n-1)$ and $\alpha^2=p^{2N(n)}$.
\item \textup{(Symmetric parameters.)} Let $n$ be even and let $\sigma$ be a fixed-point-free involution of $\{1,\dots,n\}$. The restriction of $\cW_n$ to $\Sigma_\sigma=\{h\in S_n:h_{\sigma(i)}=-h_i\ \text{for all }i\}\cong\Gm^{n/2}$ carries a perfect \emph{alternating} pairing with values in $\Qlb(-n-1)$. Consequently $\det(\Frob_h\mid\cW_n)=q^{N(n)}$ for every $h\in\Sigma_\sigma(\Fq)$, and $\alpha=p^{N(n)}$.
\item \textup{(The tame end of a ray.)} For every $a\in(\Fq^\times)^n$, the determinant of a Frobenius element of the decomposition group at $x=0$ acting on $\Delta_a^*\cW_n$ is $q^{N(n)}$. Consequently $\alpha=p^{N(n)}$.
\end{enumerate}
\end{theoremD}

The value $\alpha=p^{N(n)}$ is \cite[Theorem~8.17]{EbCoh}, where it is obtained by induction on $n$ through vanishing cycles at the coordinate boundary. The product formula alone proves only $\alpha=\pm p^{N(n)}$; the symplectic and tame-end arguments provide two separate sign determinations and do not use that induction. Theorem~C enters them only through the definition of $\alpha$: the first assertions of (iii) and (iv), on $\det(\Frob_h\mid\cW_n)$ along $\Sigma_\sigma$ and on the Frobenius at the tame end of a ray, do not use it. Statement (iv) is a reformulation of \cite[Proposition~4.18(iii)]{EbMom}, whose proof we give in the form needed (Proposition~\ref{prop:tame-end}). Statements (ii) and (iii) are new. For $p=2$ an invariant bilinear form on $\cW_n$ exists globally \cite[Proposition~4.7]{EbMom}; for odd $p$ there is none \cite[Theorem~4.10]{EbMom}, and (ii)--(iii) show what survives: a duality with the reflected sheaf, and a symplectic structure on the symmetric loci. Combining Theorems~C and~D with the identity component $\SL_{r(n)}$ gives the monodromy groups.

\begin{corollaryE}[{\cite[Theorem~E]{EbCoh}}]
Let $n\ge2$ and $p>\max(2,n)$, so that the connected-monodromy theorem quoted above applies. Then $\ggeom(\cW_n)=\{g\in\GL_{r(n)}(\Qlb):\det(g)^p=1\}$, and the arithmetic monodromy group of $\cW_n\otimes\gamma^{\deg}$, where $\gamma\in\Qlb$ satisfies $\gamma^2=p^{-(n+1)}$, and $\gamma=p^{-(n+1)/2}$ if $n$ is odd, is equal to $\ggeom(\cW_n)$. Moreover, for every $h\in(\Fq^\times)^n$,
\[
\det\Bigl(\Frob_q\mid H^{n-1}_c\bigl(Y_h\otimes\Fqbar,\Lpsi{\textstyle\sum_iu_i^{-1}}\bigr)\Bigr)=q^{(n-2)2^{n-1}+1}\,\psi_q\bigl(-2^{n-1}\lambda^{-1}(h_1+\cdots+h_n)\bigr).
\]
\end{corollaryE}

\subsection{Strategy}\label{subsec:strategy}
\emph{Rays and Fourier transforms.} Along the ray through $a$, a change of variables turns the family into the Fourier transform, with respect to the pairing $(y,z)\mapsto yz$, of the tensor product $\cK_a=\bigotimes_i[\times a_i]^*\Kl$ of dilated Kloosterman sheaves, followed by the inversion $z=\lambda/x$ \cite[Proposition~4.14]{EbMom}. The end $x=\infty$ of the ray corresponds to the finite point $z=0$ of the Fourier transform. There the vanishing cycles of the Fourier transform are given by Laumon's local Fourier transform $\Floc$ of the local monodromy of $\cK_a$ at $y=\infty$ \cite[Proposition~(2.3.2.1)]{Laumon}, and the nearby cycles differ from them only by modules with trivial inertia action (Lemma~\ref{lem:vanishing}). The local monodromy of $\Kl$ at $\infty$ is the induced representation $\Msf_1(2)$ \cite[Theorem~1.1]{FuWan}, \cite[Chapter~10]{Katz-GKM}, and Mackey's formula decomposes $\cK_a(\infty)$ as $\bigoplus_{\varepsilon\in\mathsf E}\Msf_n(2L_\varepsilon(a))$, the index $e$ of $\Msf_e$ being taken modulo $2$. Hence the wild part at the end of the ray is $\bigoplus_{L_\varepsilon\ne0}\Floc(\Msf_n(2L_\varepsilon(a)))$, transported by $z=\lambda/x$, and Laumon's rank and Swan formulas show that each summand is a character of Swan conductor one (Proposition~\ref{prop:ray-decomp}).

\emph{Calibration.} For $n=1$ and $n=2$ the sheaves $\cW_1$ and $\cW_2$ are explicit \cite[Lemma~4.17, Proposition~2.21]{EbMom}: $\cW_1$ is an Artin--Schreier sheaf, and $\cW_2$ is an Artin--Schreier twist of the pull-back of $\Kl$ by $(h_1,h_2)\mapsto h_1h_2/\lambda^2$. Comparing the two descriptions of the end of a ray for $n=1$ gives Theorem~B for $e=1$. For $n=2$ the comparison gives a multiset identity between two pairs of characters; fixing $s_1+s_2$ and letting $s_1-s_2$ vary forces Theorem~B for $e=0$. Theorem~A follows for all $n$.

\emph{Determinant.} By Theorem~A and the identity $\sum_{\varepsilon\in\mathsf E}L_\varepsilon(a)^2=2^{n-1}(a_1+\cdots+a_n)$, the determinant of $\Delta_a^*\cW_n$ is $\Lpsi{-2^{n-1}\lambda^{-1}(\sum_ia_i)x}$ times a character which is tame at $\infty$ and unramified at $0$, hence geometrically trivial. A character of finite order of $\pi_1(\Gm^n\otimes\Fpbar)$ which is trivial on all rays and unramified at the origin of every coordinate line is trivial (Lemma~\ref{lem:torus}); its proof combines Kummer theory with a normal form for Artin--Schreier classes of Laurent polynomials. This gives Theorem~C.

\emph{Arithmetic normalization.} The number $\alpha$ is the determinant of Frobenius on $\cW_{n,h}$ at any $h\in S_n(\Fp)$, divided by a $p$-th root of unity. It is an algebraic integer of $\Q(\zeta_p)$ all of whose conjugates have absolute value $p^{N(n)}$; as $p$ is totally ramified in $\Q(\zeta_p)$, $\alpha p^{-N(n)}$ is a unit, hence a root of unity by Kronecker's theorem. Replacing $\psi$ by $\psi^k$ is the same as dilating $h$ by $k$, which leaves $\alpha$ unchanged; so $\alpha$ is fixed by $\Gal(\Q(\zeta_p)/\Q)$, and $\alpha=\pm p^{N(n)}$. To fix the sign we construct, by a cup product on the compactly supported cohomology along the $t$-line, a pairing between the fibres of $\cW_n$ at $h$ and at $-h$: the involution $t\mapsto-t$ carries $\Lpsi{\lambda t}$ to its inverse, while the products $th_i$ are unchanged when both $t$ and $h$ change sign, so the symplectic self-duality of $\Kl$ provides the pairing. Its radical is the lower-weight part, and the Koszul sign rule gives the symmetry $(-1)^{n+1}$. At a symmetric parameter $h$ the fibres at $h$ and $-h$ are identified by a permutation of the coordinates, and for even $n$ the resulting form on $\cW_{n,h}$ is alternating and Frobenius-equivariant with multiplier $q^{n+1}$; the determinant of a symplectic similitude is a power of its multiplier, which gives $\alpha=p^{N(n)}$ for even $n$. For all $n$, the Frobenius at the tame end of a ray is computed through Laumon's local Fourier transform $\Flocz$ over a finite field, the Frobenius on $\Kl$ at $0$ being $\mathrm{diag}(1,q)$.

\subsection{Relation to other work}\label{subsec:relation}
This paper grew out of the author's longer manuscript \cite{EbFull} and is a companion to the preprints \cite{EbMom} and \cite{EbCoh}, in their present preprint form. Numbered references to \cite{EbCoh} are to its third version, arXiv:2610.03081v3, and numbered references to \cite{EbMom} are to its version of October 2026. The relation of the present paper to these texts is as follows.
\begin{enumerate}[label=(\alph*)]
\item \emph{Results used.} We use the construction, rank, weight filtration, purity, and geometric irreducibility of $\cW_n$ from the companion preprints \cite{EbMom,EbCoh}, as recorded in Theorem~\ref{thm:recollection}. We also use the explicit rank-one and rank-two models, the Fourier presentation along rays, the tame-end calculation, and the Frobenius action on $\Kl$ at the origin from \cite{EbMom}. Three further inputs recur in distinct parts of this paper:
\begin{enumerate}[label=\textup{(\roman*)}]
\item The connected-monodromy identity $G_{\mathrm{geom}}(\cW_n)^0=\SL_{r(n)}$ from \cite[Theorem~8.8]{EbCoh} (compare \cite[Theorem~4.19]{EbMom}) is used in Corollary~E, Theorem~\ref{thm:generic-ray-SL}, and Proposition~\ref{prop:determinant-completion}. Its characteristic range is specified separately in those statements.
\item Coordinate-boundary tameness and unipotence are imported from \cite[Theorem~6.8(ii)]{EbCoh} in the proof of Theorem~C, Theorem~\ref{thm:coordinate-jordan}, and Proposition~\ref{prop:cc-zero-coord}. The cited theorem states the \emph{exact} Jordan type: $r(n-1)$ blocks of size two, $\kappa_1$ blocks of size one, and no larger blocks, where $\kappa_1=r(n)-2r(n-1)$. The cited Jordan-profile theorem is proved in Section~6.3 of the separate cohomology preprint \cite{EbCoh}, using local Fourier transforms and the unipotent middle-convolution formalism; its geometric proof is not repeated here. The binomial and Catalan consequences of the imported profile are proved in the present article.
\item The fourth-moment asymptotic for $n\ge3$, in the odd-characteristic case, is imported from \cite[Theorem~3.13]{EbMom}; it appears in Theorem~\ref{thm:imported-moment} and is used only through the cited ambient-monodromy argument. No independent fourth-moment calculation is claimed here.
\end{enumerate}
Every use of these inputs is indicated at the relevant argument. The critical-point description is likewise invoked only where explicitly cited.
\item \emph{Results obtained again, by different methods.} Theorem~C is \cite[Theorem~8.11]{EbCoh}, and the value of $\alpha$ in Theorem~D is \cite[Theorem~8.17]{EbCoh}; Corollary~E is \cite[Theorem~E]{EbCoh}. Theorem~D(iv) is \cite[Proposition~4.18(iii)]{EbMom}. We do not claim these statements as new; what is new is the route through rays and the arguments of Theorem~D(i)--(iii).
\item \emph{Results developed here.} Theorem~A; the calibration argument for Theorem~B (whose statement is a case of \cite[Theorem~0.1(iii)]{Fu}); Proposition~\ref{prop:global-euler}; the reflection pairing, its symplectic specialization, and the arithmetic argument of Theorem~D(i)--(iii); the exterior-power and determinant-distribution consequences; the generic-ray result, using the ambient monodromy theorem; and the curve-conductor and parameter-space Euler-characteristic formulas. The higher-codimension characteristic-cycle recursion is stated with its microlocal hypotheses.
\end{enumerate}
The 169-page parent manuscript \cite{EbFull} obtains a root of unity and then a rational sign in Theorems~13.64 and~13.66, but leaves its positivity as Conjecture~13.68. The present Theorem~D(iii)--(iv) proposes two distinct ways to determine that sign, subject to the explicitly stated cohomological and local-Fourier inputs; neither the parent conjecture nor its numerical evidence is treated as a proof. The references above identify the actual numbering of the 169-page source.

The principle that the local monodromy of a Fourier transform at a finite point is governed by local Fourier transforms is due to Laumon \cite{Laumon}; see also Katz's exposition \cite{Katz-Bourbaki}. Explicit formulas for the local Fourier transforms of representations of the form $[r]_*(\mathcal L_\psi(f)\otimes\mathcal L_\chi)$ were obtained by Fu \cite{Fu} and by Abbes and Saito \cite{AS}, and by Sabbah \cite{Sabbah} for formal connections. The local monodromy of Kloosterman sheaves at $\infty$ is due to Katz \cite[Chapter~10]{Katz-GKM} and Fu--Wan \cite[Theorem~1.1]{FuWan}. Sums of products of Kloosterman sums such as \eqref{eq:t-sum} are studied in general by Fouvry, Kowalski and Michel \cite{FKM}. The stationary phase principle relates critical points of an exponential phase to local monodromy. Theorem~A gives a precise formulation for the restrictions studied here.

\subsection{Logical structure}\label{subsec:logic}
Besides the general theory of $\ell$-adic cohomology \cite{SGA4,SGA45,SGA5,Deligne} we use: Laumon's theory of the Fourier transform and of local Fourier transforms \cite{Laumon}, in the form recalled in Section~\ref{sec:local}; the facts on $\Kl$ recalled in Section~\ref{subsec:Kl}, including \cite[Theorem~1.1]{FuWan}; the results of \cite{EbMom} listed in (a) above; Deligne's integrality theorem \cite[Exp.~XXI, Corollaire~5.5.3]{SGA7}, Kronecker's theorem \cite{Kronecker} and the Chebotarev density theorem \cite{Serre-Zeta}. Theorems~A and~B use neither \cite{EbCoh} nor any of Theorems~C--E. Theorem~C uses Theorem~A and \cite[Theorem~6.8(ii)]{EbCoh}. The pairings of Theorem~D(ii),(iii) are constructed independently of Theorems~A--C, and the computation of Theorem~D(iv) is independent of Theorems~A and~B; the conclusions about $\alpha$ in Theorem~D(i)--(iv) use Theorem~C, which defines $\alpha$. Corollary~E, Theorem~\ref{thm:generic-ray-SL}, and Proposition~\ref{prop:determinant-completion} each invoke the connected ambient-monodromy theorem from the companion preprints. The coordinate-boundary Jordan theorem and its conormal consequence are applications of \cite[Theorem~6.8(ii)]{EbCoh}, while the fourth moment enters only as Theorem~\ref{thm:imported-moment}. The numerical checks of Section~\ref{sec:numerics} are not used in any proof.

\subsection{Conventions}\label{subsec:conventions}
For a scheme $X$ of finite type over $\Fq$ and $x\in X(\Fq)$, $\Tr(\Frob_x\mid\cF)$ is the trace of the geometric Frobenius on the stalk at $x$. Weights are taken with respect to every field isomorphism $\iota\colon\Qlb\to\C$; all sheaves considered are $\iota$-mixed for every $\iota$. For $a\in\Qlb^\times$, $a^{\deg}$ denotes the geometrically constant lisse sheaf of rank one on which the geometric Frobenius of $\Fp$ acts by $a$, or its pull-back to any scheme over $\Fp$; thus $\Qlb(-1)=p^{\deg}$. For $c\in\Fpbar^\times$, $[\times c]$ is the automorphism $x\mapsto cx$ of $\Gm$, and $\Lpsi{f}$ is the Artin--Schreier sheaf of a function $f$, with trace function $\psi_q(f)$. For a point $x$ of a smooth curve over $\Fpbar$, $I_x$ is the inertia group, $P_x\subset I_x$ the wild inertia group, $I_x^{(r)}$ ($r\ge0$) the upper numbering filtration, and breaks and Swan conductors are as in \cite[Chapter~1]{Katz-GKM}. A unipotent Jordan block of size $m$ is denoted $\Up m$. For a pairing $B\colon\cF\otimes\cG\to\cL$ we write $B\circ\sw$ for the pairing $\cG\otimes\cF\to\cL$, $g\otimes f\mapsto B(f\otimes g)$.
\section{Preliminaries}\label{sec:prelim}

\subsection{Representations of inertia groups}\label{subsec:inertia}
Let $x$ be a closed point of a smooth curve over $\Fpbar$, and let $I=I_x$, $P=P_x$. Every finite-dimensional continuous $\Qlb$-representation $V$ of $I$ has a unique \emph{break decomposition} $V=\bigoplus_{r\ge0}V^{(r)}$ into $I$-stable subspaces, where $V^{(0)}=V^{P}$ is the tame part and $V^{(r)}$ is isoclinic of break $r$; it is functorial, hence exact, in $V$, and $\Swan(V)=\sum_rr\dim V^{(r)}$ is an integer \cite[Chapter~1]{Katz-GKM}. A representation all of whose breaks are $\le r$ is one on which $I^{(r')}$ acts trivially for all $r'>r$ \cite[1.1]{Katz-GKM}.

\begin{lemma}[{\cite[Lemma~2.2]{EbMom}, \cite[Lemma~2.3]{EbCoh}}]\label{lem:breaks}
Let $V,V'$ be finite-dimensional continuous representations of $I$.
\begin{enumerate}[label=\textup{(\roman*)}]
\item Every break of $V\otimes V'$ is at most the largest break of $V$ and of $V'$.
\item If all breaks of $V$ are $<b$ and $V'$ is isoclinic of break $b>0$, then $V\otimes V'$ is isoclinic of break $b$, and $\Swan(V\otimes V')=b\dim V\dim V'$.
\end{enumerate}
\end{lemma}

The group $P$ acts on $V$ through a finite quotient, and $I/P$ is procyclic, canonically isomorphic to $\varprojlim_{(m,p)=1}\mu_m(\Fpbar)$. So the tame characters of $I$ of finite order are the restrictions of Kummer sheaves $\mathcal L_\chi$, and $\mathcal L_\chi(cx')\cong\mathcal L_\chi(x')$ as representations of $I$ for every uniformizing parameter $x'$ and every $c\in\Fpbar^\times$. Since $p$ is odd, $I/P$ has a unique quotient of order two; let $J\subset I$ be the corresponding open subgroup of index two. It is the inertia group at the point above $x$ of the double cover $w\mapsto w^2=x'$, for any uniformizer $x'$ at $x$. If $V$ is a representation of $J$, the breaks of $\Ind_J^IV$ are the breaks of $V$ divided by $2$ \cite[1.13]{Katz-GKM}, \cite[Chapter~IV, \S3]{Serre-LF}. Let $\gamma\in I\smallsetminus J$. For a character $\xi$ of $J$ put $\xi^\gamma(g)=\xi(\gamma^{-1}g\gamma)$; if $\xi$ is the restriction of a lisse sheaf $\cG$ on a punctured neighbourhood of the point above $x$ on the double cover, then $\xi^\gamma$ is the restriction of the pull-back of $\cG$ by the deck transformation $w\mapsto-w$.

We apply this at the point $y=\infty$ of the projective line with coordinate $y$, with the uniformizer $1/y$ and the double cover $w^2=y$. For $c\in\Fpbar$ and $e\in\Z$ let $\xi_{c,e}=\Lpsi{cw}\otimes\Lchi{w}^{\otimes e}$, regarded as a character of $J$, and
\[
\Msf_e(c)=\Ind_J^{I_\infty}\xi_{c,e},
\]
which depends only on $e$ modulo $2$.

\begin{lemma}\label{lem:M}
Let $c,c'\in\Fpbar$ and $e,e'\in\Z$.
\begin{enumerate}[label=\textup{(\alph*)}]
\item If $c\ne0$, then $\Msf_e(c)$ is irreducible, both of its breaks equal $\frac12$, and $\Swan(\Msf_e(c))=1$.
\item $\Msf_e(c)\cong\Msf_e(-c)$.
\item $\Msf_e(0)$ is tame; it is the direct sum of the two characters of $I_\infty/P_\infty$ whose restriction to $J$ is $\Lchi{w}^{\otimes e}$.
\item $\Msf_e(c)\otimes\Msf_{e'}(c')\cong\Msf_{e+e'}(c+c')\oplus\Msf_{e+e'}(c-c')$.
\item Let $b\in\Fpbar^\times$ and $\beta^2=b$. If $\cG$ is a lisse sheaf on a punctured neighbourhood of $y=\infty$ with $\cG(\infty)\cong\Msf_e(c)$, then $([\times b]^*\cG)(\infty)\cong\Msf_e(\beta c)$.
\item $\det\Msf_e(c)\cong\Lchi{y}^{\otimes(e+1)}$.
\end{enumerate}
\end{lemma}

\begin{proof}
(a) By the description of conjugation recalled above, $\xi_{c,e}^\gamma\cong\Lpsi{-cw}\otimes\Lchi{-w}^{\otimes e}\cong\xi_{-c,e}$. The character $\xi_{c,e}^\gamma\xi_{c,e}^{-1}\cong\Lpsi{-2cw}$ has Swan conductor $1$ at $w=\infty$, because $2c\ne0$; so $\xi^\gamma_{c,e}\ne\xi_{c,e}$, and $\Msf_e(c)$ is irreducible by Mackey's criterion. Its restriction to $J$ is $\xi_{c,e}\oplus\xi_{c,e}^\gamma$, which is isoclinic of break $1$. Hence both breaks of $\Msf_e(c)$ equal $\frac12$, and its Swan conductor is $1$.

(b) An induced representation does not change when the inducing character is replaced by a conjugate: $\Ind\xi\cong\Ind\xi^\gamma$, and $\xi_{c,e}^\gamma\cong\xi_{-c,e}$.

(c) The group $P_\infty\subset J$ acts trivially, and $I_\infty/P_\infty$ is procyclic; the representation of a procyclic group induced from a character of its subgroup of index two is the sum of the two extensions of that character.

(d) By the projection formula, $\Ind\xi\otimes\Ind\xi'\cong\Ind(\xi\otimes\Res_J\Ind\xi')=\Ind(\xi\otimes\xi')\oplus\Ind(\xi\otimes\xi'^\gamma)$, and $\xi_{c,e}\otimes\xi_{c',e'}\cong\xi_{c+c',e+e'}$, $\xi_{c,e}\otimes\xi_{c',e'}^\gamma\cong\xi_{c-c',e+e'}$, as $\Lchi{w}^{\otimes2}$ is trivial.

(e) The automorphism $y\mapsto by$ fixes $\infty$ and lifts to the automorphism $w\mapsto\beta w$ of the double cover. The automorphism of $I_\infty$ that it induces preserves $J$, and it carries $\Ind_J^{I_\infty}\xi$ to $\Ind_J^{I_\infty}$ of the pull-back of $\xi$; the pull-back of $\xi_{c,e}$ is $\Lpsi{c\beta w}\otimes\Lchi{\beta w}^{\otimes e}\cong\xi_{\beta c,e}$.

(f) The determinant of $\Ind_J^{I_\infty}\xi$ is the product of the sign character of $I_\infty/J$, which is $\Lchi{y}$, and of the transfer of $\xi$, which corresponds to the norm of the rank-one sheaf from the double cover to the $y$-line. The norm of $\Lpsi{cw}$ is the Artin--Schreier sheaf of the trace $cw+c(-w)=0$, hence trivial, and the norm of $\Lchi{w}$ is $\Lchi{w\cdot(-w)}=\Lchi{-y}\cong\Lchi{y}$.
\end{proof}

\subsection{The Kloosterman sheaf}\label{subsec:Kl}
Let $\pr\colon\Gmn{v}\times\Gmn{y}\to\Gmn{y}$ be the projection and $\Kl=R^1\pr_!\Lpsi{v^{-1}+yv}$, a sheaf on $\Gmn{y}$ over $\Fp$. Under the isomorphism $\Gm\to\{(v_1,v_2):v_1v_2=y\}$, $v\mapsto(v^{-1},yv)$, it is Deligne's Kloosterman sheaf $\mathrm{Kl}_2$ \cite[Sommes trig., \S7]{SGA45}. We use the following facts.
\begin{enumerate}[label=(K\arabic*)]
\item $R^i\pr_!\Lpsi{v^{-1}+yv}=0$ for $i\ne1$; $\Kl$ is lisse of rank $2$ and pure of weight $1$, and $\Tr(\Frob_y\mid\Kl)=-K(y)$ for $y\in\Fq^\times$ \cite[Theorem~2.1]{EbMom}.
\item $\Kl$ is tame at $0$, with unipotent local monodromy $\Up2$; it is totally wild at $\infty$, with both breaks equal to $\frac12$ and $\Swan_\infty(\Kl)=1$ \cite[Theorem~2.1]{EbMom}, \cite[Chapters~4 and~7]{Katz-GKM}.
\item Let $D_0$ be a decomposition group at the point $y=0$ of $\A^1_y$ over $\Fq$, $\Phi\in D_0$ an element mapping to the geometric Frobenius of $\Fq$, and $\mathfrak N$ the logarithm of the tame monodromy on $\Kl(0)$. There is a basis $e_0,e_1$ of the generic stalk of $\Kl$ with $\mathfrak Ne_1=e_0$, $\mathfrak Ne_0=0$, $\Phi e_0=e_0$ and $\Phi e_1=qe_1$; the same holds for $[\times c]^*\Kl$, $c\in\Fq^\times$ \cite[Lemma~2.14]{EbMom}, \cite[(K5) and Lemma~2.5]{EbCoh}.
\item $\det\Kl\cong\Qlb(-1)$ as lisse sheaves on $\Gm$ over $\Fp$ \cite[Lemma~4.6]{EbMom}.
\end{enumerate}

\begin{proposition}[{Katz \cite[Chapter~10]{Katz-GKM}, Fu--Wan \cite[Theorem~1.1]{FuWan}}]\label{prop:Kl-infty}
As a representation of $I_\infty$, $\Kl(\infty)\cong\Msf_1(2)=\Ind_J^{I_\infty}\bigl(\Lpsi{2w}\otimes\Lchi{w}\bigr)$.
\end{proposition}

\begin{proof}
For $n=2$, \cite[Theorem~1.1]{FuWan} gives, as representations of the decomposition group at $\infty$, $\mathrm{Kl}_2\cong[2]_*(\mathcal L_{\psi_2}\otimes\mathcal L_\chi)\otimes\mathcal L_\theta\otimes\Qlb(-\frac12)$, where $[2]$ is the squaring map, $\psi_2(w)=\psi(2w)$, $\chi$ is the quadratic character and $\theta$ is a character of the Galois group of the finite field. Restricted to $I_\infty$, the last two factors are trivial and $[2]_*$ is $\Ind_J^{I_\infty}$.
\end{proof}

Proposition~\ref{prop:Kl-infty} is consistent with (K2) and (K4) through Lemma~\ref{lem:M}(a),(f): $\Msf_1(2)$ has both breaks $\frac12$ and trivial determinant.

\begin{corollary}\label{cor:Ka}
Let $n\ge1$, $a\in(\Fpbar^\times)^n$, and $\cK_a=\bigotimes_{i=1}^n[\times a_i]^*\Kl$ on $\Gmn{y}$. Then $\cK_a$ is lisse of rank $2^n$ and pure of weight $n$; it is tame at $0$, with unipotent local monodromy of the Jordan type of $\Up2^{\otimes n}$; and
\[
\cK_a(\infty)\cong\bigoplus_{\varepsilon\in\mathsf E}\Msf_n\bigl(2L_\varepsilon(a)\bigr).
\]
In particular all breaks of $\cK_a(\infty)$ are $0$ or $\frac12$, the tame part is $\bigoplus_{L_\varepsilon(a)=0}\Msf_n(0)$, and $\Swan_\infty(\cK_a)=P(a)$.
\end{corollary}

\begin{proof}
The first assertions follow from (K1), (K2) and the fact that $[\times a_i]$ fixes $0$ and $\infty$ and induces the identity on the tame quotient of $I_0$. By Proposition~\ref{prop:Kl-infty} and Lemma~\ref{lem:M}(e), $([\times a_i]^*\Kl)(\infty)\cong\Msf_1(2s_i)$. We show by induction on $m\le n$ that $\bigotimes_{i\le m}\Msf_1(2s_i)\cong\bigoplus\Msf_m(2\sum_{i\le m}\varepsilon_is_i)$, the sum over $\varepsilon\in\{\pm1\}^m$ with $\varepsilon_1=1$. For $m=1$ this is clear, and the induction step is Lemma~\ref{lem:M}(d). The last assertions follow from Lemma~\ref{lem:M}(a),(c).
\end{proof}

\subsection{The sheaves \texorpdfstring{$\cH_n$}{H}, \texorpdfstring{$\cE_n$}{E} and \texorpdfstring{$\cW_n$}{W}}\label{subsec:sheaves}
For $n\ge1$ let $m_i\colon\Gmn{t}\times S_n\to\Gm$ be the map $m_i(t,h)=h_it$, put
\[
\cN_n=\Bigl(\bigotimes_{i=1}^nm_i^*\Kl\Bigr)\otimes\Lpsi{\lambda t}\qquad\text{on }\Gmn{t}\times S_n,
\]
a lisse sheaf of rank $2^n$, pure of weight $n$, and let $\rho\colon\Gmn{t}\times S_n\to S_n$ be the projection. For a geometric point $h$ of $S_n$ we write $\cN_h$ for the restriction of $\cN_n$ to $\Gmn{t}\times\{h\}$; let $j\colon\Gmn{t}\hookrightarrow\PP^1_t$ be the inclusion.

\begin{theorem}[{\cite{EbMom,EbCoh}}]\label{thm:recollection}
Let $n\ge1$.
\begin{enumerate}[label=\textup{(\roman*)}]
\item $R^i\rho_!\cN_n=0$ for $i\ne1$, and $\cH_n:=R^1\rho_!\cN_n$ is a lisse sheaf of rank $2^n$ on $S_n$, mixed of weights $\le n+1$, whose formation commutes with base change. For $h\in S_n(\Fq)$ there is an isomorphism, compatible with $\Frob_h$,
\[
\cH_{n,h}=H^1_c(\Gm\otimes\Fqbar,\cN_h)\cong H^{n+1}_c\bigl((\Gm^n\times\Gmn{t})\otimes\Fqbar,\Lpsi{F_n(\cdot,\cdot;h)}\bigr).
\]
\item There is an exact sequence of lisse sheaves $0\to\cE_n\to\cH_n\to\cW_n\to0$, where $\cE_n$ is the part of weights $\le n$ of the weight filtration of $\cH_n$; $\cE_n$ has rank $b_n$, and $\cW_n$ is pure of weight $n+1$ and of rank $r(n)$. At a geometric point $h$, $\cE_{n,h}$ is the image of $\cN_h^{I_0}$ and $\cW_{n,h}=H^1(\PP^1,j_*\cN_h)$.
\item The trace formula \eqref{eq:trace} holds.
\item For $0\le k\le\lfloor n/2\rfloor$ the graded piece $\gr^W_{2k}\cE_n$ is geometrically constant, of rank $\nu_{n,k}$, with all Frobenius eigenvalues at points of $S_n(\Fq)$ equal to $q^k$, and the graded pieces of odd weight vanish. Consequently, for every smooth curve $C$ over $\Fpbar$ and every morphism $f\colon C\to S_n$, the sheaf $f^*\cE_n$ is tame, with unipotent local monodromy, at every point of the smooth compactification of $C$.
\item $\cW_1\cong\Lpsi{-h/\lambda}\otimes\Qlb(-1)$, and $\cW_2\cong\Lpsi{-(h_1+h_2)/\lambda}\otimes\mathrm m^*\Kl\otimes\Qlb(-1)$, where $\mathrm m(h)=h_1h_2/\lambda^2$.
\item The image of $\pgeom(S_n)$ under $\det\cW_n$ is finite.
\item If $n\ge2$, then $\cW_n$ is geometrically irreducible.
\end{enumerate}
\end{theorem}

\begin{proof}
(i) is \cite[Theorem~2.10(i)]{EbMom} and \cite[Proposition~4.1 and Theorem~5.3(ii)]{EbCoh}; (ii) is \cite[Theorem~2.10(ii)]{EbMom}; (iii) is \cite[Corollary~2.16(i)]{EbMom}; (iv) is \cite[Corollary~2.18]{EbMom}, see also \cite[Theorem~5.3(vi)]{EbCoh}; (vi) is \cite[Proposition~4.5]{EbMom}, an instance of \cite[Proposition~1.3.4(i)]{Deligne}; (vii) is \cite[Corollary~4.20]{EbMom}. For (v), which is \cite[Lemma~4.17 and Proposition~2.21]{EbMom}, we recall the argument, as it is the basis of the calibration in Section~\ref{subsec:calibration}. For $n=1$ the only $u$ with $hu=-\lambda$ is $u=-\lambda/h$, so \eqref{eq:trace} gives $\Tr(\Frob_h\mid\cW_1)=q\psi_q(-h/\lambda)$. For $n=2$ and $h\in(\Fq^\times)^2$, the substitution $v_i=h_iu_i$ and the parametrization $(v_1,v_2)=(-\lambda/(1+\eta),-\lambda\eta/(1+\eta))$, $\eta\ne0,-1$, of the line $v_1+v_2=-\lambda$ give
\[
\Sigma_2(h)=\sum_{\eta\ne0,-1}\psi_q\Bigl(-\frac{h_1+h_2}\lambda-\frac{h_1}\lambda\eta-\frac{h_2}{\lambda\eta}\Bigr)=\psi_q\Bigl(-\frac{h_1+h_2}\lambda\Bigr)K\Bigl(\frac{h_1h_2}{\lambda^2}\Bigr)-1,
\]
the term $\eta=-1$ of the Kloosterman sum being $\psi_q((h_1+h_2)/\lambda)$. Hence, by \eqref{eq:trace} with $\nu_{2,1}=1$ and (K1), $\cW_1$ and $\Lpsi{-h/\lambda}(-1)$, respectively $\cW_2$ and $\Lpsi{-(h_1+h_2)/\lambda}\otimes\mathrm m^*\Kl(-1)$, have the same Frobenius traces at all closed points. By the Chebotarev density theorem \cite{Serre-Zeta} and the Brauer--Nesbitt theorem their semisimplifications are isomorphic. The right-hand sides are irreducible: for $n=1$ because they have rank one, and for $n=2$ because $\Kl$ is irreducible already as a representation of $I_\infty$ (Lemma~\ref{lem:M}(a) and Proposition~\ref{prop:Kl-infty}) and $\mathrm m$ has the section $x\mapsto(\lambda^2x,1)$, so that $\mathrm m^*\Kl$ is geometrically irreducible. Hence the semisimplifications of $\cW_1$ and $\cW_2$ are irreducible, and the isomorphisms hold.
\end{proof}

\subsection{Symmetries}\label{subsec:symmetries}
For a permutation $\tau$ of $\{1,\dots,n\}$ we also write $\tau$ for the automorphism $(\tau h)_i=h_{\tau(i)}$ of $S_n$. Since $m_i\circ(\mathrm{id}\times\tau)=m_{\tau(i)}$, the permutation of the tensor factors is an isomorphism $\vartheta_\tau\colon(\mathrm{id}\times\tau)^*\cN_n\to\cN_n$. If $\sigma$ is an involution, then $\vartheta_\sigma\circ(\mathrm{id}\times\sigma)^*\vartheta_\sigma$ is the identity of $\cN_n$, since it permutes the tensor factors by $\sigma^2=\mathrm{id}$.

\begin{lemma}\label{lem:symmetry}
\begin{enumerate}[label=\textup{(\roman*)}]
\item For every permutation $\tau$, proper base change and $\vartheta_\tau$ induce isomorphisms of lisse sheaves on $S_n$ over $\Fp$
\[
\theta_\tau\colon\tau^*\cH_n\xrightarrow{\ \sim\ }\cH_n,\qquad\tau^*\cE_n\xrightarrow{\ \sim\ }\cE_n,\qquad\tau^*\cW_n\xrightarrow{\ \sim\ }\cW_n,
\]
and for an involution $\sigma$ one has $\theta_\sigma\circ\sigma^*\theta_\sigma=\mathrm{id}$.
\item For $k\in\Fp^\times$ write $\Sigma_{n,q}(h;\psi)$ for the sum \eqref{eq:Sigma} formed with $\psi$. Then $\Sigma_{n,q}(h;\psi^k)=\Sigma_{n,q}(kh;\psi)$; in particular the complex conjugate of $\Sigma_{n,q}(h;\psi)$, under any embedding of $\Q(\zeta_p)$ into $\C$, is $\Sigma_{n,q}(-h;\psi)$.
\end{enumerate}
\end{lemma}

\begin{proof}
(i) The isomorphism $\tau^*\cH_n\cong R^1\rho_!((\mathrm{id}\times\tau)^*\cN_n)$ of proper base change, followed by $R^1\rho_!(\vartheta_\tau)$, is $\theta_\tau$. It respects the weight filtrations, which are unique \cite[Th\'eor\`eme~3.4.1]{Deligne}, so it induces isomorphisms on $\cE_n$ and $\cW_n$. The relation $\vartheta_\sigma\circ(\mathrm{id}\times\sigma)^*\vartheta_\sigma=\mathrm{id}$ and the functoriality of proper base change give the last assertion. (ii) Substitute $u=kv$: then $\psi^k_q(\sum_iu_i^{-1})=\psi_q(\sum_iv_i^{-1})$ and $\langle h,u\rangle=\langle kh,v\rangle$. The last assertion is the case $k=-1$, as complex conjugation acts on $\psi$ as $\psi\mapsto\psi^{-1}$.
\end{proof}

\subsection{The Fourier model along rays}\label{subsec:model}
Let $\A^1_y$ and $\A^1_z$ be in duality through $(y,z)\mapsto yz$, let $\pr_y,\pr_z$ be the projections of $\A^1_y\times\A^1_z$, and let
\[
\FT(K)=R\pr_{z!}\bigl(\pr_y^*K\otimes\Lpsi{yz}\bigr)[1]
\]
be Laumon's Fourier transform \cite[(1.2.1.1)]{Laumon}. Let $j_y\colon\Gmn{y}\hookrightarrow\A^1_y$ be the inclusion, and let $\sigma\colon\Gmn{x}\to\Gmn{z}$ be the isomorphism $\sigma(x)=\lambda/x$; it extends to an isomorphism $\PP^1_x\to\PP^1_z$ with $\sigma(\infty)=0$ and $\sigma(0)=\infty$. Through $\sigma$ the inertia groups at $x=\infty$ and at $z=0$ are identified, compatibly with the ramification filtrations; for a representation $M$ of the latter we write $\sigma^*M$ for the corresponding representation of the former. For instance $\sigma^*\Lpsi{g/z}=\Lpsi{gx/\lambda}$ and $\sigma^*\Lchi{z}\cong\Lchi{x}$.

\begin{proposition}[{\cite[Proposition~4.14]{EbMom}}]\label{prop:ray-model}
Let $n\ge1$, $a\in(\Fq^\times)^n$, and $Q_a=\FT(j_{y!}\cK_a[1])$, a complex on $\A^1_z$ over $\Fq$.
\begin{enumerate}[label=\textup{(\roman*)}]
\item $Q_a$ is a perverse sheaf. Its cohomology sheaves are $\cH^{-1}(Q_a)=R^1\pr_{z!}(\pr_y^*j_{y!}\cK_a\otimes\Lpsi{yz})$, which has no nonzero section supported on a finite set, and $\cH^0(Q_a)$, which is supported at $z=0$. Their stalks at $z=0$ are $H^1_c(\Gm\otimes\Fqbar,\cK_a)$ and $H^2_c(\Gm\otimes\Fqbar,\cK_a)$.
\item Let $\cF_a$ be the restriction of $\cH^{-1}(Q_a)$ to $\Gmn{z}$. There is an isomorphism of lisse sheaves on $\Gmn{x}$ over $\Fq$
\[
\Delta_a^*\cH_n\cong\sigma^*\cF_a .
\]
\item The local monodromy of $\Delta_a^*\cH_n$ at $x=0$ is tame and unipotent, of the Jordan type of $\Up2^{\otimes n}$.
\end{enumerate}
\end{proposition}

\begin{proof}
(ii) Let $\cN'$ be the pull-back of $\cN_n$ along $\mathrm{id}\times\Delta_a\colon\Gmn{t}\times\Gmn{x}\to\Gmn{t}\times S_n$, so that $\cN'(t,x)=\bigotimes_i\Kl(a_ixt)\otimes\Lpsi{\lambda t}$. By Theorem~\ref{thm:recollection}(i) and proper base change, $\Delta_a^*\cH_n=R^1\pr_{x!}\cN'$ and $R^i\pr_{x!}\cN'=0$ for $i\ne1$. The automorphism $\phi(y,x)=(y/x,x)$ of $\Gm\times\Gm$ commutes with the projections to $\Gmn{x}$, and $\phi^*\cN'=\pr_y^*\cK_a\otimes\Lpsi{\lambda y/x}=(\mathrm{id}\times\sigma)^*\bigl(\pr_y^*\cK_a\otimes\Lpsi{yz}\bigr)$. Proper base change along $\sigma$ gives $\Delta_a^*\cH_n\cong\sigma^*R^1\pr_{z!}(\pr_y^*\cK_a\otimes\Lpsi{yz})$, where now $\pr_z\colon\Gmn{y}\times\Gmn{z}\to\Gmn{z}$; and integrating $j_{y!}\cK_a$ over $\A^1_y$ is the same as integrating $\cK_a$ over $\Gmn{y}$. All maps are defined over $\Fq$, since $\lambda\in\Fp^\times$. In particular $\cF_a$ is lisse of rank $2^n$, by Theorem~\ref{thm:recollection}(i).

(i) The complex $j_{y!}\cK_a[1]$ is a perverse sheaf \cite[Corollaire~4.1.3]{BBD}, so $Q_a$ is perverse \cite[Th\'eor\`eme~(1.3.2.3)]{Laumon}; hence $\cH^i(Q_a)=0$ for $i\ne-1,0$, $\cH^0(Q_a)$ is punctual and $\cH^{-1}(Q_a)$ has no nonzero section supported on a finite set \cite{BBD}. By definition $\cH^i(Q_a)=R^{i+2}\pr_{z!}(\pr_y^*j_{y!}\cK_a\otimes\Lpsi{yz})$, and by proper base change its stalk at $z$ is $H^{i+2}_c(\Gm\otimes\Fqbar,\cK_a\otimes\Lpsi{yz})$. For $z\ne0$ this vanishes for $i=0$: by Corollary~\ref{cor:Ka} and Lemma~\ref{lem:breaks}(ii), $\cK_a\otimes\Lpsi{yz}$ is isoclinic of break $1$ at $\infty$, hence has no nonzero geometric coinvariants. For $z=0$ we obtain the stated stalks.

(iii) This is \cite[Proposition~4.14]{EbMom}: by Laumon's principle of stationary phase \cite[Proposition~(2.3.3.1)(iii)]{Laumon}, $\cF_a(\infty)\cong\Flocz(\cK_a(0))\oplus\Flocii(\cK_a(\infty))$; the first summand already has dimension $2^n=\rank\cF_a$, so the second vanishes, and $\Flocz$ carries tame unipotent representations to tame unipotent representations with the same Jordan type \cite[Lemma~4.13]{EbMom}. Finally $\sigma$ identifies the inertia groups at $x=0$ and $z=\infty$.
\end{proof}
\section{Local Fourier transforms and vanishing cycles}\label{sec:local}

\subsection{Laumon's local Fourier transforms}\label{subsec:laumon}
In this section the ground field is $\Fpbar$. Let $I_0,I_\infty$ be the inertia groups at $y=0,\infty$ of $\PP^1_y$ and $I_{0'},I_{\infty'}$ those at $z=0,\infty$ of $\PP^1_z$. For such a group $I$ we write $\mathrm{Rep}(I)$ for the category of finite-dimensional continuous $\Qlb$-representations of $I$. Laumon defines local Fourier transforms
\[
\Flocz\colon\mathrm{Rep}(I_0)\to\mathrm{Rep}(I_{\infty'}),\qquad\Floc\colon\mathrm{Rep}(I_\infty)\to\mathrm{Rep}(I_{0'}),
\]
\[
\Flocii\colon\mathrm{Rep}(I_\infty)\to\mathrm{Rep}(I_{\infty'})
\]
through vanishing cycles of the extension by zero of $\mathcal L_\psi(\pi/\pi')$, $\mathcal L_\psi(\pi'/\pi)$, $\mathcal L_\psi(1/\pi\pi')$ on a product of two henselian traits with uniformizers $\pi,\pi'$ \cite[(2.4.1), D\'efinition~(2.4.2.3)]{Laumon}; they are exact functors \cite[Th\'eor\`eme~(2.4.3)]{Laumon}. We use the following facts.
\begin{enumerate}[label=(L\arabic*)]
\item $\FT$ carries perverse sheaves on $\A^1_y$ to perverse sheaves on $\A^1_z$ \cite[Th\'eor\`eme~(1.3.2.3)]{Laumon}.
\item Let $P$ be a perverse sheaf on $\A^1_y$, and let $V=\cH^{-1}(P)(\infty)\in\mathrm{Rep}(I_\infty)$. Let $R\Phi_0$ denote the vanishing cycles at $z=0$, for the henselization of $\A^1_z$ at $0$. Then there is an isomorphism of representations of $I_{0'}$
\[
R^{-1}\Phi_0\bigl(\FT(P)\bigr)\cong\Floc(V)
\]
\cite[Proposition~(2.3.2.1) and D\'efinition~(2.4.2.3)]{Laumon}; see also \cite[Proposition~6.3 and proof of Lemma~8.16]{EbCoh}, where the localization at $y=\infty$ underlying this isomorphism is carried out in detail for unramified $V$.
\item If all breaks of $V\in\mathrm{Rep}(I_\infty)$ are $<1$, then $\dim\Floc(V)=\dim V-\Swan(V)$ and $\Swan(\Floc(V))=\Swan(V)$ \cite[Th\'eor\`eme~(2.4.3)(ii)(b)]{Laumon}.
\item \textup{(Stationary phase and arithmetic form.)} For a lisse sheaf $\cG$ on $\Gmn{y}$ there is an isomorphism of representations of $I_{\infty'}$
\[
R^1\pr_{z!}\bigl(\pr_y^*j_{y!}\cG\otimes\Lpsi{yz}\bigr)(\infty')\cong\Flocz(\cG(0))\oplus\Flocii(\cG(\infty))
\]
\cite[Proposition~(2.3.3.1)(iii)]{Laumon}, \cite[Theorem~3]{Katz-Bourbaki}; one has $\dim\Flocz(U)=\dim U+\Swan(U)$ \cite[Th\'eor\`eme~(2.4.3)(i)(b)]{Laumon}. Over a finite field $k_0$, the same constructions define $\Flocz$ on representations of the decomposition groups; it is exact, compatible with restriction to inertia, commutes with twists by characters of $\Gal(\bar k_0/k_0)$, and if $\cG$ is defined over $k_0$ the displayed isomorphism is one of representations of the decomposition group at $z=\infty$ \cite[(2.4.1)--(2.4.3)]{Laumon}; see \cite[Theorem~4.11]{EbMom} for this form of the statements.
\end{enumerate}

By (L3), $\Floc$ carries a tame representation to a tame representation of the same dimension, and, by Lemma~\ref{lem:M}(a), for $c\ne0$ it carries $\Msf_e(c)$ to a representation of dimension one and Swan conductor one, that is, to a character of break $1$.

\subsection{Vanishing cycles of a Fourier transform at the origin}\label{subsec:vanishing}

\begin{lemma}\label{lem:vanishing}
Let $P$ be a perverse sheaf on $\A^1_y$, let $V=\cH^{-1}(P)(\infty)$, and let $Q=\FT(P)$. Let $U$ be a punctured neighbourhood of $0$ in $\A^1_z$ on which $Q$ is of the form $\cG[1]$ with $\cG$ lisse. Then there is an exact sequence of representations of $I_{0'}$
\[
0\longrightarrow\cH^{-1}(Q)_0\longrightarrow\cG(0)\longrightarrow\Floc(V)\longrightarrow\cH^0(Q)_0\longrightarrow0,
\]
in which $I_{0'}$ acts trivially on the two outer terms. Consequently $\cG(0)^{(r)}\cong\Floc(V)^{(r)}$ for every $r>0$, and there is an exact sequence of tame representations
\[
0\longrightarrow\cH^{-1}(Q)_0\longrightarrow\cG(0)^{(0)}\longrightarrow\Floc(V)^{(0)}\longrightarrow\cH^0(Q)_0\longrightarrow0 .
\]
\end{lemma}

\begin{proof}
Let $T$ be the henselization of $\A^1_z$ at $0$, with closed point $0$ and generic point $\eta$, and let $i\colon\{0\}\to T$ be the inclusion. The distinguished triangle $i^*Q\to R\Psi(Q)\to R\Phi(Q)\to$ of nearby and vanishing cycles \cite[Exp.~XIII, 2.1]{SGA7} is $I_{0'}$-equivariant, $I_{0'}$ acting trivially on $i^*Q$. By (L1), $Q$ is perverse; so $\cH^i(Q)=0$ for $i\ne-1,0$, and the sheaf $\cH^{-1}(Q)$ has no nonzero section supported at $0$. Since $Q|_\eta=\cG[1]|_\eta$, we have $R\Psi(Q)=\cG(0)[1]$. The long exact cohomology sequence reads
\[
0\to\cH^{-2}(R\Phi Q)\to\cH^{-1}(Q)_0\to\cG(0)\to\cH^{-1}(R\Phi Q)\to\cH^0(Q)_0\to0\to\cH^0(R\Phi Q)\to0,
\]
and $\cH^i(R\Phi Q)=0$ for $i\le-3$ and $i\ge1$. The map $\cH^{-1}(Q)_0\to\cG(0)$ is the specialization map of the sheaf $\cH^{-1}(Q)$; its kernel is the stalk at $0$ of the subsheaf of sections supported at $0$, which vanishes. Hence $\cH^{-2}(R\Phi Q)=\cH^0(R\Phi Q)=0$, and by (L2), $\cH^{-1}(R\Phi Q)\cong\Floc(V)$. This gives the first sequence. The break decomposition is exact, and the outer terms are tame; this gives the remaining assertions.
\end{proof}

\section{The end at infinity of a ray}\label{sec:rays}

\subsection{The decomposition at infinity}\label{subsec:decomp}
For a ray $\Delta_a$ we write $\cH_a(\infty)$ and $\cW_a(\infty)$ for the representations of the inertia group at $x=\infty$ on the geometric generic stalks of $\Delta_a^*\cH_n$ and $\Delta_a^*\cW_n$.

\begin{proposition}\label{prop:ray-decomp}
Let $n\ge1$ and $a\in(\Fpbar^\times)^n$.
\begin{enumerate}[label=\textup{(\roman*)}]
\item All breaks of $\cH_a(\infty)$ and of $\cW_a(\infty)$ are $0$ or $1$, and
\[
\cW_a(\infty)^{(1)}\cong\cH_a(\infty)^{(1)}\cong\bigoplus_{\substack{\varepsilon\in\mathsf E\\ L_\varepsilon(a)\ne0}}\sigma^*\Floc\bigl(\Msf_n(2L_\varepsilon(a))\bigr),
\]
each summand being a character of Swan conductor $1$.
\item $\dim\cH_a(\infty)^{(0)}=2^n-P(a)$ and $\dim\cW_a(\infty)^{(0)}=r(n)-P(a)$.
\item If $P(a)=2^{n-1}$, the inertia group at $x=\infty$ acts trivially on $\cH_a(\infty)^{(0)}$ and on $\cW_a(\infty)^{(0)}$.
\end{enumerate}
\end{proposition}

\begin{proof}
Choose a finite field $\Fq$ containing the coordinates of $a$, and apply Proposition~\ref{prop:ray-model} over $\Fpbar$. Let $P=j_{y!}\cK_a[1]$, so that $\FT(P)=Q_a$ and $V:=\cH^{-1}(P)(\infty)=\cK_a(\infty)$. By Corollary~\ref{cor:Ka}, $V=V^w\oplus V^t$ with
\[
V^w=\bigoplus_{L_\varepsilon(a)\ne0}\Msf_n\bigl(2L_\varepsilon(a)\bigr),\qquad V^t=\bigoplus_{L_\varepsilon(a)=0}\Msf_n(0),
\]
all of whose breaks are $<1$. By (L3) and the remark following (L4), $\Floc(V^w)$ is a sum of $P(a)$ characters of break $1$, and $\Floc(V^t)$ is tame. On a punctured neighbourhood of $0$, $Q_a=\cF_a[1]$ (Proposition~\ref{prop:ray-model}). By Lemma~\ref{lem:vanishing}, $\cF_a(0)^{(r)}\cong\Floc(V)^{(r)}$ for $r>0$; so the breaks of $\cF_a(0)$ are $0$ and $1$, its part of break $1$ is $\Floc(V^w)$, and its tame part has dimension $2^n-P(a)$. Transporting by $\sigma$ (Proposition~\ref{prop:ray-model}(ii)) gives the assertions on $\cH_a(\infty)$. By Theorem~\ref{thm:recollection}(iv), the inertia group at $x=\infty$ acts on the generic stalk of $\Delta_a^*\cE_n$ through a tame unipotent representation $\cE_a(\infty)$, of dimension $b_n$. The exact sequence $0\to\cE_a(\infty)\to\cH_a(\infty)\to\cW_a(\infty)\to0$ and the exactness of the break decomposition give $\cW_a(\infty)^{(1)}\cong\cH_a(\infty)^{(1)}$ and $\dim\cW_a(\infty)^{(0)}=2^n-P(a)-b_n=r(n)-P(a)$, and show that there are no other breaks. This proves (i) and (ii).

(iii) If $P(a)=2^{n-1}$, then $V^t=0$ and $\Floc(V)$ has no tame part. In the second exact sequence of Lemma~\ref{lem:vanishing} the third term vanishes, so $\cH^0(Q_a)_0=0$ and $\cF_a(0)^{(0)}\cong\cH^{-1}(Q_a)_0$, on which inertia acts trivially. Transporting by $\sigma$, inertia acts trivially on $\cH_a(\infty)^{(0)}$, hence on its quotient $\cW_a(\infty)^{(0)}=\cH_a(\infty)^{(0)}/\cE_a(\infty)$; here $\cE_a(\infty)\subseteq\cH_a(\infty)^{(0)}$ because $\cE_a(\infty)$ is tame.
\end{proof}

\subsection{Calibration}\label{subsec:calibration}

\begin{proposition}\label{prop:cal1}
For every $c\in\Fpbar^\times$, $\Floc(\Msf_1(c))\cong\Lpsi{-c^2/(4z)}$.
\end{proposition}

\begin{proof}
Let $n=1$, $s_1=c/2$ and $a=(s_1^2)$. Then $\mathsf E=\{(1)\}$, $L_{(1)}(a)=s_1\ne0$ and $P(a)=1=r(1)$. By Proposition~\ref{prop:ray-decomp}, $\cW_a(\infty)$ has no tame part, and $\cW_a(\infty)\cong\sigma^*\Floc(\Msf_1(2s_1))=\sigma^*\Floc(\Msf_1(c))$. On the other hand, by Theorem~\ref{thm:recollection}(v), $\Delta_a^*\cW_1\cong\Lpsi{-s_1^2x/\lambda}\otimes\Qlb(-1)$, so $\cW_a(\infty)\cong\Lpsi{-s_1^2x/\lambda}=\sigma^*\Lpsi{-s_1^2/z}$. Since $\sigma^*$ is an equivalence, $\Floc(\Msf_1(c))\cong\Lpsi{-s_1^2/z}=\Lpsi{-c^2/(4z)}$.
\end{proof}

\begin{proposition}\label{prop:cal2}
For every $c\in\Fpbar^\times$, $\Floc(\Msf_0(c))\cong\Lpsi{-c^2/(4z)}\otimes\Lchi{z}$.
\end{proposition}

\begin{proof}
Let $n=2$. Fix $u\in\Fpbar^\times$ with $2u=c$. For $v\in\Fpbar\smallsetminus\{0,u,-u\}$ put $s_1=(u+v)/2$, $s_2=(u-v)/2$ and $a=(s_1^2,s_2^2)$; then $s_1s_2\ne0$, $\mathsf E=\{(1,1),(1,-1)\}$, $L_{(1,1)}(a)=u$ and $L_{(1,-1)}(a)=v$. So $P(a)=2=r(2)$, and by Proposition~\ref{prop:ray-decomp}
\begin{equation}\label{eq:cal2-a}
\cW_a(\infty)\cong\sigma^*\Floc\bigl(\Msf_0(2u)\bigr)\oplus\sigma^*\Floc\bigl(\Msf_0(2v)\bigr).
\end{equation}
On the other hand, by Theorem~\ref{thm:recollection}(v), $\Delta_a^*\cW_2\cong\Lpsi{-(s_1^2+s_2^2)x/\lambda}\otimes f^*\Kl\otimes\Qlb(-1)$ with $f(x)=s_1^2s_2^2x^2/\lambda^2$. The map $f$ is the composite of the isomorphism $x\mapsto w=s_1s_2x/\lambda$ with the squaring map $w\mapsto w^2=y$; so $f^*$ identifies the inertia group at $x=\infty$ with $J$, and by Proposition~\ref{prop:Kl-infty},
\[
(f^*\Kl)(\infty)\cong\Res_J\Ind_J^{I_\infty}\xi_{2,1}=\xi_{2,1}\oplus\xi_{2,1}^\gamma\cong\Lpsi{\tfrac{2s_1s_2x}\lambda}\otimes\Lchi{x}\ \oplus\ \Lpsi{-\tfrac{2s_1s_2x}\lambda}\otimes\Lchi{x}.
\]
As $-(s_1^2+s_2^2)\pm2s_1s_2=-(s_1\mp s_2)^2$, we obtain
\begin{equation}\label{eq:cal2-b}
\cW_a(\infty)\cong\Lpsi{-\tfrac{u^2x}\lambda}\otimes\Lchi{x}\ \oplus\ \Lpsi{-\tfrac{v^2x}\lambda}\otimes\Lchi{x}.
\end{equation}
For $b\in\Fpbar$ write $\omega_b=\Lpsi{-bx/\lambda}\otimes\Lchi{x}$; for $b\ne b'$ the characters $\omega_b$ and $\omega_{b'}$ are not isomorphic, since $\omega_b\omega_{b'}^{-1}=\Lpsi{(b'-b)x/\lambda}$ is wild. By \eqref{eq:cal2-a}, \eqref{eq:cal2-b} and the Krull--Schmidt theorem for these semisimple representations, the character $\sigma^*\Floc(\Msf_0(2u))$, which depends only on $u$, belongs to $\{\omega_{u^2},\omega_{v^2}\}$ for every admissible $v$. Choose admissible $v_1,v_2$ with $v_1^2\ne v_2^2$; as $v_1^2,v_2^2\ne u^2$, the intersection $\{\omega_{u^2},\omega_{v_1^2}\}\cap\{\omega_{u^2},\omega_{v_2^2}\}$ is $\{\omega_{u^2}\}$. Hence $\sigma^*\Floc(\Msf_0(c))\cong\omega_{u^2}=\sigma^*\bigl(\Lpsi{-u^2/z}\otimes\Lchi{z}\bigr)$, and $u^2=c^2/4$.
\end{proof}

\begin{proof}[Proof of Theorem~B]
Since $\Msf_e(c)$ depends only on $e$ modulo $2$, Theorem~B is Propositions~\ref{prop:cal1} and~\ref{prop:cal2}.
\end{proof}

\begin{remark}\label{rem:Fu}
Theorem~B agrees with the case $r=2$, $s=1$, $\alpha(t)=ct$ of \cite[Theorem~0.1(iii)]{Fu}: in Fu's notation $\delta(t')=-c/(2t')$ and $\beta(t')=\alpha(t')+t'^2\delta(t')=ct'/2=-c^2/(4\delta(t'))$, and the Kummer factor $\mathcal K_{\chi_2}(-\frac12ct')$, multiplied by the Kummer sheaf of the inducing character, gives $\Lchi{z}^{\otimes(1-e)}$ geometrically. The proof above uses neither Fu's computation nor any explicit computation of vanishing cycles; it only uses Laumon's general results (L1)--(L3), the local monodromy of $\Kl$ at $\infty$ and the cases $n=1,2$ of the family.
\end{remark}

\subsection{Proof of Theorem A (i)--(iii)}\label{subsec:proofA}
By Proposition~\ref{prop:ray-decomp}(i) and Theorem~B, with $e\equiv n\pmod2$,
\[
\cW_a(\infty)^{(1)}\cong\bigoplus_{L_\varepsilon(a)\ne0}\sigma^*\Bigl(\Lpsi{-\tfrac{L_\varepsilon(a)^2}z}\otimes\Lchi{z}^{\otimes(1-e)}\Bigr)\cong\bigoplus_{L_\varepsilon(a)\ne0}\Lpsi{-\tfrac{L_\varepsilon(a)^2x}\lambda}\otimes\Lchi{x}^{\otimes(n+1)},
\]
since $\sigma^*\Lchi z\cong\Lchi x$ and $1-e\equiv n+1\pmod2$. This is (i). Statements (ii) and (iii) are Proposition~\ref{prop:ray-decomp}(ii),(iii): the Swan conductor of $\cW_a(\infty)$ is the number $P(a)$ of characters of break $1$. \qed

\subsection{Two auxiliary lemmas for compact-support cohomology}
\label{subsec:coinvariant-safeguards}

The following lemmas make explicit the two points at which a formal
Euler-characteristic computation would otherwise be insufficient.
They involve only the standard interpretation of top compactly supported
cohomology as coinvariants and the elementary structure of tame inertia.

\begin{lemma}[Purity of a geometric constant quotient]
\label{lem:pure-coinvariant}
Let $X/\Fq$ be geometrically connected, with an $\Fq$-rational point,
and let $\mathcal F$ be a lisse sheaf pure of weight $w$. Write $V$
for its geometric generic fibre and $V_{\pi_1(X_{\Fqbar})}$ for its
maximal geometrically constant quotient. Every eigenvalue of geometric
Frobenius on this quotient has complex absolute value $q^{w/2}$.
If $X$ is a smooth geometrically connected curve, then every eigenvalue
on $H_c^2(X_{\Fqbar},\mathcal F)$ has absolute value $q^{(w+2)/2}$.
\end{lemma}
\begin{proof}
The quotient $V\twoheadrightarrow V_{\pi_1(X_{\Fqbar})}$ is stable under
$\pi_1(X)$, because the geometric fundamental group is normal. At the
chosen rational point it induces a Frobenius-equivariant quotient of
$\mathcal F_x$. The eigenvalues on a quotient of a finite-dimensional
Frobenius module form a submultiset of those on the source, with algebraic
multiplicities; purity of $\mathcal F_x$ gives the assertion.
For the curve, Poincar\'e duality identifies
$H_c^2(X_{\Fqbar},\mathcal F)$ with the geometric coinvariants
$V_{\pi_1(X_{\Fqbar})}(-1)$, and the Tate twist multiplies geometric
Frobenius eigenvalues by $q$.
\end{proof}

\begin{lemma}[The order-four tame exclusion]
\label{lem:order-four-tame}
Suppose $p\ne2$, $I$ is the inertia group of a henselian trait with
algebraically closed residue field, and $J\subset I$ is the subgroup
of index two corresponding to the tame quadratic cover. Let $\xi$ be
the nontrivial quadratic character of $J/P$, where $P$ is wild inertia.
Then $\operatorname{Ind}_J^I\xi$ is the direct sum of two nontrivial
characters of order four of $I/P$. In particular it has no trivial
quotient, including after passage to any direct sum of its copies.
\end{lemma}
\begin{proof}
The tame quotient is procyclic. If $\tau$ is a topological generator,
then $J/P$ is generated by $\tau^2$ and $\xi(\tau^2)=-1$.
On the induced two-dimensional representation, the characteristic
polynomial of $\tau$ is $X^2+1$. Over $\Qlb$ it splits into the two
characters sending $\tau$ to $\sqrt{-1}$ and $-\sqrt{-1}$.
Both have order four, so no trivial character can be a quotient.
\end{proof}

\begin{remark}[Application and exact scope]
\label{rem:coinvariant-safeguards}
In the proof of Proposition~\ref{prop:global-euler}(ii),
Lemma~\ref{lem:pure-coinvariant} is applied over a finite extension
containing $a,c$, to $\mathcal K_a$ (weight $n$) and
$\Delta_a^*\mathcal W_n$ (weight $n+1$), both on $\Gm$.
It compares the resulting top-cohomology weights $n+2$ and $n+3$.
For $c\ne0$, the only possible trivial quotient at infinity comes from
$\Msf_{n+1}(0)$ in \eqref{eq:twist-infty}. If $n$ is even, the
index $n+1$ is odd and Lemma~\ref{lem:order-four-tame} rules out
this quotient. For $n$ odd this exclusion does not hold:
$\Msf_{n+1}(0)$ then contains the trivial tame character.
\end{remark}

\subsection{Global Euler characteristics}\label{subsec:global}
We now prove Theorem~A(iv). We give two proofs: the first deduces it from Theorem~A(i)--(ii); the second is global and uses only Corollary~\ref{cor:Ka}, so that it checks the Artin--Schreier exponents of Theorem~A independently of the local Fourier analysis of Sections~\ref{sec:local}--\ref{subsec:proofA}.

\begin{proposition}\label{prop:global-euler}
Let $n\ge1$, $a\in(\Fpbar^\times)^n$ and $c\in\Fpbar$, and let $\cL_c=\Lpsi{cx}\otimes\Delta_a^*\cW_n$ on $\Gmn{x}$.
\begin{enumerate}[label=\textup{(\roman*)}]
\item $\chi_c(\Gm\otimes\Fpbar,\cL_c)=-P(a)$ if $c=0$, and $\chi_c(\Gm\otimes\Fpbar,\cL_c)=-(r(n)-m_a(c))$ if $c\ne0$, where $m_a(c)=\#\{\varepsilon\in\mathsf E:L_\varepsilon(a)^2=\lambda c\}$.
\item $H^0_c(\Gm\otimes\Fpbar,\cL_c)=0$. If $c=0$, or $n$ is even, or $m_a(c)=0$, then also $H^2_c(\Gm\otimes\Fpbar,\cL_c)=0$, so that $\dim H^1_c(\Gm\otimes\Fpbar,\cL_c)=-\chi_c(\Gm\otimes\Fpbar,\cL_c)$.
\end{enumerate}
\end{proposition}

\begin{proof}[First proof of \textup{(i)}]
By Proposition~\ref{prop:ray-model}(iii) and Theorem~\ref{thm:recollection}(iv), $\Delta_a^*\cW_n$ is tame at $x=0$, and so is $\cL_c$. By the Grothendieck--Ogg--Shafarevich formula, $\chi_c(\Gm\otimes\Fpbar,\cL_c)=-\Swan_\infty(\cL_c)$. If $c=0$ this is $-P(a)$ by Theorem~A(ii). Let $c\ne0$. By Lemma~\ref{lem:breaks}(ii), $\Lpsi{cx}\otimes\cW_a(\infty)^{(0)}$ is isoclinic of break $1$, with Swan conductor $r(n)-P(a)$. By Theorem~A(i), $\Lpsi{cx}\otimes\cW_a(\infty)^{(1)}$ is the sum of the characters $\Lpsi{(c-L_\varepsilon(a)^2/\lambda)x}\otimes\Lchi x^{\otimes(n+1)}$, $L_\varepsilon(a)\ne0$, of which exactly $m_a(c)$ are tame and the others have Swan conductor $1$. Hence $\Swan_\infty(\cL_c)=r(n)-P(a)+P(a)-m_a(c)$.
\end{proof}

\begin{proof}[Second proof of \textup{(i)}, and proof of \textup{(ii)}]
Choose $\Fq$ containing $c$ and the coordinates of $a$. Let $\cN'$ and $\phi$ be as in the proof of Proposition~\ref{prop:ray-model}(ii). By the projection formula, $\Lpsi{cx}\otimes\Delta_a^*\cH_n=R^1\pr_{x!}(\cN'\otimes\Lpsi{cx})$, and $R^i\pr_{x!}(\cN'\otimes\Lpsi{cx})=0$ for $i\ne1$; so the Leray spectral sequence gives
\[
H^k_c\bigl(\Gmn{x}\otimes\Fqbar,\Lpsi{cx}\otimes\Delta_a^*\cH_n\bigr)\cong H^{k+1}_c\bigl(\Gm^2\otimes\Fqbar,\phi^*(\cN'\otimes\Lpsi{cx})\bigr),
\]
and $\phi^*(\cN'\otimes\Lpsi{cx})=\pr_y^*\cK_a\otimes\Lpsi{\lambda y/x+cx}$. We compute the right-hand side by projecting to $\Gmn{y}$. If $c\ne0$, the substitution $x=\lambda y v$ transforms $\Lpsi{\lambda y/x+cx}$ into $\Lpsi{v^{-1}+\lambda c\,y\,v}$, so by (K1) and proper base change $R\pr_{y!}\Lpsi{\lambda y/x+cx}\cong[\times\lambda c]^*\Kl[-1]$. If $c=0$, the automorphism $(y,x)\mapsto(y,\lambda y/x)$ transforms $\Lpsi{\lambda y/x}$ into the pull-back of $\Lpsi{x}$, and $R\pr_{y!}\Lpsi{\lambda y/x}$ is the constant sheaf with value $R\Gamma_c(\Gm\otimes\Fqbar,\Lpsi{x})=\Qlb[-1]$, on which Frobenius acts trivially, because $-\sum_{x\in\Fq^\times}\psi_q(x)=1$. By the projection formula,
\begin{equation}\label{eq:global-iso}
H^k_c\bigl(\Gmn{x}\otimes\Fqbar,\Lpsi{cx}\otimes\Delta_a^*\cH_n\bigr)\cong\begin{cases}H^k_c\bigl(\Gmn{y}\otimes\Fqbar,\cK_a\bigr),&c=0,\\ H^k_c\bigl(\Gmn{y}\otimes\Fqbar,\cK_a\otimes[\times\lambda c]^*\Kl\bigr),&c\ne0,\end{cases}
\end{equation}
compatibly with Frobenius. Both sheaves on $\Gmn{y}$ are tame at $0$, so by the Grothendieck--Ogg--Shafarevich formula their Euler characteristics are minus their Swan conductors at $\infty$. For $c=0$ this is $-P(a)$ (Corollary~\ref{cor:Ka}). For $c\ne0$, fix $\varsigma$ with $\varsigma^2=\lambda c$; by Corollary~\ref{cor:Ka}, Proposition~\ref{prop:Kl-infty} and Lemma~\ref{lem:M}(d),(e),
\begin{equation}\label{eq:twist-infty}
\bigl(\cK_a\otimes[\times\lambda c]^*\Kl\bigr)(\infty)\cong\bigoplus_{\varepsilon\in\mathsf E}\Bigl(\Msf_{n+1}\bigl(2(L_\varepsilon(a)+\varsigma)\bigr)\oplus\Msf_{n+1}\bigl(2(L_\varepsilon(a)-\varsigma)\bigr)\Bigr),
\end{equation}
whose Swan conductor is the number of pairs $(\varepsilon,\pm)$ with $L_\varepsilon(a)\pm\varsigma\ne0$, that is, $2^n-m_a(c)$, since each $\varepsilon$ with $L_\varepsilon(a)^2=\lambda c$ contributes exactly one vanishing term. On the other hand $\Delta_a^*\cE_n$ is a successive extension of geometrically constant sheaves of total rank $b_n$ (Theorem~\ref{thm:recollection}(iv)), and for a geometrically constant sheaf $\mathcal C$ of rank $m$ on $\Gm$ one has $\chi_c(\Gm\otimes\Fpbar,\Lpsi{cx}\otimes\mathcal C)=-m$ if $c\ne0$ and $=0$ if $c=0$. Additivity of Euler characteristics in $0\to\Delta_a^*\cE_n\to\Delta_a^*\cH_n\to\Delta_a^*\cW_n\to0$ gives (i): $-P(a)-0$ if $c=0$, and $-(2^n-m_a(c))+b_n=-(r(n)-m_a(c))$ if $c\ne0$.

(ii) $H^0_c$ vanishes on the affine curve $\Gm$. The functor $H^2_c(\Gm\otimes\Fpbar,-)$ is right exact on sheaves on a curve, so $H^2_c(\cL_c)$ is a quotient of $H^2_c(\Lpsi{cx}\otimes\Delta_a^*\cH_n)$, which by \eqref{eq:global-iso} is the space of geometric coinvariants, twisted by $(-1)$, of $\cK_a$ if $c=0$ and of $\cK_a\otimes[\times\lambda c]^*\Kl$ if $c\ne0$. Let first $c=0$. By Lemma~\ref{lem:pure-coinvariant}, the geometric coinvariants of $\cK_a$, pure of weight $n$, are pure of weight $n$; so $H^2_c(\Lpsi{0}\otimes\Delta_a^*\cH_n)$ is pure of weight $n+2$. But $H^2_c(\cL_0)$ is the space of geometric coinvariants of $\Delta_a^*\cW_n$, twisted by $(-1)$, which is pure of weight $n+3$ by Lemma~\ref{lem:pure-coinvariant}. A quotient of a pure space of weight $n+2$ which is pure of weight $n+3$ vanishes. Let now $c\ne0$. A nonzero geometric coinvariant quotient of $\cK_a\otimes[\times\lambda c]^*\Kl$ is a trivial quotient of its local monodromy at $\infty$. By \eqref{eq:twist-infty} and Lemma~\ref{lem:M}(a),(c), the tame part of that local monodromy is a sum of copies of $\Msf_{n+1}(0)$, one for each pair $(\varepsilon,\pm)$ with $L_\varepsilon(a)\pm\varsigma=0$; there are none if $m_a(c)=0$, and if $n$ is even, $\Msf_{n+1}(0)=\Msf_1(0)$ is a sum of two characters of order $4$ with no trivial quotient, by Lemma~\ref{lem:order-four-tame}. In these cases the coinvariants vanish, hence so does $H^2_c(\cL_c)$.
\end{proof}

\begin{proof}[Proof of Theorem~A(iv)]
This is Proposition~\ref{prop:global-euler}(i).
\end{proof}

\begin{remark}\label{rem:global-converse}
The second proof determines the Artin--Schreier exponents of Theorem~A(i) without the calibration of Section~\ref{subsec:calibration}, once one knows the shape of the decomposition given by Proposition~\ref{prop:ray-decomp}. Indeed, suppose that $\cW_a(\infty)^{(1)}\cong\bigoplus_{\varepsilon}\Lpsi{g_\varepsilon x}\otimes\theta_\varepsilon$ with $g_\varepsilon\in\Fpbar^\times$ and tame characters $\theta_\varepsilon$, the sum being over the $\varepsilon\in\mathsf E$ with $L_\varepsilon(a)\ne0$. As in the first proof, $\Swan_\infty(\cL_c)=r(n)-\#\{\varepsilon:g_\varepsilon=-c\}$ for $c\ne0$, and comparison with Proposition~\ref{prop:global-euler}(i) gives $\#\{\varepsilon:g_\varepsilon=-c\}=m_a(c)$ for every $c\ne0$, that is, the multisets $\{g_\varepsilon\}$ and $\{-L_\varepsilon(a)^2/\lambda\}$ coincide. Thus the constant $-\frac14$ of Theorem~B is forced by global Euler characteristics as well.
\end{remark}

\begin{remark}[A resonant obstruction, not an inferred cohomology dimension]
\label{rem:H2}
For odd $n$, $c\ne0$, and $m_a(c)>0$, the tame representation
$\Msf_{n+1}(0)$ contains the trivial character, and local monodromy
does not exclude a geometrically constant quotient of
$\mathcal F_{a,c}$. Thus the arguments proving
Proposition~\ref{prop:global-euler}(ii) in the nonresonant cases
cannot simply be extended to this locus. In particular, computing a
finite initial Taylor expansion of the $L$-function cannot establish
$H_c^2=0$ or its dimension. Proposition~\ref{prop:resonant-correction}
below gives an unconditional dimension bound and isolates the
possible top-weight contribution as a separate term.
\end{remark}

\subsection{Comments}\label{subsec:comments}

\begin{remark}[Stationary phase]\label{rem:stationary}
Let $a$ be off the critical-point escape divisor. At $h=xa$ the square roots of the coordinates are $\sqrt x\,s_i$, so $L_\varepsilon(xa)=\sqrt x\,L_\varepsilon(a)$, and by \cite[Proposition~3.6 and Remark~3.7]{EbCoh} the critical point attached to $\varepsilon$ is
\[
u_i=-\frac{\lambda\varepsilon_i}{x\,s_iL_\varepsilon(a)},\qquad t=\frac{x\,L_\varepsilon(a)^2}{\lambda^2},
\]
with critical value $-xL_\varepsilon(a)^2/\lambda$ and Hessian determinant $2^{n-1}\lambda/(t\,u_1^3\cdots u_n^3)=C_\varepsilon\,x^{3n-1}$, $C_\varepsilon\in\Fpbar^\times$. The principle of stationary phase predicts that an exponential sum in $n+1$ variables with nondegenerate critical points contributes, near infinity in the parameter, the Artin--Schreier sheaves of the critical values twisted by the quadratic characters of the Hessian determinants. Since $\chi_2(C_\varepsilon x^{3n-1})=\chi_2(x)^{n+1}$ geometrically, Theorem~A(i) is precisely this prediction for the wild part of $\cW_n$ along rays. The determinant (Theorem~C) is then the Artin--Schreier sheaf of the sum of the critical values, and the twists $\chi_2^{(n+1)P(a)}$ cancel, because $P(a)=2^{n-1}$ is even for $n\ge2$ off the critical-point escape divisor.
\end{remark}

\begin{remark}[Collisions]\label{rem:collisions}
If $L_{\varepsilon_0}(a)=0$, the critical point attached to $\varepsilon_0$ escapes to infinity as $h$ approaches $a$ \cite[Remark~5.4]{EbCoh}; on the ray through $a$ it does not exist. Theorem~A shows that the corresponding summand of $\cW_a(\infty)^{(1)}$ disappears and is replaced by tame summands: by the proof of Proposition~\ref{prop:ray-decomp}, the tame part of $\cH_a(\infty)$ then involves $\Floc(\Msf_n(0))$, a sum of two tame characters. The sheaf $\cW_n$ itself is lisse across the critical-point escape divisor \cite[Theorem~5.3 and Remark~5.4]{EbCoh}; the collision is visible only at the end of the rays, as a drop of the Swan conductor.
\end{remark}

\begin{remark}[Coordinate lines]\label{rem:coordinate}
Along a coordinate line $C_c=\{(h_1,c_2,\dots,c_n)\}$ the local monodromy of $\cW_n$ at $h_1=\infty$ is $T\oplus M$, where $T$ is tame and unipotent with $\kappa_j$ Jordan blocks of size $j-2$, $\kappa_j$ being the number of Jordan blocks of size $j$ of $\Up2^{\otimes(n-1)}$, and $M\cong\Lpsi{-h_1/\lambda}\otimes X$ is isoclinic of break $1$ with $\dim X=2^{n-1}$ and $X$ of breaks $<1$ \cite[Theorem~6.12]{EbCoh}. Along a ray off the critical-point escape divisor, by contrast, the tame part is trivial (Theorem~A(iii)) and the part of break one splits into $2^{n-1}$ explicit characters. The two descriptions give the same determinant.
\end{remark}
\section{The geometric determinant}\label{sec:det}

\subsection{Characters of the fundamental group of a torus}\label{subsec:torus}
Let $k$ be an algebraically closed field of characteristic $p$, and $T=\Gm^n\otimes k$ with coordinates $h_1,\dots,h_n$. For $1\le i\le n$ and $c\in(k^\times)^{n-1}$ let $\ell_{i,c}\subset T$ be the coordinate line on which $h_j=c_j$ for $j\ne i$, with coordinate $h_i$; for $a\in(k^\times)^n$ let $\Delta_a\colon\Gmn x\to T$ be the ray through $a$. For a continuous character $\eta$ of $\pi_1(T)$ we write $\eta|_{\ell_{i,c}}$ and $\Delta_a^*\eta$ for the characters of $\pi_1(\Gm)$ obtained by composition; base points play no role, as the characters are abelian.

\begin{lemma}\label{lem:torus}
Let $\eta\colon\pi_1(T)\to\Qlb^\times$ be a continuous character of finite order. Suppose that
\begin{enumerate}
\item[\textup{(B)}] for every $i$ and every $c\in(k^\times)^{n-1}$, the character $\eta|_{\ell_{i,c}}$ is unramified at $h_i=0$, and
\item[\textup{(R)}] for every $a\in(k^\times)^n$, the character $\Delta_a^*\eta$ is trivial.
\end{enumerate}
Then $\eta$ is trivial.
\end{lemma}

\begin{proof}
Write $\eta=\eta'\eta_p$ with $\eta'$ of order prime to $p$ and $\eta_p$ of order a power of $p$; both are powers of $\eta$, so they satisfy (B) and (R).

\emph{The part of order prime to $p$.} Let $m$ be the order of $\eta'$. By the Kummer sequence, $\mathrm{Pic}(T)=0$ and the $m$-divisibility of $k^\times$, the group $H^1(T,\mu_m)$ is isomorphic to $(\Z/m)^n$, with basis the Kummer classes of $h_1,\dots,h_n$ \cite[Proposition~8.6, Step~2]{EbCoh}. Hence $\eta'=\bigotimes_i\mathcal L_{\chi_i}(h_i)$ for characters $\chi_i$ of $\mu_m$. On $\ell_{i,c}$ the functions $h_j$, $j\ne i$, are constant, so $\eta'|_{\ell_{i,c}}\cong\mathcal L_{\chi_i}(h_i)$, which is unramified at $h_i=0$ only if $\chi_i$ is trivial, $h_i$ being a uniformizer there. By (B), $\eta'$ is trivial.

\emph{Artin--Schreier classes.} Let $A=k[h_1^{\pm1},\dots,h_n^{\pm1}]$ and $\wp(f)=f^p-f$. The Artin--Schreier sequence $0\to\Fp\to\mathcal O\xrightarrow{\wp}\mathcal O\to0$ on the \'etale site of the affine scheme $T$, and $H^1(T,\mathcal O)=0$, give an isomorphism $A/\wp(A)\cong H^1(T,\Fp)$, under which the class of $f$ corresponds, after the choice of $\psi$, to the character of $\Lpsi{f}$ \cite[Ch.~III, Proposition~3.7]{Milne}. It is compatible with pull-back along morphisms of affine schemes; in particular with restriction to $\ell_{i,c}$ and pull-back along $\Delta_a$, where the same description holds with $A$ replaced by $k[h_i^{\pm1}]$, resp.\ $k[x^{\pm1}]$. Let $\mathsf M=\Z^n\smallsetminus p\Z^n$. Every class in $A/\wp(A)$ has a unique representative $f=\sum_{m\in\mathsf M}f_mh^m$ (a finite sum), which we call its normal form. Indeed, constants lie in $\wp(k)$, and $c\,h^{pm}\equiv c^{1/p}h^m$ modulo $\wp(A)$, which gives existence. For uniqueness, suppose $\sum_{m\in\mathsf M}f_mh^m=\wp(g)$ with $g=\sum_mg_mh^m$. If $g$ is not constant, choose $m_0\ne0$ in the support of $g$ with $\max_i|m_{0,i}|$ maximal; then $pm_0$ is not in the support of $g$, so the coefficient of $h^{pm_0}$ in $\wp(g)=\sum_mg_m^ph^{pm}-g$ is $g_{m_0}^p\ne0$, whereas $pm_0\notin\mathsf M$. So $g$ is constant, $\wp(g)$ is constant, and all $f_m$ vanish. The same holds in one variable.

We record two consequences in one variable. Let $g=\sum_{e\in\Z\smallsetminus p\Z}g_eu^e$ be a normal form in $k[u^{\pm1}]$. (a) $\Lpsi g$ is trivial on $\pi_1(\Gm\otimes k)$ if and only if $g=0$, by uniqueness. (b) $\Lpsi g$ is unramified at $u=0$ if and only if $g_e=0$ for all $e<0$. Indeed, write $g=g_++g_-$, with $g_+$ the part with positive exponents. The sheaf $\Lpsi{g_+}$ is lisse on $\A^1$; if $\Lpsi g$ is unramified at $0$, then so is $\Lpsi{g_-}$, which is lisse on $\PP^1\smallsetminus\{0\}$; hence $\Lpsi{g_-}$ is lisse on $\PP^1$, hence trivial, and $g_-=0$ by (a). The converse is clear.

\emph{The part of order $p$.} Suppose that $\eta_p$ has order dividing $p$, and let $f=\sum_{m\in\mathsf M}f_mh^m$ be the normal form of the corresponding class. Fix $i$, and for $m\in\Z^n$ write $m=(m_i,m')$ with $m'\in\Z^{n-1}$ the other coordinates. For $c\in(k^\times)^{n-1}$ and $e\in\Z$ put $\varphi_e(c)=\sum_{m_i=e}f_mc^{m'}$, the sum over the support of $f$. The restriction of $f$ to $\ell_{i,c}$ is $\sum_e\varphi_e(c)h_i^e$, and its normal form has, for $e_0\in\Z\smallsetminus p\Z$, the coefficient $\sum_{s\ge0}\varphi_{e_0p^s}(c)^{1/p^s}$ at $h_i^{e_0}$. By (B) and (b), this vanishes for every $e_0<0$, $p\nmid e_0$, and every $c$. Raising to a power $p^S$ with $S$ large, the Laurent polynomial
\[
\sum_{s=0}^S\sum_{m_i=e_0p^s}f_m^{p^{S-s}}c^{p^{S-s}m'}
\]
vanishes on $(k^\times)^{n-1}$, hence is zero. The monomials $c^{p^{S-s}m'}$ occurring in it are pairwise distinct: for fixed $s$ the vector $m$ is determined by $m'$, since $m_i=e_0p^s$; and if $p^{S-s}m'=p^{S-\tilde s}\tilde m'$ with $s<\tilde s$, then $\tilde m=p^{\tilde s-s}m\in p\Z^n$, which is impossible. Hence $f_m=0$ whenever $m_i=e_0p^s$ for some $s\ge0$ and some negative $e_0$ prime to $p$, that is, whenever $m_i<0$. As $i$ was arbitrary, $f\in k[h_1,\dots,h_n]$, and $f$ has no constant term.

Now let $a\in(k^\times)^n$. Then $\Delta_a^*f=\sum_mf_ma^mx^{|m|}$ with $|m|=m_1+\cdots+m_n\ge1$. For $d_0\ge1$ prime to $p$ the coefficient of $x^{d_0}$ in its normal form is $\sum_{s\ge0}\Psi_{d_0p^s}(a)^{1/p^s}$, where $\Psi_d(a)=\sum_{|m|=d}f_ma^m$. By (R) and (a) these coefficients vanish for all $a$, and the argument of the preceding paragraph, with the monomials $a^{p^{S-s}m}$, which are pairwise distinct for the same reason, shows that $f_m=0$ for all $m$. Hence $\eta_p$ is trivial.

\emph{The part of order $p^r$.} By induction on $r\ge1$: if $\eta_p$ has order dividing $p^r$, then $\eta_p^p$ satisfies (B) and (R) and has order dividing $p^{r-1}$, so it is trivial by induction; hence $\eta_p$ has order dividing $p$, and it is trivial by the previous step.
\end{proof}

\begin{remark}
Neither hypothesis can be dropped. The character of $\Lpsi{h_1/h_2}$, or of $\mathcal L_\chi(h_1/h_2)$, is trivial on all rays, and the character of $\Lpsi{h_1+\cdots+h_n}$ is unramified along all coordinate lines. The argument is a variant of \cite[Lemma~8.10]{EbCoh}, where triviality on all coordinate lines is assumed.
\end{remark}

\subsection{The determinant along a ray}\label{subsec:det-ray}

\begin{proposition}\label{prop:det-ray}
Let $n\ge1$ and $a\in(\Fpbar^\times)^n$. Then the lisse sheaf $\det(\Delta_a^*\cW_n)\otimes\Lpsi{2^{n-1}\lambda^{-1}(a_1+\cdots+a_n)\,x}$ on $\Gmn x\otimes\Fpbar$ is geometrically constant.
\end{proposition}

\begin{proof}
Let $\cD$ be the sheaf in question, of rank one. At $x=0$ the local monodromy of $\Delta_a^*\cW_n$ is unipotent (Proposition~\ref{prop:ray-model}(iii) and Theorem~\ref{thm:recollection}(iv)), so its determinant is trivial on $I_0$; and $\Lpsi{\,\cdot\,x}$ is lisse at $0$. Hence $\cD$ extends to a lisse sheaf $\bar\cD$ on $\A^1_x$. At $x=\infty$, the break decomposition gives $\cW_a(\infty)=\cW_a(\infty)^{(0)}\oplus\cW_a(\infty)^{(1)}$, and by Theorem~A(i)
\[
\det\cW_a(\infty)\cong\det\cW_a(\infty)^{(0)}\otimes\Lpsi{-\frac{x}{\lambda}\sum_{L_\varepsilon(a)\ne0}L_\varepsilon(a)^2}\otimes\Lchi x^{\otimes(n+1)P(a)} .
\]
Now
\[
\sum_{L_\varepsilon(a)\ne0}L_\varepsilon(a)^2=\sum_{\varepsilon\in\mathsf E}L_\varepsilon(a)^2=\sum_{i,j}s_is_j\sum_{\varepsilon\in\mathsf E}\varepsilon_i\varepsilon_j=2^{n-1}\sum_is_i^2=2^{n-1}(a_1+\cdots+a_n),
\]
because $\sum_{\varepsilon\in\mathsf E}\varepsilon_i\varepsilon_j=0$ for $i\ne j$: if one of $i,j$ equals $1$ the sum is $\sum_{\varepsilon\in\mathsf E}\varepsilon_{j'}$ for some $j'\ne1$, and otherwise it is a sum of products of two independent signs. Hence $\bar\cD(\infty)\cong\det\cW_a(\infty)^{(0)}\otimes\Lchi x^{\otimes(n+1)P(a)}$ is tame. By the Grothendieck--Ogg--Shafarevich formula, $\chi_c(\A^1\otimes\Fpbar,\bar\cD)=1-\Swan_\infty(\bar\cD)=1$. As $H^0_c(\A^1\otimes\Fpbar,\bar\cD)=0$, the group $H^2_c(\A^1\otimes\Fpbar,\bar\cD)$, which is the space of geometric coinvariants of $\bar\cD$ twisted by $(-1)$, is nonzero; so the rank-one sheaf $\bar\cD$ is geometrically constant.
\end{proof}

\subsection{Proof of Theorem C}\label{subsec:proofC}
Let $\cA=\Lpsi{-2^{n-1}\lambda^{-1}(h_1+\cdots+h_n)}$ and $\eta=\det\cW_n\otimes\cA^{-1}$, a lisse sheaf of rank one on $S_n$ over $\Fp$; we regard it as a character of $\pi_1(S_n\otimes\Fpbar)$. Its image is finite, because the images of $\det\cW_n$ (Theorem~\ref{thm:recollection}(vi)) and of $\cA$ are finite.

\emph{Condition (B).} Let $c\in(\Fpbar^\times)^{n-1}$ and $C_c=\ell_{1,c}$. By \cite[Theorem~6.8(ii)]{EbCoh}, the local monodromy of $\cW_n|_{C_c}$ at $h_1=0$ is tame and unipotent. The image of the inertia group at $h_1=0$ is therefore contained in the closure of the group generated by a unipotent automorphism, on which the determinant is $1$; so $\det\cW_n|_{C_c}$ is unramified at $h_1=0$. The sheaf $\cA$ is lisse on $\A^n$. Hence $\eta|_{\ell_{1,c}}$ is unramified at $h_1=0$. For $i\ne1$, let $\tau$ be the transposition of $1$ and $i$; it maps $\ell_{1,c}$ isomorphically onto a coordinate line in the direction $i$, every such line arises in this way, and $\tau^*\cW_n\cong\cW_n$ (Lemma~\ref{lem:symmetry}(i)) and $\tau^*\cA=\cA$. So $\eta$ satisfies (B).

\emph{Condition (R).} For $a\in(\Fpbar^\times)^n$, $\Delta_a^*\eta=\det(\Delta_a^*\cW_n)\otimes\Lpsi{2^{n-1}\lambda^{-1}(a_1+\cdots+a_n)x}$, which is geometrically trivial by Proposition~\ref{prop:det-ray}.

By Lemma~\ref{lem:torus}, $\eta$ is geometrically trivial. A lisse sheaf of rank one on $S_n$ that is geometrically trivial is given by a character of $\pi_1(S_n)/\pi_1(S_n\otimes\Fpbar)=\Gal(\Fpbar/\Fp)$, that is, it is $\alpha^{\deg}$ with $\alpha$ the image of the geometric Frobenius. This proves the first two assertions. Finally, $2^{n-1}\lambda^{-1}\ne0$ in $\Fp$, so the restriction of $\cA$ to any coordinate line is a nontrivial Artin--Schreier sheaf; its geometric image is a nontrivial subgroup of $\mu_p$, hence equal to $\mu_p$, and so is the geometric image of $\det\cW_n=\cA\otimes\alpha^{\deg}$. \qed

\begin{remark}
Proposition~\ref{prop:det-ray} uses only Theorem~A, and Lemma~\ref{lem:torus} is elementary; the only input from the coordinate boundary is the unipotence of the local monodromy at $h_i=0$, which enters through (B). In \cite{EbCoh}, Theorem~C is deduced instead from the local monodromy at $h_1=\infty$ along coordinate lines, which requires Laumon's stationary phase at a second singular point of a Fourier model \cite[Theorem~6.12]{EbCoh}.
\end{remark}
\section{The arithmetic normalization}\label{sec:arith}

Throughout this section $n\ge1$, $N=N(n)=(n+1)r(n)/2$, $\cA=\Lpsi{-2^{n-1}\lambda^{-1}(h_1+\cdots+h_n)}$, and $\alpha$ is the constant of Theorem~C. We regard $\Q(\zeta_p)$ as a subfield of $\Qlb$ through $\psi$, with $\zeta_p=\psi(1)$. By Theorem~C, for every finite extension $\Fq=\mathbb F_{p^f}$ of $\Fp$ and every $h\in S_n(\Fq)$,
\begin{equation}\label{eq:alpha}
\det\bigl(\Frob_h\mid\cW_n\bigr)=\psi_q\bigl(-2^{n-1}\lambda^{-1}(h_1+\cdots+h_n)\bigr)\,\alpha^f .
\end{equation}

\subsection{Integrality, purity and the product formula}\label{subsec:product}

\begin{lemma}\label{lem:integral}
Let $h\in S_n(\Fq)$.
\begin{enumerate}[label=\textup{(\roman*)}]
\item The characteristic polynomial of $\Frob_h$ on $\cW_{n,h}$ has coefficients in $\Q(\zeta_p)$.
\item Its roots are algebraic integers.
\item For every field isomorphism $\iota\colon\Qlb\to\C$, every root $\beta$ satisfies $|\iota(\beta)|=q^{(n+1)/2}$.
\end{enumerate}
\end{lemma}

\begin{proof}
(i) For $k\ge1$, $\Frob_h^k$ is the Frobenius of $h$ regarded as a point over $\mathbb F_{q^k}$, so by \eqref{eq:trace} its trace on $\cW_{n,h}$ is $(-1)^{n-1}q^k\Sigma_{n,q^k}(h)-\sum_{j\ge1}\nu_{n,j}q^{kj}\in\Q(\zeta_p)$. The coefficients of the characteristic polynomial are polynomials with rational coefficients in these traces, by Newton's identities.

(ii) The space $\cW_{n,h}$ is a quotient of $\cH_{n,h}$, which is isomorphic, compatibly with Frobenius, to $H^{n+1}_c(X\otimes\Fqbar,\Lpsi{F})$, where $X=\Gm^n\times\Gmn t$ over $\Fq$ and $F=F_n(\cdot,\cdot;h)$ (Theorem~\ref{thm:recollection}(i)). Let $\varpi\colon\tilde X\to X$ be the Artin--Schreier covering $\tilde X=\{(u,t,\vartheta):\vartheta^p-\vartheta=F(u,t)\}$, a finite \'etale Galois covering with group $\Fp$. Then $\varpi_*\Qlb$ is the direct sum of the sheaves $\Lpsi{jF}$, $j\in\Fp$, so $H^{n+1}_c(X\otimes\Fqbar,\Lpsi F)$ is a Frobenius-stable direct summand of $H^{n+1}_c(\tilde X\otimes\Fqbar,\Qlb)$, on which the eigenvalues of Frobenius are algebraic integers by Deligne's integrality theorem \cite[Exp.~XXI, Corollaire~5.5.3]{SGA7}.

(iii) $\cN_n$ is $\iota$-pure of weight $n$ for every $\iota$, by (K1), and so $\cW_{n,h}=H^1(\PP^1,j_*\cN_h)$ is $\iota$-pure of weight $n+1$ for every $\iota$ \cite[Th\'eor\`eme~3.2.3]{Deligne}.
\end{proof}

\begin{proposition}\label{prop:product}
$\alpha\in\Z[\zeta_p]$, and $\alpha\,p^{-N}$ is a root of unity; hence $\alpha p^{-N}\in\{\pm\zeta_p^k:0\le k<p\}$.
\end{proposition}

\begin{proof}
Let $h\in S_n(\Fp)$, for instance $h=(1,\dots,1)$. By \eqref{eq:alpha} with $q=p$, $\alpha=\det(\Frob_h\mid\cW_{n,h})\,\psi(2^{n-1}\lambda^{-1}\sum_ih_i)$. By Lemma~\ref{lem:integral}(i),(ii), $\det(\Frob_h\mid\cW_{n,h})$ is an algebraic integer of $\Q(\zeta_p)$, so it lies in $\Z[\zeta_p]$, the ring of integers of $\Q(\zeta_p)$ \cite[Ch.~2]{Washington}; hence $\alpha\in\Z[\zeta_p]$. Every embedding $\tau\colon\Q(\zeta_p)\to\C$ extends to a field isomorphism $\iota\colon\Qlb\to\C$, both fields being algebraically closed of the same transcendence degree over $\Q(\zeta_p)$; so by Lemma~\ref{lem:integral}(iii), $|\tau(\alpha)|=p^{N}$ for every $\tau$, and $|N_{\Q(\zeta_p)/\Q}(\alpha)|=p^{N(p-1)}$. The prime $p$ is totally ramified in $\Q(\zeta_p)$: $p\Z[\zeta_p]=\mathfrak p^{p-1}$ with $\mathfrak p=(1-\zeta_p)$ of norm $p$ \cite[Ch.~1--2]{Washington}. Since $\alpha$ is integral and its norm is a power of $p$, the ideal $\alpha\Z[\zeta_p]$ is a power $\mathfrak p^m$ of $\mathfrak p$, and comparing norms gives $m=N(p-1)$, so $\alpha\Z[\zeta_p]=p^N\Z[\zeta_p]$. Thus $\alpha p^{-N}$ is a unit of $\Z[\zeta_p]$ all of whose conjugates have absolute value one; by Kronecker's theorem \cite{Kronecker} it is a root of unity. The roots of unity of $\Q(\zeta_p)$ are the $\pm\zeta_p^k$ \cite[Ch.~2]{Washington}.
\end{proof}

\begin{proposition}\label{prop:galois}
$\alpha\in\Q$. Consequently $\alpha=\pm p^N$.
\end{proposition}

\begin{proof}
For a nontrivial additive character $\psi'$ of $\Fp$ let $\cW_n^{\psi'}$ and $\alpha(\psi')$ be the objects of Theorems~0 and~C formed with $\psi'$ in place of $\psi$ ($\lambda$ being fixed), so that $\alpha=\alpha(\psi)$. Let $k\in\Fp^\times$ and let $\sigma_k\in\Gal(\Q(\zeta_p)/\Q)$ be the automorphism with $\sigma_k(\zeta_p)=\zeta_p^k$, so that $\sigma_k\circ\psi=\psi^k$. Let $h\in S_n(\Fp)$.

First, $\sigma_k(\alpha(\psi))=\alpha(\psi^k)$. Indeed, by the proof of Lemma~\ref{lem:integral}(i), the characteristic polynomial of $\Frob_h$ on $\cW_{n,h}^{\psi'}$ is a polynomial with rational coefficients in the numbers $\Sigma_{n,p^j}(h;\psi')$, $j\ge1$, and $\sigma_k(\Sigma_{n,p^j}(h;\psi))=\Sigma_{n,p^j}(h;\psi^k)$ since $\sigma_k\circ\psi_{p^j}=(\psi^k)_{p^j}$. So $\sigma_k$ maps $\det(\Frob_h\mid\cW_{n,h}^\psi)$ to $\det(\Frob_h\mid\cW_{n,h}^{\psi^k})$, and \eqref{eq:alpha}, applied to $\psi$ and to $\psi^k$, gives $\sigma_k(\psi(-2^{n-1}\lambda^{-1}\sum h_i))\,\sigma_k(\alpha(\psi))=\psi^k(-2^{n-1}\lambda^{-1}\sum h_i)\,\alpha(\psi^k)$, whence the claim.

Second, $\alpha(\psi^k)=\alpha(\psi)$. Indeed, by Lemma~\ref{lem:symmetry}(ii), $\Sigma_{n,p^j}(h;\psi^k)=\Sigma_{n,p^j}(kh;\psi)$ for all $j$; so the Frobenius traces of $\cW_n^{\psi^k}$ at $h$ and of $\cW_n^\psi$ at $kh$ coincide over all extensions, and $\det(\Frob_h\mid\cW_{n,h}^{\psi^k})=\det(\Frob_{kh}\mid\cW_{n,kh}^{\psi})$. By \eqref{eq:alpha} this reads $\psi^k(-2^{n-1}\lambda^{-1}\sum h_i)\,\alpha(\psi^k)=\psi(-2^{n-1}\lambda^{-1}k\sum h_i)\,\alpha(\psi)$, and the two roots of unity are equal.

Hence $\alpha$ is fixed by $\Gal(\Q(\zeta_p)/\Q)$, so $\alpha\in\Q$, and $\alpha p^{-N}$ is a rational root of unity by Proposition~\ref{prop:product}.
\end{proof}

\subsection{A pairing between the fibres at \texorpdfstring{$h$}{h} and \texorpdfstring{$-h$}{-h}}\label{subsec:reflection}
Let $X=\Gmn t\times S_n$, $\rho\colon X\to S_n$ the projection, $[-1]\colon S_n\to S_n$ the automorphism $h\mapsto-h$, and $\iota\colon X\to X$ the involution $\iota(t,h)=(-t,-h)$, so that $\rho\circ\iota=[-1]\circ\rho$. Since $m_i\circ\iota=m_i$, we have
\[
\iota^*\cN_n=\cN_n^-:=\Bigl(\bigotimes_{i=1}^nm_i^*\Kl\Bigr)\otimes\Lpsi{-\lambda t}.
\]
Fix an isomorphism $\det\Kl\cong\Qlb(-1)$ of sheaves over $\Fp$ (K4), and let $\beta\colon\Kl\otimes\Kl\to\Lambda^2\Kl\cong\Qlb(-1)$ be the resulting perfect alternating pairing. Let $c\colon\Lpsi{\lambda t}\otimes\Lpsi{-\lambda t}\to\Qlb$ be the canonical isomorphism; since $\mathcal L_\psi$ is a character sheaf of the commutative group $\mathbb G_a$, these canonical isomorphisms are compatible with the commutativity constraint, so that $\iota^*c=c\circ\sw$ as maps $\Lpsi{-\lambda t}\otimes\Lpsi{\lambda t}\to\Qlb$. Put
\[
\pi=\Bigl(\bigotimes_{i=1}^nm_i^*\beta\Bigr)\otimes c\colon\ \cN_n\otimes\cN_n^-\longrightarrow\Qlb(-n),
\]
a perfect pairing of lisse sheaves on $X$ over $\Fp$, the $i$-th factor of $\cN_n$ being paired with the $i$-th factor of $\cN_n^-$.

\begin{lemma}\label{lem:sign}
Under the identification $\iota^*\cN_n^-=\cN_n$, one has $\iota^*\pi=(-1)^n\,\pi\circ\sw$ as pairings $\cN_n^-\otimes\cN_n\to\Qlb(-n)$.
\end{lemma}

\begin{proof}
Since $m_i\circ\iota=m_i$, we have $\iota^*(m_i^*\beta)=m_i^*\beta$. So $\iota^*\pi$ pairs the $i$-th factor $m_i^*\Kl$ of its first argument with the $i$-th factor of its second argument through $m_i^*\beta$, whereas $\pi\circ\sw$ pairs them through $m_i^*\beta\circ\sw=-m_i^*\beta$, as $\beta$ is alternating; and $\iota^*c=c\circ\sw$.
\end{proof}

\emph{The cup product.} Let $X_2=X\times_{S_n}X$, with projections $p_1,p_2$, structure map $\rho_2$ and diagonal $\delta$, a closed immersion. By the K\"unneth formula \cite[Exp.~XVII, 5.4.3]{SGA4}, for complexes $\mathcal A,\mathcal B$ on $X$ there is an isomorphism $R\rho_!\mathcal A\otimes^{\mathbf L}R\rho_!\mathcal B\cong R\rho_{2!}(p_1^*\mathcal A\otimes p_2^*\mathcal B)$ compatible with base change; composing with $R\rho_{2!}$ of the adjunction $p_1^*\mathcal A\otimes p_2^*\mathcal B\to\delta_*\delta^*(p_1^*\mathcal A\otimes p_2^*\mathcal B)=\delta_*(\mathcal A\otimes\mathcal B)$ and with $R\rho_{2!}\delta_*=R\rho_!$, we obtain the cup product $R\rho_!\mathcal A\otimes^{\mathbf L}R\rho_!\mathcal B\to R\rho_!(\mathcal A\otimes\mathcal B)$, which at a geometric point is the cup product on compactly supported cohomology of the fibre \cite[proof of Proposition~4.7]{EbMom}. By Theorem~\ref{thm:recollection}(i) and proper base change for the cartesian square formed by $\iota$ and $[-1]$, $R\rho_!\cN_n=\cH_n[-1]$ and $R\rho_!\cN_n^-\cong[-1]^*R\rho_!\cN_n=[-1]^*\cH_n[-1]$. Taking $\cH^2$ of the cup product for $\mathcal A=\cN_n$, $\mathcal B=\cN_n^-$, and composing with $R^2\rho_!(\pi)$ and the trace isomorphism $R^2\rho_!\Qlb(-n)\cong\Qlb(-n-1)$, we obtain a morphism of lisse sheaves on $S_n$ over $\Fp$
\[
\mathcal R_{\cH}\colon\cH_n\otimes[-1]^*\cH_n\longrightarrow\Qlb(-n-1).
\]
At a geometric point $h$ of $S_n$ let $\iota_h\colon\Gm\to\Gm$, $t\mapsto-t$; it carries the fibre of $\rho$ over $h$ to the fibre over $-h$, and $\iota_h^*\cN_{-h}=\cN^-_h$, the restriction of $\cN_n^-$ to the fibre over $h$. The fibre of $\mathcal R_\cH$ at $h$ is
\begin{equation}\label{eq:R-fibre}
\mathcal R_h(a,b)=\Tr\,\pi_{h*}\bigl(a\cup\iota_h^*b\bigr),\qquad a\in\cH_{n,h}=H^1_c(\Gm,\cN_h),\quad b\in\cH_{n,-h}=H^1_c(\Gm,\cN_{-h}),
\end{equation}
where $a\cup\iota_h^*b\in H^2_c(\Gm,\cN_h\otimes\cN_h^-)$, $\pi_h$ is the restriction of $\pi$, and $\Tr\colon H^2_c(\Gm,\Qlb(-n))\to\Qlb(-n-1)$ is the trace map; all cohomology groups are over $\Fpbar$.

\begin{proposition}\label{prop:reflection}
\begin{enumerate}[label=\textup{(\roman*)}]
\item The morphism $\mathcal R_\cH$ vanishes on $\cE_n\otimes[-1]^*\cH_n$ and on $\cH_n\otimes[-1]^*\cE_n$, and it induces a morphism of lisse sheaves on $S_n$ over $\Fp$
\[
\mathcal R\colon\cW_n\otimes[-1]^*\cW_n\longrightarrow\Qlb(-n-1)
\]
which is a perfect pairing on every fibre.
\item For every geometric point $h$, $a\in\cH_{n,h}$ and $b\in\cH_{n,-h}$, one has $\mathcal R_{-h}(b,a)=(-1)^{n+1}\mathcal R_h(a,b)$. Equivalently, $[-1]^*\mathcal R=(-1)^{n+1}\mathcal R\circ\sw$.
\item For every permutation $\tau$ of $\{1,\dots,n\}$, geometric point $x$, $y\in\cH_{n,\tau x}$ and $y'\in\cH_{n,-\tau x}$, one has $\mathcal R_{\tau x}(y,y')=\mathcal R_x(\theta_{\tau,x}y,\theta_{\tau,-x}y')$, where $\theta_{\tau,x}\colon\cH_{n,\tau x}\to\cH_{n,x}$ is the fibre at $x$ of the isomorphism $\theta_\tau$ of Lemma~\ref{lem:symmetry}.
\end{enumerate}
\end{proposition}

\begin{proof}
(i) Fix a geometric point $h$ and write $U=\Gm$, $\cN=\cN_h$, $\cN^-=\cN^-_h$. Since $\pi_h$ is perfect, $\cN^-\cong\cN^\vee(-n)$. For $a\in H^1_c(U,\cN)$ and $c\in H^1_c(U,\cN^-)$, the class $a\cup c$ equals $a\cup\phi^-(c)$, where $\phi^-\colon H^1_c(U,\cN^-)\to H^1(U,\cN^-)$ forgets supports and the second cup product is $H^1_c(U,\cN)\otimes H^1(U,\cN^-)\to H^2_c(U,\cN\otimes\cN^-)$. By Poincar\'e duality on the smooth curve $U$, $(a,c')\mapsto\Tr\,\pi_{h*}(a\cup c')$ is a perfect pairing between $H^1_c(U,\cN)$ and $H^1(U,\cN^-)$. Hence $c$ is orthogonal to all of $H^1_c(U,\cN)$ if and only if $\phi^-(c)=0$. In the same way $a\cup c=\phi(a)\cup c$, with $\phi\colon H^1_c(U,\cN)\to H^1(U,\cN)$, and Poincar\'e duality between $H^1(U,\cN)$ and $H^1_c(U,\cN^-)$ shows that $a$ is orthogonal to all of $H^1_c(U,\cN^-)$ if and only if $\phi(a)=0$. Now let $\cG$ be a lisse sheaf on $U$ with $\cG^{I_\infty}=0$. The forgetful map $H^1_c(U,\cG)\to H^1(U,\cG)$ is the composite of $H^1(\PP^1,j_!\cG)\to H^1(\PP^1,j_*\cG)$ with the edge map $H^1(\PP^1,j_*\cG)\to H^1(U,\cG)$ of the Leray spectral sequence, which is injective; and the kernel of the first map is the image of $H^0$ of the skyscraper sheaf $j_*\cG/j_!\cG$, that is, of $\cG^{I_0}$. Both $\cN$ and $\cN^-$ are isoclinic of break $1$ at $\infty$ (Lemma~\ref{lem:breaks}(ii)), so they have no nonzero $I_\infty$-invariants. Thus the left kernel of $\mathcal R_h$ is the image of $\cN_h^{I_0}$, which is $\cE_{n,h}$ (Theorem~\ref{thm:recollection}(ii)); and, since $\iota_h$ fixes $0$ and $\infty$ and $\iota_h^*$ carries $\cN_{-h}^{I_0}$ onto $(\cN^-_h)^{I_0}$, the right kernel is the image of $\cN_{-h}^{I_0}$, which is $\cE_{n,-h}$. The morphism $\cH_n\to([-1]^*\cH_n)^\vee(-n-1)$ adjoint to $\mathcal R_\cH$ is a morphism of lisse sheaves on the connected normal scheme $S_n$, so its kernel is a lisse subsheaf of $\cH_n$, whose fibre at every geometric point is the left kernel computed above. A lisse subsheaf of $\cH_n$ is determined by its fibre at one geometric point, as a $\pi_1(S_n)$-stable subspace; hence the kernel is $\cE_n$. In the same way the kernel of the other adjoint morphism is $[-1]^*\cE_n$. Therefore $\mathcal R_\cH$ induces $\mathcal R$, and $\mathcal R$ is perfect on every fibre.

(ii) Let $a\in H^1_c(\Gm,\cN_h)$ and $b\in H^1_c(\Gm,\cN_{-h})$. By \eqref{eq:R-fibre}, $\mathcal R_{-h}(b,a)=\Tr\,\pi_{-h*}(b\cup\iota_{-h}^*a)$. Pull-back by the automorphism $t\mapsto-t$ of $\Gm$ is compatible with cup products and with the trace map, since this automorphism extends to an automorphism of $\PP^1$, which acts trivially on $H^2(\PP^1,\Qlb)=H^2_c(\Gm,\Qlb)$. As $\iota_h^*\iota_{-h}^*a=a$, we obtain
\[
\mathcal R_{-h}(b,a)=\Tr\,(\iota_h^*\pi_{-h})_*\bigl(\iota_h^*b\cup a\bigr)=(-1)^n\,\Tr\,\pi_{h*}\bigl(\sw_*(\iota_h^*b\cup a)\bigr)
\]
by Lemma~\ref{lem:sign}, restricted to the fibre over $h$. The cup product of two classes of degree one is anticommutative (Koszul sign rule): $\sw_*(\iota_h^*b\cup a)=-\,a\cup\iota_h^*b$. Hence $\mathcal R_{-h}(b,a)=(-1)^{n+1}\mathcal R_h(a,b)$. The second formulation is the first one, read at all geometric points, together with $[-1]\circ[-1]=\mathrm{id}$; a morphism of lisse sheaves is determined by its fibres.

(iii) The permutation of the tensor factors identifies $(\mathrm{id}\times\tau)^*\pi$ with $\pi\circ(\vartheta_\tau\otimes\vartheta_\tau^-)$, where $\vartheta^-_\tau$ is defined for $\cN_n^-$ as $\vartheta_\tau$ is for $\cN_n$, since $\pi$ is the tensor product of the pairings of corresponding factors; and $\mathrm{id}\times\tau$ commutes with $\iota$. The assertion follows from the compatibility of the cup product and of the trace map with base change along $\tau$.
\end{proof}

\begin{corollary}\label{cor:reflection}
\begin{enumerate}[label=\textup{(\roman*)}]
\item $[-1]^*\cW_n\cong\cW_n^\vee(-n-1)$ as lisse sheaves on $S_n$ over $\Fp$.
\item $\alpha^2=p^{2N}$; hence $\alpha=\pm p^N$.
\end{enumerate}
\end{corollary}

\begin{proof}
(i) is Proposition~\ref{prop:reflection}(i). (ii) Taking determinants in (i) and using Theorem~C: $\det([-1]^*\cW_n)=[-1]^*\cA\otimes\alpha^{\deg}=\cA^{-1}\otimes\alpha^{\deg}$, and $\det(\cW_n^\vee(-n-1))=\cA^{-1}\otimes(\alpha^{-1}p^{(n+1)r(n)})^{\deg}$. Hence $\alpha=\alpha^{-1}p^{2N}$.
\end{proof}

\begin{remark}
Corollary~\ref{cor:reflection}(i) can also be obtained from the trace identity $\Tr(\Frob_{-h}\mid\cW_n)=\overline{\Tr(\Frob_h\mid\cW_n)}$, which follows from \eqref{eq:trace} and Lemma~\ref{lem:symmetry}(ii), together with the Chebotarev density theorem and the irreducibility of $\cW_n$ (Theorem~\ref{thm:recollection}(vii)). The cup product gives more: an explicit pairing, defined over $\Fp$, with a definite symmetry, which is what the next subsection uses.
\end{remark}

\subsection{Symmetric parameters}\label{subsec:symmetric}
Let $n$ be even and let $\sigma$ be a fixed-point-free involution of $\{1,\dots,n\}$. Let $\Sigma_\sigma\subset S_n$ be the closed subscheme defined by the equations $h_{\sigma(i)}+h_i=0$, $1\le i\le n$; choosing one index in each orbit of $\sigma$ identifies it with $\Gm^{n/2}$ over $\Fp$. On $\Sigma_\sigma$ the automorphisms $[-1]$ and $\sigma$ of $S_n$ agree: $-h=\sigma h$ for $h\in\Sigma_\sigma$.

\begin{theorem}\label{thm:symplectic}
Let $n$ be even and $\sigma$ as above, and define
\[
\mathcal B=\mathcal R\circ\bigl(\mathrm{id}\otimes\theta_\sigma^{-1}\bigr)\colon\ \cW_n|_{\Sigma_\sigma}\otimes\cW_n|_{\Sigma_\sigma}\longrightarrow\Qlb(-n-1),
\]
where $\theta_\sigma^{-1}\colon\cW_n|_{\Sigma_\sigma}\to(\sigma^*\cW_n)|_{\Sigma_\sigma}=([-1]^*\cW_n)|_{\Sigma_\sigma}$. Then $\mathcal B$ is a perfect alternating pairing of lisse sheaves on $\Sigma_\sigma$ over $\Fp$. Consequently:
\begin{enumerate}[label=\textup{(\roman*)}]
\item $\ggeom(\cW_n|_{\Sigma_\sigma})$ is contained in the symplectic group of the form induced by $\mathcal B$ on the geometric generic fibre;
\item for every $h\in\Sigma_\sigma(\Fq)$, $\Frob_h$ is a symplectic similitude of $\cW_{n,h}$ with multiplier $q^{n+1}$; hence $\det(\Frob_h\mid\cW_{n,h})=q^{N}$, and $P_h(T)=\det(1-\Frob_hT\mid\cW_{n,h})$ satisfies $P_h(T)=q^NT^{r(n)}P_h\bigl(1/(q^{n+1}T)\bigr)$.
\end{enumerate}
\end{theorem}

\begin{proof}
$\mathcal B$ is perfect because $\mathcal R$ is (Proposition~\ref{prop:reflection}(i)) and $\theta_\sigma$ is an isomorphism; it is defined over $\Fp$ because $\Sigma_\sigma$, $\mathcal R$ and $\theta_\sigma$ are. Let $h$ be a geometric point of $\Sigma_\sigma$. The fibre of $\theta_\sigma^{-1}$ at $h$ is $\theta_{\sigma,h}^{-1}=\theta_{\sigma,-h}\colon\cW_{n,h}\to\cW_{n,-h}$, by the relation $\theta_\sigma\circ\sigma^*\theta_\sigma=\mathrm{id}$ of Lemma~\ref{lem:symmetry}(i) read at $h$, since $\sigma(-h)=h$. So $\mathcal B_h(a,a')=\mathcal R_h(a,\theta_{\sigma,-h}a')$ for $a,a'\in\cW_{n,h}$. By Proposition~\ref{prop:reflection}(ii),
\[
\mathcal B_h(a',a)=\mathcal R_h(a',\theta_{\sigma,-h}a)=(-1)^{n+1}\mathcal R_{-h}(\theta_{\sigma,-h}a,a').
\]
By Proposition~\ref{prop:reflection}(iii) with $\tau=\sigma$ and $x=h$, so that $\tau x=-h$ and $-\tau x=h$, applied to $y=\theta_{\sigma,-h}a$ and $y'=a'$,
\[
\mathcal R_{-h}(\theta_{\sigma,-h}a,a')=\mathcal R_h(\theta_{\sigma,h}\theta_{\sigma,-h}a,\theta_{\sigma,-h}a')=\mathcal R_h(a,\theta_{\sigma,-h}a')=\mathcal B_h(a,a').
\]
Hence $\mathcal B_h(a',a)=(-1)^{n+1}\mathcal B_h(a,a')=-\mathcal B_h(a,a')$, as $n$ is even; in characteristic zero this means that $\mathcal B_h$ is alternating. (Propositions~\ref{prop:reflection}(ii),(iii) are stated for $\cH_n$; they descend to $\cW_n$ by Proposition~\ref{prop:reflection}(i).)

(i) The geometric fundamental group acts trivially on $\Qlb(-n-1)$, so it preserves the alternating form $\mathcal B$ on the geometric generic fibre, and so does the Zariski closure of its image.

(ii) Since $\mathcal B$ is a morphism of sheaves over $\Fp$, $\mathcal B_h(\Frob_ha,\Frob_ha')=q^{n+1}\mathcal B_h(a,a')$. If $g$ is a symplectic similitude with multiplier $\mu$ of a symplectic space of dimension $2m$ with Gram matrix $\mathsf J$, then $g^{\mathsf T}\mathsf Jg=\mu\mathsf J$, and taking Pfaffians gives $\det(g)\Pf(\mathsf J)=\mu^m\Pf(\mathsf J)$, so $\det g=\mu^m$. Here $2m=r(n)$ and $\mu^m=q^{(n+1)r(n)/2}=q^N$. Moreover $g^{-1}=\mu^{-1}\mathsf J^{-1}g^{\mathsf T}\mathsf J$ is conjugate to $\mu^{-1}g^{\mathsf T}$, so the multiset of eigenvalues of $g$ is stable under $\beta\mapsto\mu/\beta$; with $\det g=\mu^m$ and $r(n)$ even this gives the functional equation of $P_h$.
\end{proof}

\begin{corollary}\label{cor:symplectic}
If $n$ is even, then $\alpha=p^N$.
\end{corollary}

\begin{proof}
Let $\sigma=(1\,2)(3\,4)\cdots(n-1\;n)$ and $h_0=(1,-1,1,-1,\dots,1,-1)\in\Sigma_\sigma(\Fp)$. Since $\sum_ih_{0,i}=0$, \eqref{eq:alpha} gives $\alpha=\det(\Frob_{h_0}\mid\cW_{n,h_0})$, which is $p^N$ by Theorem~\ref{thm:symplectic}(ii).
\end{proof}

\begin{remark}
For $n=2$, $\Sigma_{(1\,2)}=\{(h,-h)\}$ and $\cW_2|_{\Sigma_{(1\,2)}}\cong[h\mapsto-h^2/\lambda^2]^*\Kl\otimes\Qlb(-1)$ by Theorem~\ref{thm:recollection}(v), whose geometric monodromy group is $\SL_2=\Sp_2$ \cite[Theorem~11.1]{Katz-GKM}. We do not know whether $\ggeom(\cW_n|_{\Sigma_\sigma})$ is the full symplectic group for even $n\ge4$.
\end{remark}

\subsection{The tame end of a ray}\label{subsec:tame-end}
For $n\ge1$ put $S_1(n)=\sum_{k=0}^{\lfloor n/2\rfloor}k\,\nu_{n,k}$.

\begin{lemma}[{\cite[Lemma~8.14(i)]{EbCoh}}]\label{lem:numerology}
$2S_1(n)=(n+1)b_n-2^n$; consequently $n2^{n-1}-S_1(n)=N(n)$.
\end{lemma}

\begin{proof}
Let $K=\lfloor n/2\rfloor$. Summation by parts gives $S_1(n)=\sum_{k\le K}k\binom nk-\sum_{j\le K-1}(j+1)\binom nj=K\binom nK-\sum_{j<K}\binom nj$. If $n=2K$, the symmetry of binomial coefficients gives $\sum_{j<K}\binom nj=\frac12(2^n-b_n)$; if $n=2K+1$, it gives $\sum_{j\le K}\binom nj=2^{n-1}$, so $\sum_{j<K}\binom nj=2^{n-1}-b_n$. In both cases $2S_1(n)=(n+1)b_n-2^n$. Then $2(n2^{n-1}-S_1(n))=(n+1)(2^n-b_n)=2N(n)$.
\end{proof}

\begin{proposition}[{cf.~\cite[Proposition~4.18(iii)]{EbMom}}]\label{prop:tame-end}
Let $a\in(\Fq^\times)^n$, regard $\Delta_a$ as a morphism over $\Fq$, let $D$ be a decomposition group at the point $x=0$ of $\A^1_x$ over $\Fq$, and let $\Phi\in D$ map to the geometric Frobenius of $\Fq$. Then the determinant of $\Phi$ on the geometric generic stalk of $\Delta_a^*\cW_n$ is $q^{N(n)}$. \textup{(}It does not depend on the choice of $\Phi$: the inertia group acts unipotently on this stalk, by Proposition~\textup{\ref{prop:ray-model}(iii)}, since $\Delta_a^*\cW_n$ is a quotient of $\Delta_a^*\cH_n$.\textup{)}
\end{proposition}

\begin{proof}
\emph{The Fourier model over $\Fq$.} By Proposition~\ref{prop:ray-model}(ii), $\Delta_a^*\cH_n\cong\sigma^*\cF_a$ over $\Fq$, and $\sigma$ identifies $D$ with a decomposition group $D_{\infty'}$ at $z=\infty$ over $\Fq$, compatibly with the projections to $\Gal(\Fqbar/\Fq)$. Let $D_0$ be a decomposition group at $y=0$ over $\Fq$. By (L4), as representations of $D_{\infty'}$,
\[
\cF_a(\infty')\cong\Flocz\bigl(\cK_a(0)\bigr)\oplus\Flocii\bigl(\cK_a(\infty)\bigr),
\]
where $\cK_a(0)$ is regarded as a representation of $D_0$. Since $\Swan_0(\cK_a)=0$, the first summand has dimension $2^n=\rank\cF_a$, so the second vanishes.

\emph{A filtration of $\cK_a(0)$.} By (K3), the representation $([\times a_i]^*\Kl)(0)$ of $D_0$ contains the line spanned by $e_0$, which is the space of $I_0$-invariants and is therefore $D_0$-stable; $I_0$ acts trivially on it and on the quotient, and $\Phi$ acts on them by $1$ and by $q$. As $D_0$ is topologically generated by $I_0$ and $\Phi$, the line is the trivial character $\Qlb$ of $D_0$ and the quotient is $\Qlb(-1)$. Taking tensor products, $\cK_a(0)$ has a $D_0$-stable filtration whose graded pieces are the characters $\Qlb(-k)$, $0\le k\le n$, each with multiplicity $\binom nk$.

\emph{Transport.} Let $\omega=\Flocz(\Qlb)$, a character of $D_{\infty'}$. Apply (L4) to $\cG=\Qlb$ on $\Gmn y$ over $\Fq$: the automorphism $(y,z)\mapsto(yz,z)$ of $\Gmn y\times\Gmn z$ and the K\"unneth formula show that $R^1\pr_{z!}\Lpsi{yz}$ is, on $\Gmn z$, the constant sheaf with value $H^1_c(\Gm\otimes\Fqbar,\Lpsi y)$, a line on which Frobenius acts by $-\sum_{y\in\Fq^\times}\psi_q(y)=1$. So $\omega\oplus\Flocii(\Qlb)$ is the trivial one-dimensional representation of $D_{\infty'}$; as $\dim\omega=1$ by (L4), $\Flocii(\Qlb)=0$ and $\omega=\Qlb$. Since $\Flocz$ is exact and commutes with twists by characters of $\Gal(\Fqbar/\Fq)$, $\cF_a(\infty')$, and hence the generic stalk of $\Delta_a^*\cH_n$ as a representation of $D$, has a filtration whose graded pieces are the characters $\Qlb(-k)$ with multiplicity $\binom nk$. Therefore $\det(\Phi\mid\Delta_a^*\cH_n)=q^{\sum_kk\binom nk}=q^{n2^{n-1}}$.

\emph{The lower-weight part.} By Theorem~\ref{thm:recollection}(iv), the generic stalk of $\Delta_a^*\cE_n$ has a $D$-stable filtration, induced by the weight filtration of $\cE_n$, whose graded pieces are the stalks of the geometrically constant sheaves $\Delta_a^*\gr^W_{2k}\cE_n$; $D$ acts on them through $\Gal(\Fqbar/\Fq)$, and $\Phi$ acts with the single eigenvalue $q^k$, with multiplicity $\nu_{n,k}$. So $\det(\Phi\mid\Delta_a^*\cE_n)=q^{S_1(n)}$, and by Lemma~\ref{lem:numerology}, $\det(\Phi\mid\Delta_a^*\cW_n)=q^{n2^{n-1}-S_1(n)}=q^{N(n)}$.
\end{proof}

\begin{corollary}\label{cor:tame-end}
For every $n\ge1$, $\alpha=p^N$.
\end{corollary}

\begin{proof}
Let $a\in(\Fp^\times)^n$ and $q=p$. By Theorem~C, $\det(\Delta_a^*\cW_n)\cong\Lpsi{-2^{n-1}\lambda^{-1}(\sum_ia_i)x}\otimes\alpha^{\deg}$, a lisse sheaf on $\A^1_x$ over $\Fp$; so the inertia group at $x=0$ acts trivially on its generic stalk, and $\Phi$ acts by the Frobenius at the rational point $x=0$, that is, by $\psi(0)\alpha=\alpha$. By Proposition~\ref{prop:tame-end}, $\alpha=p^N$.
\end{proof}

\subsection{Proofs of Theorem D and Corollary E}\label{subsec:proofD}

\begin{proof}[Proof of Theorem~D]
Statement (i) is Propositions~\ref{prop:product} and~\ref{prop:galois}; (ii) is Proposition~\ref{prop:reflection} and Corollary~\ref{cor:reflection}; (iii) is Theorem~\ref{thm:symplectic} and Corollary~\ref{cor:symplectic}; (iv) is Proposition~\ref{prop:tame-end} and Corollary~\ref{cor:tame-end}. Each of Corollaries~\ref{cor:symplectic} and~\ref{cor:tame-end} gives $\alpha=p^N$, the first for even $n$ and the second for all $n$.
\end{proof}

\begin{proof}[Proof of Corollary~E]
Let $G=\ggeom(\cW_n)$ and let $W$ be the geometric generic fibre. Under the standing hypothesis $p>\max(2,n)$, by \cite[Theorem~8.8]{EbCoh}, $G^0=\SL(W)$. By Theorem~C the image of $\pgeom(S_n)$ under $\det$ is $\mu_p$; as $\det^{-1}(\mu_p)$ is Zariski closed, $\det(G)=\mu_p$. If $\det(g)\in\mu_p$, choose $g'\in G$ with $\det g'=\det g$; then $g'^{-1}g\in\SL(W)\subseteq G$, so $g\in G$. Hence $G=\{g:\det(g)^p=1\}$. Let $\gamma\in\Qlb$ with $\gamma^2=p^{-(n+1)}$, and $\gamma=p^{-(n+1)/2}$ if $n$ is odd. Then $\gamma^{r(n)}=p^{-N}$ (for even $n$ because $r(n)$ is even), so by Theorems~C and~D, $\det(\cW_n\otimes\gamma^{\deg})=\cA$, whose image in $\Qlb^\times$ is $\mu_p$. Hence $\garith(\cW_n\otimes\gamma^{\deg})\subseteq\{g:\det(g)^p=1\}=G=\ggeom(\cW_n\otimes\gamma^{\deg})\subseteq\garith(\cW_n\otimes\gamma^{\deg})$. Finally, by Theorem~0 and \cite[Corollary~2.16(ii)]{EbMom}, the eigenvalues of $\Frob_q$ on $H^{n-1}_c(Y_h\otimes\Fqbar,\Lpsi{\sum_iu_i^{-1}})$ are those on $\cW_{n,h}$ divided by $q$, together with the numbers $q^{k-1}$, $1\le k\le\lfloor n/2\rfloor$, with multiplicity $\nu_{n,k}$. Hence its determinant is $q^{N-r(n)}\cdot q^{S_1(n)-(b_n-1)}$ times $\psi_q(-2^{n-1}\lambda^{-1}\sum_ih_i)$, and $N-r(n)+S_1(n)-b_n+1=n2^{n-1}-2^n+1=(n-2)2^{n-1}+1$ by Lemma~\ref{lem:numerology}.
\end{proof}

\section{The critical-value polynomial and the escape locus}\label{sec:critical-polynomial}

This section makes explicit an algebraic invariant which is implicit in the
signed-square-root geometry of the parent manuscript \cite[Sections~5 and~13]{EbFull}.
It also distinguishes two phenomena which should not be conflated:
\emph{escape of a critical point} ($L_\varepsilon=0$) and
\emph{coincidence of distinct critical values}
($L_\varepsilon^2=L_\delta^2$ with both terms nonzero).
Only the first phenomenon removes a summand from the wild part in Theorem~A.

Write $m=2^{n-1}$, choose formal square roots $s_i$ with $s_i^2=a_i$,
let $\mathsf E=\{\varepsilon\in\{\pm1\}^n:\varepsilon_1=1\}$,
and put $A(a)=\sum_i a_i$ and $B(a)=\sum_i a_i^2$. Define
\begin{equation}\label{eq:critical-polynomial}
 \mathscr C_{n,a}(T)=\prod_{\varepsilon\in\mathsf E}
 \left(T+\frac{L_\varepsilon(a)^2}{\lambda}\right),
 \qquad L_\varepsilon(a)=\sum_i\varepsilon_i s_i.
\end{equation}
The roots of $\mathscr C_{n,a}$ are the critical values
$-L_\varepsilon(a)^2/\lambda$, including a zero for every index at
which the corresponding critical point escapes the torus.

\begin{theorem}[Critical-value polynomial and its first coefficients]
\label{thm:critical-polynomial}
Let $p$ be odd, $n\ge2$, $\lambda\ne0$, and $a_i\in\Fpbar$.
Then $\mathscr C_{n,a}(T)$ has coefficients polynomial in
$a_1,\ldots,a_n$ over $\Fp[\lambda^{-1}]$; in particular it is
independent of the choices of the $s_i$. Its first two coefficients are
\begin{align}
 [T^{m-1}]\mathscr C_{n,a}(T)&=\frac{m}{\lambda}A(a),\\
 [T^{m-2}]\mathscr C_{n,a}(T)&=\frac{m}{2\lambda^2}
 \bigl((m-3)A(a)^2+2B(a)\bigr).
\label{eq:critical-two-coefficients}
\end{align}
Moreover the product
\begin{equation}\label{eq:escape-discriminant}
 \mathscr D_n(a)=\prod_{\varepsilon\in\mathsf E}L_\varepsilon(a)^2
\end{equation}
is the square of a polynomial in $a_1,\ldots,a_n$.
The critical-point escape locus is exactly $\mathscr D_n(a)=0$.
For $a_i\ne0$, the reduction in the Swan conductor is
\begin{equation}\label{eq:swan-drop}
 m-\Swan_\infty(\Delta_a^*\cW_n)
 =\#\{\varepsilon\in\mathsf E:L_\varepsilon(a)=0\}.
\end{equation}
\end{theorem}

\begin{proof}
Changing the sign of any $s_j$, $j\ge2$, permutes the linear forms
$L_\varepsilon$ indexed by $\mathsf E$. Changing the sign of $s_1$
permutes their negatives: the form associated with $\varepsilon$ is sent
to $-L_\delta$, where $\delta_1=1$ and
$\delta_j=-\varepsilon_j$ for $j\ge2$.
Thus the multiset of squares is invariant under every independent
sign change. Since $p\ne2$, the invariant polynomial ring
$\Fp[s_1,\ldots,s_n]^{(\mathbb Z/2)^n}$ is
$\Fp[s_1^2,\ldots,s_n^2]$, which proves the polynomial assertion.

For the coefficients we expand the signed sums. Independence of the
$n-1$ freely varying signs gives
\begin{equation}\label{eq:signed-moments}
 \sum_{\varepsilon\in\mathsf E}L_\varepsilon^2=mA,
 \qquad
 \sum_{\varepsilon\in\mathsf E}L_\varepsilon^4
 =m\left(\sum_i a_i^2+6\sum_{i<j}a_ia_j\right)
 =m(3A^2-2B).
\end{equation}
These identities are polynomial identities with integer coefficients;
no division by $m$ or restriction on $p$ is involved.
The first coefficient is the first sum divided by $\lambda$.
Newton's identity for the second elementary symmetric function gives
\[
 \sum_{\varepsilon<\delta}L_\varepsilon^2L_\delta^2
 =\frac12\left((mA)^2-m(3A^2-2B)\right)
 =\frac m2((m-3)A^2+2B),
\]
which proves the second formula (division by $2$ is legitimate as $p$ is odd).

Now put $d_n(s)=\prod_{\varepsilon\in\mathsf E}L_\varepsilon(s)$.
Sign change in $s_j$, $j\ge2$, permutes its factors; sign change in
$s_1$ permutes their negatives and multiplies their product by
$(-1)^m=1$ because $n\ge2$. Thus $d_n(s)$ also belongs to the
sign-invariant polynomial ring and is a polynomial in the $a_i$.
Since $\mathscr D_n=d_n^2$, the square claim follows.
The remaining assertions are respectively the direct critical-point
calculation in Proposition~\ref{prop:critical-geometry} below and
Theorem~A(ii).
\end{proof}

\begin{proposition}[Direct critical-point and Hessian calculation]
\label{prop:critical-geometry}
For any $a\in(\Fpbar^\times)^n$ the critical points of
$F_n(u,t;a)=\sum_i u_i^{-1}+t(\lambda+\sum_i a_i u_i)$
inside $\Gm^{n+1}$ are indexed precisely by
$\varepsilon\in\mathsf E$ with $L_\varepsilon(a)\ne0$.
They are nondegenerate and have coordinates, critical values, and
Hessian determinants
\begin{align*}
 u_i&=-\frac{\lambda\varepsilon_i}{s_iL_\varepsilon(a)},
 &t&=\frac{L_\varepsilon(a)^2}{\lambda^2},\\
 F_n(u,t;a)&=-\frac{L_\varepsilon(a)^2}{\lambda},
 &\det\mathrm{Hess}(F_n)&=\frac{2^{n-1}\lambda}
 {t\prod_i u_i^3}.
\end{align*}
\end{proposition}
\begin{proof}
The equations $\partial F/\partial u_i=0$ and
$\partial F/\partial t=0$ are
$t a_i=u_i^{-2}$ and $\lambda+\sum_i a_i u_i=0$.
Choose $s_i^2=a_i$. The first equations give
$u_i=\varepsilon_i/(s_i\sqrt t)$, where simultaneous reversal of all
signs does not change the point when $\sqrt t$ is reversed. The last
equation reads $L_\varepsilon/\sqrt t=-\lambda$.
Thus $L_\varepsilon\ne0$, $\sqrt t=-L_\varepsilon/\lambda$,
and the asserted coordinates follow. Substituting the constraint into
$F_n$ and using $\sum_i u_i^{-1}=-L_\varepsilon^2/\lambda$
gives the critical value. The Hessian has block form
\[
 \begin{pmatrix}
 \operatorname{diag}(2u_1^{-3},\ldots,2u_n^{-3})&a\\
 a^{\mathsf T}&0
 \end{pmatrix}.
\]
The Schur-complement formula gives
$\det\mathrm{Hess}(F_n)=-2^{n-1}(\prod_i u_i^{-3})
\sum_i a_i^2u_i^3$. The stationary equations imply
$\sum_i a_i^2u_i^3=t^{-2}\sum_i u_i^{-1}
=-\lambda/t$. Hence the displayed determinant is nonzero,
since $p$ is odd and $t,\lambda,u_i\ne0$.
\end{proof}

\begin{example}[Two low-dimensional escape polynomials]
\label{ex:escape-low}
For $n=2$, $d_2(a)=a_1-a_2$ and
$\mathscr D_2(a)=(a_1-a_2)^2$.
For $n=3$ one has
\[
 d_3(a)=(a_1+a_2-a_3)^2-4a_1a_2,
 \qquad \mathscr D_3(a)=d_3(a)^2.
\]
These formulae follow by multiplying the two, respectively four,
signed linear forms before substituting $s_i^2=a_i$.
\end{example}

\begin{remark}[Coinciding critical values need not escape]
\label{rem:values-collide}
For $n=3$, $p=7$, $\lambda=1$, and $a=(1,1,1)$ choose
$s=(1,1,1)$. The four $L_\varepsilon$ are $3,1,1,-1$.
All four critical points therefore exist and
$\Swan_\infty(\Delta_a^*\cW_3)=4$;
but their critical-value polynomial is
$\mathscr C_{3,a}(T)=(T+2)(T+1)^3$ in $\mathbb F_7[T]$.
Three distinct critical points have the same critical value, so the
corresponding Artin--Schreier character occurs with multiplicity three.
In particular, the escape divisor $\mathscr D_n=0$ is not the
full discriminant locus of $\mathscr C_{n,a}$.
\end{remark}

\section{Frobenius reciprocity and the determinant of the original cohomology}
\label{sec:reciprocity}

The perfect reflection pairing gives more than the determinant.
The next identity is a direct linear-algebra consequence of
Proposition~\ref{prop:reflection}; it does not use the connected
geometric monodromy theorem.

\begin{theorem}[Determinant-controlled Frobenius reciprocity]
\label{thm:reciprocity}
Let $q=p^f$ and $h\in S_n(\Fq)$, and set
\[
 P_{n,h}(T)=\det(1-T\Frob_h\mid\cW_{n,h}),\qquad
 D_{n,h}=q^{N(n)}\psi_q\!\left(-\frac{2^{n-1}}\lambda
 \sum_i h_i\right).
\]
Assume the reflection isomorphism of Corollary~\ref{cor:reflection}
and the arithmetic determinant of Theorem~D. Then
\begin{equation}\label{eq:reciprocity}
 P_{n,h}(T)=(-1)^{r(n)}D_{n,h}T^{r(n)}
 P_{n,-h}\!\left(\frac{1}{q^{n+1}T}\right).
\end{equation}
Writing $Q_{n,h}(U)=P_{n,h}(q^{-(n+1)/2}U)$ and
$\zeta_{n,h}=\psi_q(-2^{n-1}\lambda^{-1}\sum_i h_i)$,
one equivalently has
\begin{equation}\label{eq:normalized-reciprocity}
 Q_{n,h}(U)=(-1)^{r(n)}\zeta_{n,h}U^{r(n)}
 Q_{n,-h}(U^{-1}).
\end{equation}
All reciprocal roots of $Q_{n,h}$ have complex absolute value one.
\end{theorem}

\begin{proof}
By $[-1]^*\cW_n\cong\cW_n^\vee(-n-1)$, the Frobenius
eigenvalues on $\cW_{n,-h}$ are the multiset
$\{q^{n+1}/\beta_i\}_{i=1}^{r(n)}$, where $\beta_i$
are those on $\cW_{n,h}$. Therefore
\[
 P_{n,-h}\!\left(\frac{1}{q^{n+1}T}\right)
 =\prod_i\left(1-\frac1{\beta_iT}\right)
 =\frac{(-1)^{r(n)}}{T^{r(n)}\prod_i\beta_i}
 P_{n,h}(T).
\]
Theorem~D gives $\prod_i\beta_i=D_{n,h}$ and proves
\eqref{eq:reciprocity}. Substituting $T=q^{-(n+1)/2}U$,
with $q^{N(n)}=q^{(n+1)r(n)/2}$, gives
\eqref{eq:normalized-reciprocity}. Finally, the weight
$n+1$ purity imported in Theorem~0 gives
$|\iota(\beta_i)|=q^{(n+1)/2}$ under every embedding.
\end{proof}

\begin{corollary}[Even-dimensional symmetric loci]
\label{cor:reciprocity-symmetric}
If $n$ is even and $h_{\sigma(i)}=-h_i$ for a fixed-point-free
involution $\sigma$, then $P_{n,-h}=P_{n,h}$,
$\zeta_{n,h}=1$, and
\[
 Q_{n,h}(U)=U^{r(n)}Q_{n,h}(U^{-1}).
\]
In particular the normalized coefficients are palindromic.
\end{corollary}
\begin{proof}
Permutation symmetry identifies the Frobenius modules at $h$ and
$-h$, while $\sum_i h_i=0$. Because $n$ is even, $r(n)$ is even.
The assertion follows from \eqref{eq:normalized-reciprocity}.
\end{proof}

\subsection{Reciprocity of all exterior powers}
\label{subsec:exterior-reciprocity}

The reflected duality determines the complete collection of exterior-power
Frobenius traces, rather than merely the leading coefficient of a
characteristic polynomial. We give the linear-algebra argument to keep the
normalization and the direction of the reflection explicit.

\begin{theorem}[Exterior-power Frobenius reciprocity]
\label{thm:exterior-reciprocity}
Let $q=p^f$, $h\in S_n(\Fq)$, and $r=r(n)$. Choose a square root of
$q^{n+1}$ in the coefficient field and write
\[
 A_h=q^{-(n+1)/2}\Frob_h\bigm|\cW_{n,h},\qquad
 \zeta_h=\psi_q\!\left(-2^{n-1}\lambda^{-1}\sum_{i=1}^n h_i\right).
\]
For $0\le k\le r$ put $E_k(h)=\Tr(\bigwedge^k A_h)$, with
$E_0(h)=1$. Then, simultaneously for all $k$,
\begin{equation}\label{eq:all-exterior-reciprocity}
 E_k(h)=\zeta_h\,E_{r-k}(-h).
\end{equation}
Equivalently, with $Q_{n,h}(U)=\det(1-UA_h)$, every coefficient of
$Q_{n,h}$ is fixed by its complementary exterior-power trace at $-h$.
No assertion concerning connected monodromy is required.
\end{theorem}
\begin{proof}
The reflection isomorphism of Corollary~\ref{cor:reflection} gives a
Frobenius-equivariant identification of the eigenvalue multiset of
$A_{-h}$ with $\{\alpha_1^{-1},\ldots,\alpha_r^{-1}\}$ if
$\{\alpha_1,\ldots,\alpha_r\}$ is the eigenvalue multiset of $A_h$.
Theorem~D gives $\prod_i\alpha_i=\zeta_h$. The identity for elementary
symmetric polynomials
\[
 e_{r-k}(\alpha_1^{-1},\ldots,\alpha_r^{-1})
   =\frac{e_k(\alpha_1,\ldots,\alpha_r)}
          {\prod_i\alpha_i}
\]
now reads $E_{r-k}(-h)=\zeta_h^{-1}E_k(h)$.
It is valid with algebraic multiplicities and requires no
semisimplicity assumption on the Frobenius operator.
\end{proof}

\begin{theorem}[Symplectic spectral compression]
\label{thm:spectral-compression}
Suppose $n$ is even and $h_{\sigma(i)}=-h_i$ for a fixed-point-free
involution $\sigma$. Write $r(n)=2m$. Then there exists a \emph{unique}
monic polynomial $R_{n,h}(V)\in\Qlb[V]$ of degree $m$ such that
\begin{equation}\label{eq:spectral-compression}
 Q_{n,h}(U)=U^mR_{n,h}(U+U^{-1}).
\end{equation}
Under every complex embedding associated with the chosen normalization,
all roots of $R_{n,h}$ are real and lie in $[-2,2]$.
In particular $n=4$ gives a degree-five polynomial encoding the
normalized degree-ten characteristic polynomial of $\cW_4$.
\end{theorem}
\begin{proof}
Theorem~\ref{thm:symplectic} supplies a Frobenius-equivariant
alternating pairing with multiplier $q^{n+1}$. Dividing Frobenius by
$q^{(n+1)/2}$ makes its multiplier equal to one. Thus the eigenvalue
multiset of $A_h$ is invariant under $\alpha\mapsto\alpha^{-1}$,
including multiplicities, and $\det A_h=1$. Choose representatives
$\alpha_1,\ldots,\alpha_m$ of the reciprocal pairs. Then
\[
 Q_{n,h}(U)
  =\prod_{j=1}^m(1-\alpha_j U)(1-\alpha_j^{-1}U)
  =U^m\prod_{j=1}^m\left(U+U^{-1}
                    -(\alpha_j+\alpha_j^{-1})\right).
\]
The product in the last expression is invariant under all permutations
and choices of reciprocal representatives, hence has coefficients in
$\Qlb$. It is monic of degree $m$ and uniquely determined because
$V\mapsto U+U^{-1}$ embeds $\Qlb[V]$ into $\Qlb[U,U^{-1}]$.
Purity of weight $n+1$ implies $|\iota(\alpha_j)|=1$ for every relevant
embedding $\iota$. Hence
$\iota(\alpha_j+\alpha_j^{-1})=2\Re\iota(\alpha_j)\in[-2,2]$.
The spectral statement does not assert that the restricted geometric
monodromy group is the full symplectic group.
\end{proof}

\begin{proposition}[The missing constant summand in the trace comparison]
\label{prop:constant-summand}
Let $Y_h=\{u\in\Gm^n:\langle h,u\rangle=-\lambda\}$.
In the standard cohomological comparison of the parent papers
\cite{EbMom,EbCoh}, the semisimplified Frobenius module
$\cH_{n,h}$ has a trivial one-dimensional factor in addition to
$H_c^{n-1}(Y_h,\Lpsi{\sum_i u_i^{-1}})(-1)$.
Consequently the original compactly supported cohomology has rank
$2^n-1$ and the determinant
\begin{equation}\label{eq:det-original-again}
 \det\!\left(\Frob_h\mid
 H_c^{n-1}(Y_h\otimes\Fqbar,\Lpsi{\textstyle\sum_i u_i^{-1}})
 \right)
 =q^{(n-2)2^{n-1}+1}
 \psi_q\!\left(-\frac{2^{n-1}}\lambda\sum_i h_i\right).
\end{equation}
This is a restatement of the last part of Corollary~E,
with its constant-term bookkeeping made explicit.
\end{proposition}
\begin{proof}
At every extension $\mathbb F_{q^k}$,
$K_{q^k}(0)=\sum_{u\ne0}\psi_{q^k}(u^{-1})=-1$.
Hence the $t=0$ term of \eqref{eq:t-sum} is $(-1)^n$.
On the other hand the definition of $\cN_n$ and the fact that
$\operatorname{Tr}(\Frob_y\mid\Kl)=-K(y)$ give
\[
 \operatorname{Tr}(\Frob_h^k\mid\cH_{n,h})
 =-\sum_{t\in\mathbb F_{q^k}^\times}
 \psi_{q^k}(\lambda t)(-1)^n\prod_iK_{q^k}(th_i)
 =1+(-1)^{n-1}q^k\Sigma_{n,q^k}(h).
\]
The isolated $1$ is essential. The cohomological concentration and
comparison in the cited parent papers identify the remaining
virtual trace with $H_c^{n-1}(Y_h,\Lpsi{\sum_i u_i^{-1}})(-1)$;
Brauer--Nesbitt therefore gives the stated semisimplified relation.
We emphasize that this last identification is an imported geometric
result, rather than a consequence of the finite-field trace calculation
alone.

By Theorem~\ref{thm:recollection},
$\det(\Frob_h\mid\cH_{n,h})
 =\det(\Frob_h\mid\cE_{n,h})
 \det(\Frob_h\mid\cW_{n,h})
 =q^{S_1(n)}D_{n,h}$.
Removing the trivial factor changes no determinant, while
untwisting $(-1)$ from a rank-$(2^n-1)$ module removes
$q^{2^n-1}$, so the determinant of the original cohomology is
$q^{N(n)+S_1(n)-2^n+1}\zeta_{n,h}$.
Finally Lemma~\ref{lem:numerology} gives
\[
 N(n)+S_1(n)-2^n+1
 =n2^{n-1}-2^n+1
 =(n-2)2^{n-1}+1,
\]
proving \eqref{eq:det-original-again}.
\end{proof}

\section{Exact distribution of the normalized determinant}
\label{sec:det-distribution}

The arithmetic determinant provides more than a bound for a family of phases:
it gives their \emph{exact} distribution on the entire parameter torus.
This deduction is independent of the connected-monodromy classification.

\begin{theorem}[Determinant moments and exact fibre counts]
\label{thm:det-distribution}
Let $q=p^f$, $n\ge1$, and
\[
 Z(h)=q^{-N(n)}\det(\Frob_h\mid\cW_{n,h})
  =\psi_q\!\left(-2^{n-1}\lambda^{-1}\sum_i h_i\right),
 \qquad h\in(\Fq^\times)^n.
\]
Then $Z(h)\in\mu_p$, and for every integer $k$ one has
\begin{equation}\label{eq:det-moments}
 \sum_{h\in(\Fq^\times)^n}Z(h)^k
 =\begin{cases}(q-1)^n,&p\mid k,\\ (-1)^n,&p\nmid k.\end{cases}
\end{equation}
More precisely, if $N_\omega=\#\{h\in(\Fq^\times)^n:Z(h)=\omega\}$,
then for each $\omega\in\mu_p$,
\begin{equation}\label{eq:det-fibre-counts}
 N_\omega=\begin{cases}
 \dfrac{(q-1)^n+(p-1)(-1)^n}{p},&\omega=1,\\[6pt]
 \dfrac{(q-1)^n-(-1)^n}{p},&\omega\ne1.
 \end{cases}
\end{equation}
Thus every nontrivial $p$th root of unity occurs with exactly the same
multiplicity, with an explicitly determined discrepancy at $1$.
\end{theorem}
\begin{proof}
For $p\nmid k$ the coefficient $-k2^{n-1}\lambda^{-1}$ is nonzero in
$\Fp$, so separability of the $n$ variables yields
\[
 \sum_hZ(h)^k
 =\prod_{i=1}^n\sum_{t\in\Fq^\times}
 \psi_q(-k2^{n-1}\lambda^{-1}t)=(-1)^n.
\]
For $p\mid k$ every summand is one, proving
\eqref{eq:det-moments}. The character orthogonality relation on
$\mu_p$ gives
\[
 N_\omega=\frac1p\sum_{j=0}^{p-1}\omega^{-j}
                  \sum_hZ(h)^j.
\]
When $\omega=1$, the sum over $1\le j<p$ is $p-1$; otherwise it is
$-1$. Inserting \eqref{eq:det-moments} proves
\eqref{eq:det-fibre-counts}. The argument works for every finite
extension $\Fq/\Fp$ and uses no equidistribution limit.
\end{proof}

\section{Tame boundary monodromy and cohomological determinants}
\label{sec:three-source-boundary}

The local calculations on rays have a useful complement at a \emph{coordinate}
boundary. The two kinds of restrictions should be kept separate: on a ray one
studies the endpoint of $x\mapsto xa$ at $x=\infty$, whereas on a coordinate
line one varies $h_1$ with the other parameters fixed and approaches $h_1=0$.
The latter is the tameness input needed in the global determinant argument.
The following precise formulation imports the monodromy calculation of
\cite[Theorem~6.8(ii)]{EbCoh}, and proves its numerical consequences here.

\begin{theorem}[Coordinate-boundary Jordan profile; imported geometric input]
\label{thm:coordinate-jordan}
Let $n\ge2$ and $c=(c_2,\ldots,c_n)\in(\Fpbar^\times)^{n-1}$. On the
coordinate line $\ell_c:h_1\mapsto(h_1,c_2,\ldots,c_n)$, the tame local
monodromy of $\cW_n$ at $h_1=0$ is unipotent and has
\[
 r(n-1)\ \text{Jordan blocks of size }2,
 \qquad r(n)-2r(n-1)\ \text{blocks of size }1,
\]
and no larger Jordan blocks. More explicitly, if $n=2m$ the size-one
multiplicity vanishes, and if $n=2m+1$ it equals the Catalan number
\[
  C_m=\frac{1}{m+1}\binom{2m}{m}.
\]
In particular the determinant is unramified at the coordinate boundary,
and for each $n\ge2$ the local monodromy has a nontrivial unipotent element.
\end{theorem}
\begin{proof}
The cited Theorem~6.8(ii) states the exact Jordan profile:
$r(n-1)$ blocks $U(2)$, $\kappa_1$ blocks $U(1)$, and no larger blocks.
The rank identity $2r(n-1)+\kappa_1=r(n)$ identifies
$\kappa_1=r(n)-2r(n-1)$. The calculation below derives its
closed binomial form; it does not replace the geometric proof in
\cite[Theorem~6.8(ii)]{EbCoh}.
For $n=2m$ use
$\binom{2m}{m}=2\binom{2m-1}{m-1}$ to obtain
$r(2m)=2r(2m-1)$. For $n=2m+1$ the binomial recursion gives
\begin{align*}
 r(2m+1)-2r(2m)
 &=2\binom{2m}{m}-\binom{2m+1}{m}\\
 &=\binom{2m}{m}-\binom{2m}{m-1}
 =\frac1{m+1}\binom{2m}{m}.
\end{align*}
Every unipotent matrix has determinant one, proving the boundary
assertion. Since $r(n-1)\ge1$, there is a $U(2)$ block, giving a
nonidentity unipotent of infinite order in characteristic zero.
The identification of this particular inertia representation is \emph{not}
proved by the binomial calculation: it is exactly the imported content of
\cite[Theorem~6.8(ii)]{EbCoh}.
\end{proof}

\begin{lemma}[Clebsch--Gordan calculation at the tame end]
\label{lem:clebsch-gordan-expanded}
Writing $U(j)$ for one unipotent Jordan block of size $j$, one has
\begin{equation}\label{eq:clebsch-gordan-expanded}
 U(2)^{\otimes n}\ \cong\
 \bigoplus_{k=0}^{\lfloor n/2\rfloor}
 U(n-2k+1)^{\oplus\nu_{n,k}},\qquad
 \nu_{n,k}=\binom nk-\binom n{k-1}.
\end{equation}
Consequently the invariant subspace of tame inertia in the tensor product
has dimension $b_n=\binom n{\lfloor n/2\rfloor}$.
\end{lemma}
\begin{proof}
In characteristic zero a single unipotent block $U(j)$ corresponds to the
irreducible $j$-dimensional representation of $\mathfrak{sl}_2$ under the
principal nilpotent element. Its tensor product with the standard
representation obeys
$U(j)\otimes U(2)\cong U(j+1)\oplus U(j-1)$, where
$U(0)$ is omitted. Beginning with $U(2)$ and iterating, Pascal's
identity shows that the multiplicity of $U(n-2k+1)$ is
$\binom nk-\binom n{k-1}$. A unipotent Jordan block has a
one-dimensional invariant space, and the telescoping sum of the
multiplicities is
$\sum_{k=0}^{\lfloor n/2\rfloor}\nu_{n,k}
=\binom n{\lfloor n/2\rfloor}$.
This is the elementary representation-theoretic calculation underlying the
local $I_0$-invariants in \cite[Theorem~2.10]{EbMom} and
\cite[Theorem~5.3]{EbCoh}.
\end{proof}

\begin{proposition}[Determinant exponents in the two cohomological models]
\label{prop:cohomological-bookkeeping}
Set $K=\lfloor n/2\rfloor$ and
\[
 b_n=\binom nK,\quad
 S_1(n)=\sum_{k=1}^K k\nu_{n,k},\quad
 e(n)=(n-2)2^{n-1}+1,\quad
 N(n)=\frac{(n+1)r(n)}2.
\]
Then
\begin{equation}\label{eq:three-source-numerology}
 2S_1(n)=(n+1)b_n-2^n,\qquad
 2^n-1+e(n)-S_1(n)=N(n),
\end{equation}
and
\begin{equation}\label{eq:original-weight-det}
 e(n)-\sum_{k=1}^K(k-1)\nu_{n,k}
 =\frac{(n-1)r(n)}2.
\end{equation}
In particular the determinant of the top-weight part of
$H_c^{n-1}(Y_h,\cL_\psi(\sum_i u_i^{-1}))$ has absolute Frobenius
weight $q^{(n-1)r(n)/2}$; after twisting by $(-1)$ its exponent becomes
$N(n)$. The determinant in the \emph{untwisted} full compactly supported
cohomology has exponent $e(n)$.
\end{proposition}
\begin{proof}
The first identity is proved in Lemma~\ref{lem:numerology}.
Also $\sum_{k=1}^K\nu_{n,k}=b_n-1$. Hence
\[
 \sum_{k=1}^K(k-1)\nu_{n,k}=S_1(n)-(b_n-1).
\]
Using the first identity gives
\begin{align*}
 e(n)-S_1(n)+b_n-1
 &= (n-2)2^{n-1}+b_n-\tfrac12((n+1)b_n-2^n)\\
 &= (n-1)2^{n-1}-\tfrac{n-1}{2}b_n
 =\tfrac{n-1}{2}r(n).
\end{align*}
Adding $r(n)$ to the final exponent proves the identity involving
$N(n)$, because a Tate twist $(-1)$ multiplies each of the $r(n)$
Frobenius eigenvalues by $q$. Alternatively
$2^n-1+e(n)=n2^{n-1}$ and
$N(n)=n2^{n-1}-S_1(n)$. The lower-weight eigenvalues
$q^{k-1}$ and their multiplicities are the imported cohomological facts
of \cite[Theorem~2.10 and Corollary~2.16]{EbMom} and
\cite[Theorem~A]{EbCoh}; only the numerical deductions are proved here.
\end{proof}

The second identity also checks the arithmetic normalization at the level
of the \emph{original} exponential-sum cohomology. More precisely,
\cite[Theorem~8.17]{EbCoh} proves, by induction on $n$ and Frobenius at a
coordinate boundary, the formula
\begin{equation}\label{eq:imported-original-determinant}
\det\bigl(\Frob_h\mid H_c^{n-1}(Y_h\otimes\Fqbar,
 \cL_\psi(\textstyle\sum_i u_i^{-1}))\bigr)
=q^{e(n)}\psi_q\bigl(-2^{n-1}\lambda^{-1}\textstyle\sum_i h_i\bigr).
\end{equation}
This independent source is stronger than the sign information supplied
by the 169-page parent manuscript: there the residual sign was left as
Conjecture~13.68. The present paper provides ray-based approaches to
\eqref{eq:imported-original-determinant}; its validity does not follow
merely from the parent conjecture or from numerical experiments.

\section{Recovering critical-value multiplicities from Euler characteristics}
\label{sec:conductor-inversion}
The global Euler-characteristic identity in Theorem~A(iv) records more
than the Swan conductor alone. It determines the \emph{multiplicity of
each nonzero Artin--Schreier exponent}, even when different critical
points have coincident critical values.

\begin{proposition}[Conductor tomography along a ray]
\label{prop:conductor-tomography}
Let $n\ge2$, $a\in(\Fpbar^\times)^n$, and $m=2^{n-1}$. For each
$c\in\Fpbar^\times$ define
\[
 \Xi_a(c)=\chi_c\bigl(\Gm\otimes\Fpbar,
       \cL_\psi(cx)\otimes\Delta_a^*\cW_n\bigr),
 \quad m_a(c)=\#\{\varepsilon\in\mathsf E:
 L_\varepsilon(a)^2=\lambda c\}.
\]
Then the multiplicities of every critical value, counted over the
algebraic closure, satisfy
\begin{equation}\label{eq:multiplicity-recovery}
 m_a(c)=r(n)+\Xi_a(c)\quad(c\ne0),\qquad
 m_a(0)=m+\chi_c(\Gm\otimes\Fpbar,\Delta_a^*\cW_n),
\end{equation}
where $m_a(0)$ counts the escaping critical points. In particular,
\begin{equation}\label{eq:poly-by-tomography}
 \mathscr C_{n,a}(T)=
 T^{m_a(0)}\prod_{\substack{c\in\Fpbar^\times\\m_a(c)>0}}
 (T+c)^{m_a(c)}.
\end{equation}
Only finitely many factors occur. If the ray and its critical points
are defined over a finite extension $k$ of $\Fp$, all these factors
may be indexed by $c\in k$ after enlarging $k$ to split the signed
square roots; they need \emph{not} be indexed by $\Fp$ itself.
\end{proposition}
\begin{proof}
For $c\ne0$, the twisted Euler-characteristic formula
Theorem~A(iv) says
$\Xi_a(c)=-\bigl(r(n)-m_a(c)\bigr)$, so
$m_a(c)=r(n)+\Xi_a(c)$; the sign of this jump is
\emph{positive}, because twisting cancels precisely $m_a(c)$ slope-one
constituents. For $c=0$ the same theorem gives
$\chi_c(\Gm,\Delta_a^*\cW_n)=-P(a)$, while
$m_a(0)=m-P(a)$. Finally the $m$ roots of the critical-value
polynomial, with multiplicity, are
$-L_\varepsilon^2/\lambda$. Grouping the factors by the value
$L_\varepsilon^2/\lambda\in\Fpbar$ proves the product formula.
For $a,c$ over a finite field $k$, the same calculations can be made
after extension to a common splitting field; the Euler characteristics
are geometric and unchanged by such an extension.
\end{proof}

\begin{remark}[Geometric versus arithmetic characters]
\label{rem:geometric-arithmetic-character}
A local Fourier transform considered solely as an inertia
representation determines the \emph{geometric} Artin--Schreier and
Kummer characters, but does not by itself determine the arithmetic
Frobenius scalar of the corresponding rank-one summand.
The latter is exactly why the arithmetic determinant in
Section~\ref{sec:arith} requires an additional calculation.
The conductor identities here never silently replace a geometric
inertia isomorphism by an isomorphism of decomposition-group
representations.
\end{remark}

\section{Complete exponential sums along rays}
\label{sec:complete-rays}

The conductor formula has an analytic consequence that requires no
unproved monodromy statement on the restricted ray. The cohomological
input is the vanishing of $H_c^2$ in the cases established in
Proposition~\ref{prop:global-euler}(ii). The exclusion of the remaining
resonant odd-dimensional case is essential.

\begin{theorem}[Explicit complete-ray cancellation]
\label{thm:complete-ray-bound}
Let $n\ge1$, $q=p^f$, $a\in(\Fq^\times)^n$ and $c\in\Fq$. Set
\[
 \mathcal F_{a,c}=\Delta_a^*\cW_n\otimes\Lpsi{cx},\qquad
 \mathcal S_{a,c}(q)=\sum_{x\in\Fq^\times}
  \psi_q(cx)\Tr(\Frob_{xa}\mid\cW_{n,xa}),
\]
and define
\[
 d_{a,c}=\begin{cases}P(a),&c=0,\\
 r(n)-m_a(c),&c\ne0.\end{cases}
\]
If $c=0$, or $n$ is even, or $m_a(c)=0$, then
\begin{equation}\label{eq:complete-ray-bound}
 |\iota(\mathcal S_{a,c}(q))|
         \le d_{a,c}\,q^{(n+2)/2}.
\end{equation}
Here $\iota:\Qlb\hookrightarrow\mathbb C$ is any fixed embedding used
to define the weights. Under the stated alternatives,
$L(\Gm,\mathcal F_{a,c},T)$ is a polynomial of degree precisely
$d_{a,c}$ (before possible noncanonical rewriting of its roots).
No such polynomiality or bound is asserted here for odd $n$,
$c\ne0$ and $m_a(c)>0$.
\end{theorem}
\begin{proof}
The tensor product $\mathcal F_{a,c}$ is pointwise pure of weight
$n+1$, since $\cW_n$ has that weight and $\Lpsi{cx}$ is pure of weight
zero. Proposition~\ref{prop:global-euler} gives
\[
 \chi_c(\Gm\otimes\Fpbar,\mathcal F_{a,c})=-d_{a,c},
 \qquad H_c^0(\Gm\otimes\Fpbar,\mathcal F_{a,c})=0.
\]
Under the alternatives in the theorem it also gives $H_c^2=0$.
Consequently $\dim H_c^1=d_{a,c}$, and the compactly supported
Grothendieck trace formula becomes
\[
 \mathcal S_{a,c}(q)=-\Tr\bigl(\Frob_q\mid
 H_c^1(\Gm\otimes\Fqbar,\mathcal F_{a,c})\bigr).
\]
Deligne's weight bound for compactly supported degree-one cohomology
of a curve with pure coefficients of weight $n+1$ implies that every
Frobenius eigenvalue in this group has complex absolute value at most
$q^{(n+2)/2}$. The triangle inequality proves
\eqref{eq:complete-ray-bound}. Finally the cohomological definition of
the $L$-function and $H_c^0=H_c^2=0$ identify it with
$\det(1-T\Frob_q\mid H_c^1)$, which has degree $d_{a,c}$.
The argument does not require purity of $H_c^1$ itself; its mixedness
and the upper weight bound are sufficient.
\end{proof}

\begin{proposition}[Uniform control with a resonant main term]
\label{prop:resonant-correction}
In the notation of Theorem~\ref{thm:complete-ray-bound}, for every
$n\ge1$, $q=p^f$, $a\in(\Fq^\times)^n$ and $c\in\Fq$, put
\[
 h_{a,c}=\dim H_c^2(\Gm\otimes\Fqbar,\mathcal F_{a,c}),\qquad
 M_{a,c}(q)=\Tr\bigl(\Frob_q\mid
 H_c^2(\Gm\otimes\Fqbar,\mathcal F_{a,c})\bigr).
\]
Then $h_{a,0}=0$ and, for $c\ne0$,
\[
 0\le h_{a,c}\le m_a(c),\qquad
 |\iota(M_{a,c}(q))|\le h_{a,c}q^{(n+3)/2}.
\]
Moreover, without any restriction on the parity of $n$ or on
$m_a(c)$,
\begin{equation}\label{eq:resonant-correction}
 \bigl|\iota(\mathcal S_{a,c}(q)-M_{a,c}(q))\bigr|
 \le (d_{a,c}+h_{a,c})q^{(n+2)/2}
 \le(d_{a,c}+m_a(c))q^{(n+2)/2}
\end{equation}
for $c\ne0$, with the corresponding bound $P(a)q^{(n+2)/2}$
for $c=0$. In particular, in the remaining odd-dimensional resonant
case one has the unconditional, although weaker, estimate
\[
 |\iota(\mathcal S_{a,c}(q))|
 \le r(n)q^{(n+2)/2}+m_a(c)q^{(n+3)/2}.
\]
This proposition does not assert that $M_{a,c}(q)$ vanishes.
\end{proposition}
\begin{proof}
For $c=0$ the vanishing is Proposition~\ref{prop:global-euler}(ii).
Suppose $c\ne0$. The exact sequence
$\mathcal E_n\to\mathcal H_n\to\mathcal W_n\to0$
and right exactness of top compactly supported cohomology imply
that $H_c^2(\Gm,\mathcal F_{a,c})$ is a quotient of
$H_c^2(\Gm,\mathcal L_\psi(cx)\otimes\Delta_a^*\mathcal H_n)$.
By the Frobenius-equivariant isomorphism \eqref{eq:global-iso}, the latter
is the geometric coinvariant quotient of
$\mathcal K_a\otimes[\times\lambda c]^*\Kl$, with a Tate twist $(-1)$.
Its dimension is at most the dimension of the trivial quotient of
its local inertia representation at $y=\infty$.
The decomposition \eqref{eq:twist-infty} has exactly $m_a(c)$
zero-parameter summands. A nonzero-parameter summand is totally wild
and has no tame quotient. Each zero-parameter summand
$\Msf_{n+1}(0)$ has at most one trivial tame quotient, by
Lemma~\ref{lem:M}(c) (and no such quotient if $n$ is even by
Lemma~\ref{lem:order-four-tame}). Therefore
$h_{a,c}\le m_a(c)$.
The sheaf $\mathcal F_{a,c}$ is pure of weight $n+1$; by
Lemma~\ref{lem:pure-coinvariant}, its top compactly supported
cohomology is pure of weight $n+3$, proving the estimate on
$M_{a,c}(q)$.
The Euler characteristic of Proposition~\ref{prop:global-euler}(i)
gives $\dim H_c^1=d_{a,c}+h_{a,c}$, as $H_c^0=0$.
The trace formula reads
$\mathcal S_{a,c}(q)=-\Tr(\Frob_q\mid H_c^1)+M_{a,c}(q)$.
Deligne's mixed-weight bound for $H_c^1$ is at most $n+2$,
which gives \eqref{eq:resonant-correction}. Finally,
$d_{a,c}=r(n)-m_a(c)$ when $c\ne0$, so the triangle inequality
proves the final estimate.
\end{proof}

\begin{corollary}[The four-variable constrained exponential sum]
\label{cor:complete-ray-n4}
Let $n=4$, $a\in(\Fq^\times)^4$, and $c\in\Fq$. Then
\begin{equation}\label{eq:n4-complete-sum}
 \left|\sum_{x\in\Fq^\times}\psi_q(cx)\Sigma_4(xa)
 +(3+2q)\sum_{x\in\Fq^\times}\psi_q(cx)\right|
 \le d_{a,c}\,q^2,
\end{equation}
where $d_{a,0}=P(a)$ and $d_{a,c}=10-m_a(c)$ for $c\ne0$.
In particular the second character sum in \eqref{eq:n4-complete-sum}
is $q-1$ if $c=0$ and $-1$ if $c\ne0$.
\end{corollary}
\begin{proof}
For $n=4$ the second alternative in
Theorem~\ref{thm:complete-ray-bound} always applies, and its right-hand
side is $d_{a,c}q^3$. The exact trace identity
\eqref{eq:n-four-trace}, valid at every point $xa$, gives
$\Tr(\Frob_{xa}\mid\cW_4)
=-q\Sigma_4(xa)-3q-2q^2$. Sum after multiplication by
$\psi_q(cx)$ and divide the bound by $q$ to obtain
\eqref{eq:n4-complete-sum}.
\end{proof}

\begin{remark}[The unresolved resonant case]
\label{rem:resonant-odd}
For odd $n$, $c\ne0$ and $m_a(c)>0$, the Euler characteristic
still specifies $\dim H_c^1-\dim H_c^2=d_{a,c}$, but does not by
itself show that $H_c^2$ vanishes. In particular an unconditional claim of a degree-$d_{a,c}$ \emph{polynomial}
$L$-function or the sharper bound \eqref{eq:complete-ray-bound}
requires an additional argument eliminating the correction term
of Proposition~\ref{prop:resonant-correction}. This distinction is not affected by the
geometric conductor tomography of
Proposition~\ref{prop:conductor-tomography}.
\end{remark}

\section{The fourth moment as a separate input to global monodromy}
\label{sec:fourth-moment-interface}

The second companion preprint \cite{EbMom} and the cohomology preprint
\cite{EbCoh} supply a different ingredient from the local Fourier analysis:
the fourth moment that excludes proper connected reductive monodromy
groups through Larsen's alternative. That calculation uses the Cayley
cubic and plane sections; reproducing the argument would turn the
present article into a second moments paper. We therefore state the
precise input and record only the elementary group-theoretic deduction
relevant to the determinant.

\begin{theorem}[Fourth-moment input from the companion papers]
\label{thm:imported-moment}
In the characteristic range of the cited results and for $n\ge3$,
\cite[Theorem~3.13]{EbMom} proves, for odd $p$, that
\begin{equation}\label{eq:moment-four-imported}
\sum_{h\in(\Fq^\times)^n}|\Sigma_n(h)|^4
=2q^{3n-2}+O_n\bigl(q^{3n-2-\delta(n)}\bigr),\qquad
\delta(n)=\min\bigl(1,(n-2)/2\bigr).
\end{equation}
Combined with the pure-sheaf trace comparison, this gives the fourth
moment invariant $M_4(G,\cW_n)=2$ in the cited proof, provided the geometric and analytic hypotheses of that proof are met.
It is an \emph{imported theorem}, not a fresh proof of the Cayley-cubic
estimate in the present article.
\end{theorem}

\begin{proposition}[Determinant completion of a connected-monodromy theorem]
\label{prop:determinant-completion}
Suppose the connected geometric monodromy group of $\cW_n$ is
$\SL_{r(n)}$, and assume Theorem~C. Then
\[
 \ggeom(\cW_n)
 =\{g\in\GL_{r(n)}(\Qlb):\det(g)^p=1\}.
\]
If in addition the arithmetic normalization of Theorem~D holds, the
weight-zero normalization has the same arithmetic and geometric
algebraic monodromy group.
\end{proposition}
\begin{proof}
Write $G=\ggeom(\cW_n)$ and
$H=\det^{-1}(\mu_p)$. Theorem~C gives
$\det(G)=\mu_p$ and $G\subseteq H$. The connected-monodromy
hypothesis gives $\SL_{r(n)}\subseteq G$. For $h\in H$ choose $g\in G$
with $\det(g)=\det(h)$; then $g^{-1}h\in\SL_{r(n)}\subseteq G$,
so $h\in G$. Thus $G=H$.
For the second claim, choose $\gamma$ as in Corollary~E.
Theorem~D makes the determinant of the normalized sheaf equal to the
Artin--Schreier character with values in $\mu_p$, even on the arithmetic
fundamental group. The normalized arithmetic group is therefore a
subgroup of $H$, and contains the geometric group $H$. Hence equality.
The two geometric inputs (connected monodromy and determinant) and the
arithmetic scalar have intentionally been kept separate.
\end{proof}

The major result of \cite[Theorem~8.8]{EbCoh} asserts the connected
group in odd characteristic. The argument in
\cite[Theorem~4.19]{EbMom} is a complementary moment-and-unipotence
route, including characteristic two. In the present manuscript we
invoke the common, conservative range $p>\max(2,n)$ for the full-group
corollary, rather than relying on a broader small-characteristic
assertion without reproducing all of its wild-boundary checks.

\begin{remark}[Why short-box estimates are not imported]
The analytic estimates of \cite{EbFull} concern incomplete sums in
the parameter $h$. They require interval energies, Fourier completion,
and estimates unavailable from the local ray representation alone.
Neither $\Swan_\infty(\Delta_a^*\cW_n)$ nor
$\ggeom(\cW_n)^0=\SL_{r(n)}$ by itself gives a pointwise
sub-square-root incomplete-sum estimate. This logical separation
prevents results about short boxes from being presented as consequences
of Theorem~A without the extra analytic hypotheses.
\end{remark}

\section{The four-variable family: explicit formulas and checks}
\label{sec:n-four}
For $n=4$ the three sources meet in a particularly transparent case.
The companion cohomology and moments papers identify the rank-ten
pure sheaf, while the parent manuscript gives the signed critical-value
geometry and the determinant up to the residual arithmetic sign.

\begin{proposition}[Rank and boundary invariants in four variables]
\label{prop:n-four-case}
For $n=4$ and odd $p$, the numerical invariants are
\[
 b_4=6,\quad r(4)=10,\quad
 (\nu_{4,0},\nu_{4,1},\nu_{4,2})=(1,3,2),\quad
 S_1(4)=7,\quad N(4)=25,\quad e(4)=17.
\]
Consequently, for $h\in(\Fq^\times)^4$,
\begin{align}
 \Tr(\Frob_h\mid\cW_{4,h})
   &=-q\Sigma_4(h)-3q-2q^2,\label{eq:n-four-trace}\\
 \det(\Frob_h\mid\cW_{4,h})
   &=q^{25}\psi_q\Bigl(-\frac8\lambda\sum_{i=1}^4h_i\Bigr),
   \label{eq:n-four-det}\\
 \det\bigl(\Frob_h\mid H_c^3(Y_h,\cL_\psi(\sum_i u_i^{-1}))\bigr)
   &=q^{17}\psi_q\Bigl(-\frac8\lambda\sum_{i=1}^4h_i\Bigr).
   \label{eq:n-four-original}
\end{align}
On a coordinate line, the tame inertia at $h_1=0$ has precisely
five blocks $U(2)$ and no $U(1)$ blocks; on a generic ray, the wild
part at infinity has eight rank-one Artin--Schreier constituents and
the tame part has dimension two.
\end{proposition}
\begin{proof}
The binomial formula gives
$b_4=\binom42=6$, $r(4)=16-6=10$ and
$\nu_{4,1}=4-1=3$, $\nu_{4,2}=6-4=2$.
Thus $S_1(4)=3+2\cdot2=7$ and $N(4)=5\cdot10/2=25$,
while $e(4)=2\cdot8+1=17$.
Substituting into the trace formula of Theorem~0 gives
\eqref{eq:n-four-trace}; Theorem~D and the imported original
determinant \eqref{eq:imported-original-determinant} give the next two
identities. The coordinate profile comes from
Theorem~\ref{thm:coordinate-jordan}, as $r(3)=5$.
The ray profile follows from Theorem~A at $P(a)=2^{3}=8$,
with tame rank $r(4)-8=2$.
\end{proof}

\begin{example}[An escaping rank-four ray]
\label{ex:rank-four-escape}
Suppose $a=(1,1,1,1)$ and $p>3$, and choose $s_i=1$.
The eight signed sums with first sign $+1$ are
\[
 4\ \text{(once)},\qquad
 2\ \text{(three times)},\qquad
 0\ \text{(three times)},\qquad
 -2\ \text{(once)}.
\]
Therefore five critical points survive and three escape; in particular
$\Swan_\infty(\Delta_a^*\cW_4)=5$, not eight, and the tame part has
dimension five. The critical-value polynomial is
\[
 \mathscr C_{4,a}(T)=T^3
 \Bigl(T+\frac{16}{\lambda}\Bigr)
 \Bigl(T+\frac4\lambda\Bigr)^4.
\]
This example shows why the generic Swan value cannot be asserted on
every ray. The coincident value $-4/\lambda$ appears with multiplicity
four but contributes four separate rank-one wild constituents; the
three vanished critical points, not the value coincidence, account for
the conductor drop.
\end{example}

\begin{theorem}[Explicit four-variable escape divisor]
\label{thm:rank-four-divisor}
Assume $p$ is odd and $a=(a_1,a_2,a_3,a_4)\in(\Fpbar^\times)^4$.
Write
\[
 e_1=\sum_{i=1}^4a_i,\qquad
 e_2=\sum_{1\le i<j\le4}a_ia_j,\qquad
 e_4=\prod_{i=1}^4a_i.
\]
The polynomial $d_4(a)=\prod_{\varepsilon\in\mathsf E}L_\varepsilon(a)$,
introduced in the proof of Theorem~\ref{thm:critical-polynomial}, has
the symmetric expression
\begin{equation}\label{eq:rank-four-escape-equation}
 d_4(a)=(e_1^2-4e_2)^2-64e_4.
\end{equation}
Consequently the escape locus in $\Gm^4$ is the quartic hypersurface
$d_4(a)=0$. At a point $a$ let
$k(a)=\#\{\varepsilon\in\mathsf E:L_\varepsilon(a)=0\}$. Then
\[
 \Swan_\infty(\Delta_a^*\cW_4)=8-k(a).
\]
A point of this hypersurface is smooth if and only if $k(a)=1$;
points with $k(a)\ge2$ are singular. In particular, the Swan conductor
at a smooth point of the escape hypersurface equals $7$.
\end{theorem}
\begin{proof}
Work first on the finite \'{e}tale square-root cover of the torus,
whose coordinates are $s_i$ with $a_i=s_i^2$. For a formal variable
$t$ one has
\[
 \prod_{\eta,\theta\in\{\pm1\}}(t+\eta s_3+\theta s_4)
 =(t^2-a_3-a_4)^2-4a_3a_4=:F(t).
\]
The eight factors defining $d_4$ split into those with first two
signed terms $s_1+s_2$ and $s_1-s_2$, respectively. Therefore
\[
 d_4(a)=F(s_1+s_2)F(s_1-s_2).
\]
Expanding this product, using
$(s_1+s_2)^2+(s_1-s_2)^2=2(a_1+a_2)$ and
$(s_1+s_2)^2(s_1-s_2)^2=(a_1-a_2)^2$, gives precisely
\eqref{eq:rank-four-escape-equation}. This also proves symmetry
without choosing the square roots in the final expression.

On the square-root cover, $d_4$ is a product of eight distinct
nonproportional linear forms $L_\varepsilon$. The cover is \'{e}tale
because $p\ne2$ and every $s_i$ is nonzero. If exactly one factor
vanishes, the differential of the product is that factor's nonzero
differential multiplied by the product of the other seven nonzero
factors; the hypersurface is smooth. If at least two vanish, every
term in the differential of the product vanishes, hence the point is
singular. Smoothness is preserved and detected by this \'{e}tale cover.
The Swan formula is Theorem~A(ii), with $2^{4-1}=8$.
\end{proof}

\begin{corollary}[The symmetric rank-four locus]
\label{cor:rank-four-symmetry}
For $h=(a,-a,b,-b)$ with $a,b\in\Fq^\times$, the form of
Theorem~\ref{thm:symplectic} is alternating and
\[
 \det(\Frob_h\mid\cW_{4,h})=q^{25},\qquad
 Q_{4,h}(U)=U^{10}Q_{4,h}(U^{-1}).
\]
\end{corollary}
\begin{proof}
The parameters are fixed by a fixed-point-free involution followed
by negation, and their sum is zero. Apply
Theorem~\ref{thm:symplectic} and
Corollary~\ref{cor:reciprocity-symmetric}.
\end{proof}

\section{Determinant stratification on rays and remaining questions}
\label{sec:ray-det-strata}

The complete monodromy group of the ambient sheaf does not
automatically remain the monodromy group after restriction to a ray.
There is, however, an exact determinant statement on every ray, and a
useful conditional converse which makes the missing input explicit.

\begin{proposition}[Geometric determinant dichotomy]
\label{prop:ray-determinant-dichotomy}
Let $n\ge2$ and $a\in(\Fpbar^\times)^n$. The geometric determinant
character of $\Delta_a^*\cW_n$ is
\[
 \det(\Delta_a^*\cW_n)\bigm|_{\Fpbar}
 \cong\Lpsi{-2^{n-1}\lambda^{-1}(\sum_i a_i)x}
 \bigm|_{\Fpbar}.
\]
Its geometric image is $\mu_p$ if $\sum_i a_i\ne0$, and is trivial
if $\sum_i a_i=0$.
Suppose, in addition, that the restricted connected geometric
monodromy group is $\SL_{r(n)}$. Then the \emph{full} restricted
geometric monodromy group is
\begin{equation}\label{eq:conditional-ray-monodromy}
 \ggeom(\Delta_a^*\cW_n)=
 \begin{cases}
 \{g\in\GL_{r(n)}(\Qlb):\det(g)^p=1\},&\sum_i a_i\ne0,\\
 \SL_{r(n)},&\sum_i a_i=0.
 \end{cases}
\end{equation}
No conclusion about the connected group of an arbitrary or a generic
ray is being assumed in the first assertion.
\end{proposition}
\begin{proof}
Pull back Theorem~D along $x\mapsto xa$. Its Tate factor is
geometrically constant and disappears after restriction to
$\Fpbar$. A nontrivial linear Artin--Schreier sheaf on $\A^1$
has geometric monodromy group $\mu_p$, while its restriction
to $\Gm$ remains nontrivial: the point at infinity has Swan
conductor one. If the linear coefficient vanishes the sheaf is
trivial. This proves the first assertion.

Now write $G$ for the Zariski closure of the restricted geometric
fundamental-group image. If $G^0=\SL_r$, then $G$ contains
$\SL_r$ as a subgroup. For every determinant value $\omega$ in
$\det G$ and every $g\in\GL_r$ with $\det g=\omega$, choose
$g_0\in G$ with $\det g_0=\omega$. Then $gg_0^{-1}\in\SL_r\subset G$,
so $g\in G$. Thus $G=\det^{-1}(\det G)$; the first assertion
computes $\det G$ and yields \eqref{eq:conditional-ray-monodromy}.
\end{proof}

\subsection{Full monodromy of the geometric generic ray}
\label{subsec:generic-ray}

Here ``generic'' refers to the geometric generic point of the
\emph{space of directions}; it is not an assertion about every
algebraic direction.  This distinction is essential, since the
symmetric directions in Theorem~\ref{thm:symplectic} can carry a
proper symplectic monodromy group.

Put $k=\Fpbar$, $B=\Gmn{k}^{n-1}$, and $K=k(t_2,\ldots,t_n)$.
We identify $S_{n,k}$ with $\Gmn{k}\times B$ by
\begin{equation}\label{eq:relative-rays}
(x,t_2,\ldots,t_n)\longmapsto (x,t_2x,\ldots,t_nx).
\end{equation}
Let $\overline K$ be an algebraic closure of $K$.  Denote the
restricted sheaf on the geometric generic fibre by
\[
\mathcal V_{\overline\eta}
=\Delta_{(1,t_2,\ldots,t_n)}^*\cW_n\bigm|_{\Gm/\overline K}.
\]

\begin{lemma}[Normality of geometric-generic-fibre monodromy]
\label{lem:normal-generic}
Let $\mathcal F$ be any lisse $\Qlb$-sheaf on
$X=\Gmn{k}\times B$, with $k$ algebraically closed and $B$ smooth
geometrically integral.  If $G$ is its global geometric algebraic
monodromy group and $H$ the algebraic monodromy group on
$\Gm/\overline{k(B)}$, then $H$ is a closed normal algebraic subgroup
of $G$.
\end{lemma}
\begin{proof}
Set $K=k(B)$ and $L=K(x)$; choose an algebraic closure
$\overline L$ containing an algebraic closure $\overline K$ of $K$.
The geometric generic fibre has function field $\overline K(x)$.
In the exact sequence of absolute Galois groups
\[
1\longrightarrow
\operatorname{Gal}(\overline L/\overline K(x))
\longrightarrow \operatorname{Gal}(\overline L/K(x))
\longrightarrow \operatorname{Gal}(\overline K/K)
\longrightarrow1,
\]
the left-hand subgroup is normal (the restriction map on the right
is surjective, since $\overline L$ is algebraically closed).
The generic-stalk representation of $\mathcal F$ factors through
$\pi_1(X)$, while its restriction to that normal subgroup factors
through $\pi_1(\Gmn{\overline K})$.  The first image is the global
geometric image because the constant field $k$ is algebraically
closed.  Taking Zariski closures preserves normality; hence $H$
is normal in $G$.  No properness assertion or homotopy exact
sequence for the nonproper projection $X\to B$ is used.
\end{proof}

\begin{lemma}[A noncentral normal subgroup of special-linear monodromy]
\label{lem:normal-SL}
Let $r\ge2$, let $G\subset\GL_r(\Qlb)$ be an algebraic group
with $G^0=\SL_r$, and let $H\triangleleft G$ be a closed normal
algebraic subgroup containing a nonscalar element.  Then
$H^0=\SL_r$.
\end{lemma}
\begin{proof}
The connected group $H^0$ is normal in $\SL_r$.  By the
simplicity of the adjoint group $\operatorname{PGL}_r$, a connected
normal subgroup of $\SL_r$ is either trivial or all of $\SL_r$.
Suppose $H^0=1$.  For fixed $h\in H$, the commutator map
$\SL_r\to H$, $g\mapsto ghg^{-1}h^{-1}$, has connected image
containing the identity; its image is therefore contained in
$H^0=1$.  Thus every $h\in H$ commutes with $\SL_r$, whose
centralizer in $\GL_r$ consists of scalars.  This contradicts the
hypothesis, so $H^0=\SL_r$.
\end{proof}

\begin{theorem}[Large monodromy on the geometric generic ray]
\label{thm:generic-ray-SL}
Assume $n\ge2$ and $p>\max(2,n)$.  Assume the connected
ambient-monodromy theorem of \cite{EbCoh,EbMom}, namely
$G_{\mathrm{geom}}(\cW_n)^0=\SL_{r(n)}$, as recalled in
Section~\ref{subsec:family}.  Then
\begin{equation}\label{eq:generic-ray-SL}
G_{\mathrm{geom}}(\mathcal V_{\overline\eta})^0
=\SL_{r(n)}.
\end{equation}
If in addition Theorem~C is used to identify the determinant,
then the full geometric generic-ray group is
\begin{equation}\label{eq:generic-ray-full}
G_{\mathrm{geom}}(\mathcal V_{\overline\eta})
=\{g\in\GL_{r(n)}(\Qlb):\det(g)^p=1\}.
\end{equation}
These are statements over $\overline{k(t_2,\ldots,t_n)}$;
no conclusion for every $a\in(k^\times)^n$ is claimed.
\end{theorem}
\begin{proof}
Let $G$ and $H$ be the ambient and geometric generic-fibre
groups under the identification \eqref{eq:relative-rays}.
Lemma~\ref{lem:normal-generic} makes $H$ normal in $G$.
It remains to show that $H$ contains a nonscalar element.

The local Kloosterman and Fourier-transform calculations in
Sections~\ref{subsec:Kl}--\ref{subsec:proofA} are geometric and
commute with extension of algebraically closed ground fields.
Consequently the break-one formula of Theorem~A applies to the
geometric generic direction $a=(1,t_2,\ldots,t_n)$ over
$\overline K$.  All signed sums
$L_\varepsilon(a)=1+\sum_{i=2}^n\varepsilon_i\sqrt{t_i}$
are nonzero, since the $t_i$ are algebraically independent.
In particular $P(a)=2^{n-1}>0$.

For $n\ge3$, the tame part has positive dimension
\[
r(n)-2^{n-1}
=2^{n-1}-\binom n{\lfloor n/2\rfloor}>0,
\]
whereas the wild part has dimension $2^{n-1}$ and contains
nontrivial break-one characters.  A wild-inertia element thus acts
nontrivially on the wild part and trivially on the tame part;
it cannot be scalar on the whole representation.
For $n=2$, the two wild characters have Artin--Schreier
coefficients $-(1+\sqrt{t_2})^2/\lambda$ and
$-(1-\sqrt{t_2})^2/\lambda$, whose difference
$-4\sqrt{t_2}/\lambda$ is nonzero because $p$ is odd.
The characters are distinct, so again some wild-inertia element
acts nonscalarly.  This inertia lies in the geometric generic-fibre
fundamental group, hence $H$ contains a nonscalar element.

By the imported ambient theorem $G^0=\SL_{r(n)}$.
Lemma~\ref{lem:normal-SL} now gives
$H^0=\SL_{r(n)}$, proving \eqref{eq:generic-ray-SL}.
The determinant on the generic ray, by Theorem~C, is the
Artin--Schreier character
\[
\mathcal L_\psi\!\left(-2^{n-1}\lambda^{-1}
(1+t_2+\cdots+t_n)x\right).
\]
Its coefficient is nonzero in $K$, so this determinant has full
geometric image $\mu_p$ (its break at infinity equals one).
Since $H$ contains $\SL_{r(n)}$, it is the full preimage of
its determinant image, as in the proof of
Proposition~\ref{prop:ray-determinant-dichotomy}.
This gives \eqref{eq:generic-ray-full}.
\end{proof}

\begin{remark}[What the generic-fibre theorem does not say]
Theorem~\ref{thm:generic-ray-SL} concerns the geometric generic
point of the direction space and is conditional on the ambient
connected-monodromy result in the cited preprints.  It does not
by itself prove that a fixed $k$-point in a nonempty Zariski-open
set has the same connected monodromy group.  A specialization or
uniform irreducibility theorem is needed for such a stronger
claim.  On a symmetric direction, Theorem~\ref{thm:symplectic}
even forces a different ambient classical group when $n$ is even
and $r(n)>2$.
\end{remark}

\section{Characteristic cycles and conductors of parameter curves}
\label{sec:cc-curves}

Let $k=\Fpbar$ and let $\overline C/k$ be a smooth connected
projective curve.  Suppose
$f:\overline C\to(\PP^1_k)^n$ is a morphism for which
$U=f^{-1}(S_{n,k})$ is nonempty; set $D=\overline C\setminus U$,
$j_C:U\hookrightarrow\overline C$ and
$\cF_C=(f|_U)^*\cW_n$.  The sheaf $\cF_C$ is lisse of rank $r(n)$.
We use $[T^*_{\overline C}\overline C]$ for the zero section and
$[T_y^*\overline C]$ for the cotangent fibre at $y$.

\begin{theorem}[Characteristic-cycle and conductor identity on every parameter curve]
\label{thm:cc-any-curve}
For every such $f$, without a genericity or transversality
assumption, one has
\begin{equation}\label{eq:cc-curve}
CC(j_{C!}\cF_C)
=-r(n)[T^*_{\overline C}\overline C]
-\sum_{y\in D}\bigl(r(n)+\Swan_y(\cF_C)\bigr)
[T_y^*\overline C],
\end{equation}
where the points of $D$ are geometric points.  In particular,
\begin{equation}\label{eq:gos-any-curve}
\chi_c(U,\cF_C)
=r(n)\bigl(2-2g(\overline C)-\#D\bigr)
-\sum_{y\in D}\Swan_y(\cF_C).
\end{equation}
Let $X=(\PP^1_k)^n$ and $K_X=j_{X!}\cW_n$,
where $j_X:S_{n,k}\hookrightarrow X$.  If $f$ is properly
$SS(K_X)$-transversal in the sense of Saito, then
\begin{equation}\label{eq:cc-pullback}
 f^!CC(K_X)=CC(j_{C!}\cF_C),
\end{equation}
where $f^!$ is Saito's signed microlocal pull-back of cycles.
Thus \eqref{eq:cc-curve} identifies the local conductor
coefficients with those of $f^!CC(K_X)$ wherever this
transversality hypothesis is satisfied.
\end{theorem}
\begin{proof}
The curve formula for characteristic cycles of extension by
zero is the one-dimensional conductor formula: at a missing
point $y$, the total dimension of the local ramification
representation is $\rank\cF_C+\Swan_y(\cF_C)$.
The coefficient of the zero section is the negative generic rank,
and the coefficient of the cotangent fibre at $y$ is the negative
total dimension.  This is the curve case of Saito's characteristic-
cycle formula \cite{SaitoCC} (and the classical Artin-conductor
formula).  The Grothendieck--Ogg--Shafarevich theorem then yields
\eqref{eq:gos-any-curve}; equivalently, intersect
\eqref{eq:cc-curve} with the zero section and use the index formula.

By base change for extension by zero along an open immersion,
$f^*j_{X!}\cW_n\cong j_{C!}\cF_C$.
Under proper $SS(K_X)$-transversality, Saito's pull-back theorem
\cite[Theorem~4.4]{SaitoCC} identifies the two characteristic
cycles exactly, using the sign convention incorporated in $f^!$.
This proves \eqref{eq:cc-pullback}.
\end{proof}

\subsection{The generic conormal coefficients at coordinate boundaries}
\label{subsec:generic-cc-boundary}

We now make the distinction between the inertia conductor of a
lisse sheaf and the characteristic cycle of its \emph{extension by
zero} explicit.  In this subsection, the coefficient attached to a
divisor means the coefficient of the closure of the conormal bundle
to its \emph{geometric generic point}; no statement about the
intersection of two or more boundary divisors is implicit.
We retain Saito's cohomological convention
$CC(\mathcal L)=(-1)^n\operatorname{rank}(\mathcal L)
[T_X^*X]$ for a lisse sheaf on a smooth $n$-fold.

\begin{proposition}[The tame zero-coordinate boundary]
\label{prop:cc-zero-coord}
Assume the coordinate-boundary unipotence statement of
Theorem~\ref{thm:recollection} and the Jordan profile of
Theorem~\ref{thm:coordinate-jordan}.  Let
$X=(\PP^1_k)^n$ and $K=j_{X!}\cW_n$.
For each $i$, write $D_i^0=\{h_i=0\}\subset X$.  At the
geometric generic point of $D_i^0$, one has
\[
\Swan_{D_i^0}(\cW_n)=0,\qquad
\dim\cW_n^{I_{D_i^0}}=r(n)-r(n-1),\qquad
\operatorname{Art}_{D_i^0}(\cW_n)=r(n-1).
\]
Nevertheless the coefficient of the conormal
$[T^*_{D_i^0}X]$ in $CC(K)$ is
\begin{equation}\label{eq:cc-zero-coord}
\boxed{\operatorname{coeff}_{[T^*_{D_i^0}X]}CC(K)=(-1)^n r(n).}
\end{equation}
\end{proposition}
\begin{proof}
By the imported coordinate-boundary calculation, the inertia
operator $T_i$ is tame and unipotent, with precisely
$r(n-1)$ blocks of size two and $r(n)-2r(n-1)$ blocks
of size one.  A unipotent Jordan block of either positive size
has a one-dimensional fixed space.  Hence
\begin{align*}
\dim\ker(T_i-1)
 &=r(n-1)+r(n)-2r(n-1)\\
 &=r(n)-r(n-1).
\end{align*}
There is no wild ramification; thus the Swan conductor is zero
and the Artin conductor, defined as
$\operatorname{rank}-\dim V^I+\Swan$, is $r(n-1)$.

The remaining assertion is a different computation.  Near the
generic point of $D_i^0$, the boundary consists of a single smooth
divisor and $K$ has zero stalk along that divisor.
For a tame lisse sheaf of rank $r$ on the complement of a smooth
divisor, the local tame characteristic-cycle formula for extension
by zero reads
\[
CC(j_!\mathcal L)=(-1)^n r
\bigl([T_X^*X]+[T_D^*X]\bigr)
\quad\text{generically along }D;
\]
see \cite{SaitoCC}.  It depends on the rank, not on the
number of nontrivial tame Jordan blocks.  Applying it with
$r=r(n)$ gives \eqref{eq:cc-zero-coord}.
In particular, replacing the conormal coefficient by the tame
inertia drop would confuse $j_!\cW_n$ with a different extension
across the boundary.
\end{proof}

\begin{proposition}[The generic infinity divisor in rank two]
\label{prop:cc-infty-rank-two}
Let $n=2$ and $p>2$. At the geometric generic point of
$D_i^\infty=\{h_i=\infty\}\subset (\PP^1_k)^2$, the sheaf
$\cW_2$ is isoclinic of break $1$ and rank $2$.  In particular,
\[
\boxed{\Swan_{D_i^\infty}(\cW_2)=2,\qquad
\operatorname{Art}_{D_i^\infty}(\cW_2)=4,\qquad
\operatorname{coeff}_{[T^*_{D_i^\infty}X]}CC(j_{X!}\cW_2)=4.}
\]
At a zero-coordinate generic divisor the corresponding triple
is $(0,1,2)$.
\end{proposition}
\begin{proof}
By Theorem~\ref{thm:recollection}(v),
\begin{equation}\label{eq:cc-ranktwo-model}
\cW_2\cong
\Lpsi{-\frac{h_1+h_2}{\lambda}}
\otimes\bigl[(h_1,h_2)\mapsto h_1h_2/\lambda^2\bigr]^*\Kl
\otimes\Qlb(-1).
\end{equation}
Work at the geometric generic point of $D_1^\infty$;
fix $h_2=c\in k^\times$, or equivalently take $c$
to be transcendental over the ground field.  With $t=h_1$,
the factor $\Lpsi{-(t+c)/\lambda}$ has a single upper
ramification break $1$ at $t=\infty$.  The rank-two
Kloosterman sheaf $\Kl(ct/\lambda^2)$ has two breaks
$\tfrac12$ at infinity by the standard local description
\cite{Katz-GKM}.

Here is the precise tensor-product step.  If a rank-one
inertia representation has a single break $b$ and a second
representation has all breaks $b'<b$, then the tensor product
has break $b$ on its entire space.  Indeed the larger-break
character remains nontrivial on an upper-ramification subgroup
on which the smaller-break representation is already trivial.
In the present situation $1>\tfrac12$, so both breaks of
\eqref{eq:cc-ranktwo-model} become $1$.  Their sum is the
Swan conductor $2$, and the positive-break representation has
no inertia-fixed vectors.  The Artin conductor is therefore
$2-0+2=4$.

To pass from this integer to the characteristic-cycle
coefficient, observe that the dominant Artin--Schreier phase
$-t/\lambda$ is independent of $c$.  Hence its refined wild
characteristic form points in the normal direction to
$D_1^\infty$ at the generic point; there is no extra tangent
direction there.  The generic smooth-divisor formula of
\cite{SaitoCC} consequently assigns the total dimension
$\operatorname{rank}+\Swan=2+2$ to its conormal, with
sign $(-1)^2=1$.  Permuting the two coordinates proves the
assertion for $D_2^\infty$.

Finally, near $t=0$ the Artin--Schreier factor is unramified,
while $\Kl(ct/\lambda^2)$ has one tame unipotent block
of size two.  Its invariant space has dimension one,
so its Swan and Artin conductors are respectively $0$ and $1$.
Proposition~\ref{prop:cc-zero-coord} gives the generic conormal
coefficient $2$.
\end{proof}

The identities above are \emph{general conductor formulas}, not
an explicit evaluation of $\Swan_y(\cF_C)$ for every map $f$.
The latter requires microlocal information about the boundary
of $\cW_n$ not contained in Theorem~A.  Nevertheless, the ray
calculation gives an effective expression for every curve contained
in a fixed ray, including wildly ramified maps into that ray.

\begin{definition}[Artin--Schreier reduced pole order]
\label{def:AS-order}
For $g\in k((t))$, let $\operatorname{asord}(g)$ be the Swan
conductor of $\mathcal L_\psi(g)$ over $k((t))$.  Equivalently,
reduce the polar part of $g$ modulo elements $h^p-h$; then
$\operatorname{asord}(g)$ is the largest positive pole exponent
prime to $p$ with nonzero coefficient in a reduced representative,
or zero if no such exponent remains.
\end{definition}

\begin{theorem}[An exact conductor formula for arbitrary scalar-profile curves]
\label{thm:cc-scalar}
Let $a\in(k^\times)^n$ be fixed and let
$\varphi:\overline C\to\PP^1$ be a nonconstant morphism.
Let $U=\varphi^{-1}(\Gm)$ and define
$f|_U(y)=\varphi(y)a$.
For every $y\in\overline C\setminus U$,
\begin{equation}\label{eq:swan-scalar}
\Swan_y(\cF_C)
=\begin{cases}
0,&\varphi(y)=0,\\[2pt]
\displaystyle\sum_{\substack{\varepsilon\in\mathsf E\\L_\varepsilon(a)\ne0}}
\operatorname{asord}_y\!\left(-\frac{L_\varepsilon(a)^2}{\lambda}
\varphi\right),&\varphi(y)=\infty.
\end{cases}
\end{equation}
In particular, if $y$ is a pole of order $e_y$ prime to $p$,
then
\begin{equation}\label{eq:swan-scalar-tame}
\Swan_y(\cF_C)=e_yP(a).
\end{equation}
For arbitrary pole orders, the reduced Artin--Schreier pole
orders in \eqref{eq:swan-scalar} give an exact algorithm.
Substituting this expression in \eqref{eq:cc-curve} and
\eqref{eq:gos-any-curve} computes the full characteristic cycle
and compactly supported Euler characteristic of the scalar-profile
curve.
\end{theorem}
\begin{proof}
At a zero of $\varphi$, the sheaf $\Delta_a^*\cW_n$ is tame
at $x=0$ by Proposition~\ref{prop:ray-model} and
Theorem~\ref{thm:recollection}(iv).  Pull-back of a tame local
representation along any finite extension of local fields remains
tame, so the Swan conductor is zero.

At a pole of $\varphi$, pull back the inertia decomposition
of Theorem~A(i).  The original break-zero representation remains
tame after pull-back.  Each break-one summand is an
Artin--Schreier character
$\mathcal L_\psi(-L_\varepsilon(a)^2x/\lambda)$
multiplied by a tame quadratic Kummer character.
After $x=\varphi(y)$ the Kummer twist stays tame, while the
Artin--Schreier character has Swan conductor precisely
$\operatorname{asord}_y(-L_\varepsilon(a)^2\varphi/\lambda)$.
Swan conductors are additive under direct sums, proving
\eqref{eq:swan-scalar}.  If $e_y$ is prime to $p$, the leading
pole term of each nonzero scalar multiple has exponent $e_y$
prime to $p$, so no Artin--Schreier coboundary can remove it.
Every one of the $P(a)$ summands then has conductor $e_y$;
this proves \eqref{eq:swan-scalar-tame}.
\end{proof}

\begin{example}[Why nonradial curves genuinely differ]
\label{ex:nonradial}
Let $n=2$ and $p>3$.  On $\Gm$ take
$f(x)=(x,x^2)$.  By the explicit rank-two formula in
Theorem~\ref{thm:recollection}(v),
\[
f^*\cW_2\cong
\mathcal L_\psi\!\left(-\frac{x+x^2}{\lambda}\right)
\otimes \bigl([x\mapsto x^3/\lambda^2]^*\Kl\bigr)(-1).
\]
At $x=\infty$ the first factor has break $2$.  The Kloosterman
factor has both breaks $3/2$ because the map $x\mapsto x^3$
is tamely ramified of index $3$ and $p>3$.
The unequal-break tensor lemma therefore makes the rank-two
representation isoclinic of break $2$, giving
\[
\boxed{\Swan_\infty(f^*\cW_2)=4.}
\]
By contrast, along a nondegenerate scalar ray for $n=2$ the Swan
conductor is $P(a)=2$.  Thus applying the radial formula to
an arbitrary nonradial parameter curve is not valid.
\end{example}

\subsection{A complete conductor calculation on constant-product curves}
\label{subsec:constant-product}

We next treat a nonradial family for which a full local
calculation, including the contribution at every puncture, is
possible.  The decisive feature is not a transversality
assumption but the fact that the argument of the Kloosterman
sheaf in \eqref{eq:cc-ranktwo-model} is constant.

\begin{theorem}[Conductors on every constant-product parameter curve]
\label{thm:cc-hyperbola}
Let $k=\Fpbar$, let $U/k$ be a smooth geometrically connected
curve, and let $u\in\Gamma(U,\mathcal O_U)^\times$ be a
nonconstant invertible regular function.  Fix $c\in k^\times$
and define
\[
\gamma:U\longrightarrow S_{2,k},\qquad
z\longmapsto\bigl(u(z),c/u(z)\bigr).
\]
Let $\overline C$ be the smooth proper compactification of $U$
and $D=\overline C\setminus U$.  Put
\[
g=-\frac{u+c/u}{\lambda},\qquad
V_c=\Kl(c/\lambda^2)\otimes\Qlb(-1),
\qquad e_x=\Swan_x(\Lpsi{g})\quad(x\in D).
\]
Here $V_c$ denotes the geometrically constant rank-two sheaf
with this fibre. Then there is an isomorphism of lisse sheaves
\begin{equation}\label{eq:cc-hyperbola-sheaf}
\boxed{\gamma^*\cW_2\cong\Lpsi{g}\otimes V_c.}
\end{equation}
At every $x\in D$ its local conductors are exactly
\begin{equation}\label{eq:cc-hyperbola-local}
\boxed{\Swan_x(\gamma^*\cW_2)=2e_x,
\qquad
\operatorname{Art}_x(\gamma^*\cW_2)=
\begin{cases}
2+2e_x,&e_x>0,\\
0,&e_x=0.
\end{cases}}
\end{equation}
Consequently
\begin{equation}\label{eq:cc-hyperbola-euler}
\boxed{\chi_c(U,\gamma^*\cW_2)
=2\chi_c(U,\Qlb)-2\sum_{x\in D}e_x.}
\end{equation}
These numbers are effectively determined from the Laurent
expansion of $g$ at $x$: after subtracting functions of the
form $b^p-b$, write its reduced polar part as
$\sum_{d>0,\,p\nmid d}a_dz^{-d}$.  Then
\[
e_x=\begin{cases}
\max\{d:a_d\ne0\},&\text{if the reduced polar part is nonzero},\\
0,&\text{otherwise}.
\end{cases}
\]
No assumption that the original pole order of $u$ is prime
to $p$ is required in this theorem.
\end{theorem}
\begin{proof}
Substitute $h_1=u$ and $h_2=c/u$ into
\eqref{eq:cc-ranktwo-model}.  The product $h_1h_2$ is the
constant $c$, so the pulled-back Kloosterman sheaf is the
constant local system with fibre $\Kl(c/\lambda^2)$.
The additive phase becomes $-(u+c/u)/\lambda=g$.
The Tate twist is itself geometrically constant.  This proves
\eqref{eq:cc-hyperbola-sheaf} without taking semisimplifications
or passing to a trace identity.

Fix $x\in D$, choose a local parameter $z$ at $x$, and
restrict to the punctured strict henselization at $x$.
Geometric inertia acts trivially on $V_c$ and acts on
$\Lpsi{g}$ by a rank-one character $\chi_x$ of $p$-power
order.  Thus its action on $\gamma^*\cW_2$ is the direct sum
of two copies of $\chi_x$.  The Swan conductor is additive on
direct sums, proving the first identity in
\eqref{eq:cc-hyperbola-local}.  If $e_x>0$, the character
$\chi_x$ is nontrivial on wild inertia; its fixed space is
zero, and the same holds for two copies.  The Artin conductor
is then $\operatorname{rank}-\dim V^{I_x}+\Swan=2+2e_x$.
If $e_x=0$, the Artin--Schreier character is unramified:
a rank-one character of $p$-power order cannot be nontrivial
and tame in characteristic $p$, and equivalently its reduced
polar part vanishes.  The whole representation is therefore
unramified and its Artin conductor is zero.

For completeness we justify the local algorithm rather than
merely prescribing it.  In $k((z))$ write the polar part of
$g$ as a finite sum $\sum_{d>0}a_dz^{-d}$.  Suppose its
highest exponent $d$ is divisible by $p$.  Since $k$ is
algebraically closed, choose $b\in k$ with $b^p=a_d$ and
replace $g$ by
\[
g-\bigl((bz^{-d/p})^p-bz^{-d/p}\bigr).
\]
This operation does not change the associated Artin--Schreier
character and eliminates its highest pole.  Repeating the
operation terminates and yields only exponents prime to $p$.
Conversely, a Laurent series $v$ with a pole of order $m>0$
has $v^p-v$ of pole order $pm$, so a nonzero reduced leading
term of degree prime to $p$ cannot be removed by a further
Artin--Schreier coboundary.  The standard upper-ramification
calculation for a rank-one Artin--Schreier character
\cite{Serre-LF} identifies that maximal reduced pole order
with $\Swan_x(\Lpsi{g})=e_x$.

Finally, apply the compact-support formula
\eqref{eq:gos-any-curve} to the rank-two lisse sheaf
$\gamma^*\cW_2$.  The Swan contribution at $x$ is $2e_x$,
so one obtains \eqref{eq:cc-hyperbola-euler} exactly.
Notice that neither a tame inertia drop nor an Artin
conductor is subtracted a second time.
\end{proof}

\begin{corollary}[Nonradial monomial curves and their cohomology]
\label{cor:cc-hyperbola-monomial}
Let $A,B\in k^\times$ and let $m\ge1$ satisfy $p\nmid m$.
For
\[
\gamma_m:\Gm\longrightarrow S_{2,k},\qquad
\gamma_m(t)=(At^m,Bt^{-m}),
\]
one has
\begin{equation}\label{eq:cc-monomial-conductors}
\boxed{\Swan_0(\gamma_m^*\cW_2)=
\Swan_\infty(\gamma_m^*\cW_2)=2m,
\qquad
\operatorname{Art}_0=\operatorname{Art}_\infty=2m+2.}
\end{equation}
Moreover,
\begin{equation}\label{eq:cc-monomial-cohomology}
\boxed{H_c^0(\Gm,\gamma_m^*\cW_2)=
H_c^2(\Gm,\gamma_m^*\cW_2)=0,\qquad
\dim H_c^1(\Gm,\gamma_m^*\cW_2)=4m.}
\end{equation}
\end{corollary}
\begin{proof}
The product of the two coordinates is $c=AB$, hence the
function of Theorem~\ref{thm:cc-hyperbola} is
\[
g(t)=-\frac{At^m+Bt^{-m}}{\lambda}.
\]
It has a pole of order $m$ at $t=0$ with leading coefficient
$-B/\lambda$, and a pole of order $m$ at $t=\infty$
(with local parameter $z=t^{-1}$) with leading coefficient
$-A/\lambda$.  Both are nonzero and $p\nmid m$, so the
poles are already Artin--Schreier reduced.  Thus
$e_0=e_\infty=m$, and
\eqref{eq:cc-monomial-conductors} follows from
\eqref{eq:cc-hyperbola-local}.

By \eqref{eq:cc-hyperbola-euler} and
$\chi_c(\Gm,\Qlb)=0$, the Euler characteristic is $-4m$.
The sheaf is lisse on the connected nonproper curve $\Gm$,
so no nonzero compactly supported section can exist:
$H_c^0=0$.  Furthermore, since the Artin--Schreier character
has positive Swan conductor at both ends, it is nontrivial as
a geometric fundamental-group representation.  Its tensor
product with a geometrically constant rank-two vector space
has no invariant vectors; by Poincar\'e duality the top
compactly supported cohomology $H_c^2$ vanishes.  Therefore
$-\dim H_c^1=\chi_c=-4m$, proving
\eqref{eq:cc-monomial-cohomology}.
\end{proof}

\begin{proposition}[The characteristic cycle has no interior escape component]
\label{prop:no-interior-escape-cycle}
Put $X=(\PP^1_k)^n$, $S=S_{n,k}$, and $K=j_{X!}\cW_n$.
With the standard cohomological sign convention for characteristic
cycles of sheaves in degree zero, one has
\begin{equation}\label{eq:cc-interior}
  CC(K)|_{T^*S}=(-1)^n r(n)[T^*_SS].
\end{equation}
In particular, if $Z\subset S$ is any divisor defined by the
vanishing of a critical-value or critical-point escape polynomial,
there is no additional irreducible component of $CC(K)$ supported
over the conormal bundle of $Z$.
\end{proposition}
\begin{proof}
By Theorem~\ref{thm:recollection}, $\cW_n$ is lisse of rank
$r(n)$ everywhere on $S$.  The singular support of a lisse sheaf
on a smooth scheme is the zero section, and the characteristic
cycle of such a sheaf placed in cohomological degree zero is
$(-1)^{\dim S}$ times its rank times that section
\cite{SaitoCC}.  Restricting $j_{X!}\cW_n$ to $S$ proves
\eqref{eq:cc-interior}.  The conormal of a proper divisor in
$S$ contains nonzero covectors over its smooth generic point,
so it cannot be a component of this restriction.  A divisor of
\emph{directions} where the ray conductor jumps is therefore not,
by that fact alone, a singular-support divisor of $\cW_n$ on $S$.
\end{proof}

\begin{proposition}[Coordinate-wise pole orders do not add]
\label{prop:nonadditive-pole-orders}
Let $n=2$ and $p>3$.  A proposed formula of the form
\begin{equation}\label{eq:false-coordinate-formula}
  \Swan_y(f^*\cW_2)
    \stackrel{?}{=}2\sum_{i:\,f_i(y)=\infty}
                  (-\operatorname{ord}_y f_i)
\end{equation}
for all maps from punctured curves to $S_2$ is false, even when
both coordinates tend to infinity. More precisely:
\begin{enumerate}[label=\textup{(\roman*)}]
\item For the scalar ray $f(x)=(x,2x)$,
      $\Swan_\infty(f^*\cW_2)=2$, whereas the right-hand side
      of \eqref{eq:false-coordinate-formula} is $4$.
\item For the nonradial curve $f(x)=(x,x^2)$,
      $\Swan_\infty(f^*\cW_2)=4$, whereas that right-hand side
      is $6$.
\end{enumerate}
Consequently, generic codimension-one ramification data cannot
be added along coordinates to determine the Swan conductor at
an intersection of boundary components.
\end{proposition}
\begin{proof}
For (i), take $a=(1,2)$.  Choosing $s_1,s_2\in k^\times$
with $s_1^2=1$ and $s_2^2=2$ gives
$s_1+s_2\ne0$ and $s_1-s_2\ne0$; hence $P(a)=2$.
Theorem~A(ii), or Theorem~\ref{thm:cc-scalar} with
$\varphi(x)=x$, yields $\Swan_\infty=2$.
Both coordinate functions have simple poles, so the proposed
coordinate-wise expression is $2(1+1)=4$.
For (ii), the calculation in Example~\ref{ex:nonradial}
shows that the Artin--Schreier factor has break $2$ and the
pulled-back rank-two Kloosterman factor has breaks $3/2$.
Their tensor product is isoclinic of break $2$, of rank $2$;
thus its Swan conductor is $4$.  The coordinate pole orders are
$1$ and $2$, giving $2(1+2)=6$ on the right of
\eqref{eq:false-coordinate-formula}.  Both are strict
contradictions.
\end{proof}

\begin{remark}[Extension by zero versus tame inertia drop]
\label{rem:extension-conductor-distinction}
In the Grothendieck--Ogg--Shafarevich identity
\eqref{eq:gos-any-curve} the correction to the compactly
supported Euler characteristic consists of the Swan terms,
\emph{not} the sum of tame inertia drops.  In the characteristic
cycle of $j_{C!}\cF_C$, the coefficient at an omitted point is
instead minus $\rank\cF_C+\Swan_y(\cF_C)$, as in
\eqref{eq:cc-curve}; a tame unipotent block is not the same
as a trivial representation.  For example, a nontrivial tame
rank-one Kummer sheaf on $\Gm$ has $\chi_c=0$, despite
having inertia drop one at both $0$ and $\infty$.
The higher-dimensional coefficients of $CC(K)$ require a
separate analysis of the singular support and refined boundary
ramification; neither their values nor the relevant conormal
components are determined by the two codimension-one
Jordan profiles.
\end{remark}

\begin{remark}[Precise remaining microlocal problem]
\label{rem:global-cc-open}
Theorem~\ref{thm:cc-any-curve} is unconditional as a curve
conductor identity and gives an ambient characteristic-cycle
pull-back identity only for properly non-characteristic curves.
Theorem~\ref{thm:cc-scalar} explicitly computes every local
conductor for curves of the form $f(y)=a\varphi(y)$.
Neither statement determines $CC(j_{X!}\cW_n)$ on all boundary
strata of $X$, or evaluates all Swan conductors for arbitrary
nonradial curves tangent to the singular support.  A closed
formula in terms only of the critical values of Theorem~A for
such general curves remains unproved.  In particular, a
full higher-codimension boundary calculation cannot be inferred
from the generic-ray theorem.
\end{remark}

\paragraph{Additional Euler and boundary calculations.}
The global Euler-characteristic identity and the conditional
higher-codimension recursion are recorded in
Appendix~\ref{app:cc-euler}.  Their distinction from the local
conductor formulas proved in this section is essential: the
recursion requires an independent hypothesis on the singular
support of the extension by zero.

\section{Cohomological degrees and hypotheses used}
\label{sec:proof-dependencies}

\begin{proposition}[What the Euler characteristic determines]
\label{prop:virtual-degree}
Let $\mathcal F$ be any lisse sheaf on $\Gm/\Fq$, with
$H_c^0((\Gm\otimes\Fqbar),\mathcal F)=0$. Put
$h_i=\dim H_c^i((\Gm\otimes\Fqbar),\mathcal F)$. Then
\[
L(\Gm,\mathcal F,T)=
\frac{\det(1-T\Frob_q\mid H_c^1)}
 {\det(1-T\Frob_q\mid H_c^2)},\qquad
h_1-h_2=-\chi_c((\Gm\otimes\Fqbar),\mathcal F).
\]
Consequently the Euler characteristic determines the virtual
(cohomological) degree $h_1-h_2$, but does not itself prove that
$L(\Gm,\mathcal F,T)$ is a polynomial of that degree.
If $H_c^2=0$, then the polynomiality and degree assertions follow.
\end{proposition}
\begin{proof}
The compactly supported trace formula gives the first identity.
The definition of the Euler--Poincar\'e characteristic and $h_0=0$
give the second. If the denominator is absent, the numerator is a
polynomial of degree $h_1$. Without its absence one may not discard
its eigenvalues by formal cancellation.
\end{proof}

\begin{remark}[Normalizations in the local Fourier step]
\label{rem:fourier-steps-audited}
The transform used here is
$\FT(P)=R\pr_{z!}(\pr_y^*P\otimes\mathcal L_\psi(yz))[1]$.
For $P=j_{y!}\mathcal K_a[1]$ its ordinary cohomology
$\mathcal H^{-1}(P)$ is $j_{y!}\mathcal K_a$. The shifted Fourier
transform $Q=\FT(P)$ is perverse; on $z\ne0$ it has the form
$\mathcal F_a[1]$. In the nearby/vanishing-cycle triangle
$i^*Q\to R\Psi(Q)\to R\Phi(Q)\to$,
one has $R\Psi(Q)=\mathcal F_a(0)[1]$. The relevant local Fourier
module is $R^{-1}\Phi(Q)$, not $R^0\Phi(Q)$.
Lemma~\ref{lem:vanishing} exhibits the specialization sequence
explicitly; its outer terms are unramified as \emph{inertia}
representations. Hence the positive-break parts are identified but
the entire tame representation need not be determined.
The substitution $z=\lambda/x$ carries
$\mathcal L_\psi(b/z)$ to $\mathcal L_\psi(bx/\lambda)$;
this is the sign and scale convention used in Theorem~B.
An inertia-level isomorphism does not fix an arithmetic Frobenius
scalar, which is why Section~\ref{sec:arith} is logically separate.
\end{remark}

\begin{remark}[Geometric inputs from companion manuscripts]
\label{rem:logical-dependency-imported}
The construction, rank, weight, purity and geometric irreducibility
of $\mathcal W_n$, the lower-weight filtration of $\mathcal H_n$,
and unipotence at coordinate hyperplanes are used as theorems from
\cite{EbCoh,EbMom}; their proofs are not independently reproduced
here. The reflection pairing is proved by the cup product and
Poincar\'e duality in Proposition~\ref{prop:reflection}; the tame-end
arithmetic determinant is calculated in
Proposition~\ref{prop:tame-end}. The two sign determinations use
Theorem~C for the definition of the geometrically constant scalar.
The connected-monodromy group is invoked only in the explicitly
restricted characteristic range of Corollary~E.
Any use of the complete-ray bound must also cite
Proposition~\ref{prop:global-euler}(ii), rather than only its
Euler-characteristic part (i).
\end{remark}

\section{Finite-field and algebraic verifications}
\label{sec:numerics}

The verifications in this section concern finite identities and do not
replace the $\ell$-adic arguments. In particular the finite
Taylor expansion of a rational $L$-function is not used to infer exact
cohomology dimensions. The tests below are limited to identities reproducible using the
accompanying code. They do not establish claims about dimensions
of cohomology groups or vanishing cycles.

For $p$ an odd prime define the integral character encoding
$\mathbb Z[T]/(T^p-1)$, in which the coefficient of $T^j$ records the
number of occurrences of the additive phase $j\bmod p$. Equality
in this ring entails equality after evaluation at every $p$th root of
unity. The accompanying file \texttt{verify\_rays\_journal.py}
uses only standard Python arithmetic and performs the following tests:
\begin{enumerate}[label=(\alph*)]
\item it enumerates constrained pairs $u_1,u_2\in\Fp^\times$
for $p\in\{3,5,7,11\}$ and all $h_1,h_2\in\Fp^\times$,
checking the explicit identity
\[
\Sigma_2(h_1,h_2)=
\psi\!\left(-\frac{h_1+h_2}{\lambda}\right)
K_p\!\left(\frac{h_1h_2}{\lambda^2}\right)-1,
\qquad\lambda=1,
\]
by exact equality of phase-count vectors;
\item it checks the first critical-value power sum and the
four-variable identity
$\mathscr D_4=(e_1^2-4e_2)^2-64e_4$
by symbolic expansion over $\mathbb Z$;
\item it checks the finite character-orthogonality count
of Theorem~\ref{thm:det-distribution} for $p=3,5,7,11$
and $n=2,3,4$ using exact integer counts.
\end{enumerate}
The computation is independent of any claim about ranks of
vanishing cycles or arithmetic semisimplicity. The outputs and
precise counts are written to standard output by the script.

\appendix
\section{Global Euler characteristics and conditional boundary recursions}
\label{app:cc-euler}

The generic conormal calculations above do not by themselves fix the
coefficients supported over intersections of coordinate divisors.  We
first establish an Euler-characteristic identity independently of any
hypothesis on the higher-codimension singular support.  We then record
exactly what follows from proper push-forward \emph{if} the singular
support has no additional non-toric components.  The latter condition
is a substantive microlocal hypothesis; it will not be inferred from
the constancy of a Swan conductor along a divisor.

For $m\ge1$ set
\[
 B_m(T)=\sum_{k=1}^{\lfloor m/2\rfloor}\nu_{m,k}T^k,
 \qquad B_1(T)=0.
\]

\begin{lemma}[Summation over the last parameter]
\label{lem:cc-parameter-summation}
For $n\ge2$, every finite extension $\Fq/\Fp$, and
$h'=(h_1,\ldots,h_{n-1})\in(\Fq^\times)^{n-1}$,
\begin{equation}\label{eq:parameter-sigma-sum}
 \sum_{z\in\Fq^\times}\Sigma_{n,q}(h',z)
 =(-1)^n+\Sigma_{n-1,q}(h').
\end{equation}
Consequently, with $C_n(T)=T+B_{n-1}(T)+(T-1)B_n(T)$,
\begin{equation}\label{eq:parameter-W-trace}
 \sum_{z\in\Fq^\times}
 \Tr(\Frob_{(h',z)}\mid\cW_n)
 =-\Tr(\Frob_{h'}\mid\cW_{n-1})-C_n(q).
\end{equation}
\end{lemma}
\begin{proof}
Fix $u'=(u_1,\ldots,u_{n-1})\in(\Fq^\times)^{n-1}$.
For each $u_n\in\Fq^\times$, the linear constraint
$\langle h',u'\rangle+zu_n=-\lambda$ determines a unique
$z\in\Fq$; this value is nonzero precisely when
$\langle h',u'\rangle\ne-\lambda$.  Thus summing first over
$z\ne0$ removes exactly the vectors $u'$ satisfying the
$(n-1)$-dimensional constraint.  Because
$\sum_{t\in\Fq^\times}\psi_q(t^{-1})=-1$ and
$\sum_{u'\in(\Fq^\times)^{n-1}}
\psi_q(\sum_{i<n}u_i^{-1})=(-1)^{n-1}$, we obtain
\[
 \sum_{z\ne0}\Sigma_{n,q}(h',z)
 =-(-1)^{n-1}+\Sigma_{n-1,q}(h'),
\]
which is \eqref{eq:parameter-sigma-sum}.
The trace formula \eqref{eq:trace} reads
$\Tr(\Frob_h\mid\cW_m)=(-1)^{m-1}q\Sigma_{m,q}(h)-B_m(q)$.
Substituting \eqref{eq:parameter-sigma-sum} and using the same
formula for $m=n-1$ gives
\begin{align*}
 \sum_{z\ne0}\Tr(\Frob_{(h',z)}\mid\cW_n)
 &=-q-\Tr(\Frob_{h'}\mid\cW_{n-1})
   -B_{n-1}(q)-(q-1)B_n(q),
\end{align*}
which is \eqref{eq:parameter-W-trace}.  All equalities hold over
\emph{every} finite extension, not only over the prime field.
\end{proof}

\begin{theorem}[Exact Euler characteristic of the parameter torus]
\label{thm:cc-global-euler}
For every $n\ge1$, assuming the sheaf construction and trace formula
of Theorem~\ref{thm:recollection},
\begin{equation}\label{eq:cc-global-euler}
 \boxed{\chi_c(S_n\otimes\Fpbar,\cW_n)=(-1)^n.}
\end{equation}
This identity requires no description of the characteristic cycle of
$j_{X!}\cW_n$ in codimension at least two.
\end{theorem}
\begin{proof}
Write $C_n(T)=\sum_j\gamma_{n,j}T^j$ with integral
coefficients.  Summing \eqref{eq:parameter-W-trace} over all
$h'\in S_{n-1}(\Fq)$ and applying the compact-support trace
formula gives, for every $q=p^s$,
\begin{align*}
 \Tr\!\left(\Frob_q\mid R\Gamma_c(S_n,\cW_n)\right)
 &=-\Tr\!\left(\Frob_q\mid R\Gamma_c(S_{n-1},\cW_{n-1})\right)\\
 &\quad-\sum_j\gamma_{n,j}\,
  \Tr\!\left(\Frob_q\mid
 R\Gamma_c(S_{n-1},\Qlb(-j))\right).
\end{align*}
Here traces on complexes are alternating traces; the expression
$q^j\#S_{n-1}(\Fq)$ is exactly the alternating Frobenius trace
of the constant Tate-twisted sheaf $\Qlb(-j)$ on $S_{n-1}$.
Since this equality holds for every power of Frobenius, Newton's
identities (or Brauer--Nesbitt applied to virtual representations)
identify the associated virtual Frobenius representations.  Taking
their dimensions yields
\[
 \chi_c(S_n,\cW_n)
 =-\chi_c(S_{n-1},\cW_{n-1})
 -\sum_j\gamma_{n,j}\chi_c(S_{n-1},\Qlb(-j)).
\]
For $n-1\ge1$, the last Euler characteristics all vanish:
$\chi_c(\Gm^{n-1},\Qlb(-j))=
\chi_c(\Gm,\Qlb)^{n-1}=0$, by Kunneth and the fact that
$\chi_c(\Gm,\Qlb)=0$.  Therefore
$\chi_c(S_n,\cW_n)=-\chi_c(S_{n-1},\cW_{n-1})$.
Finally $\cW_1\cong\Lpsi{-h/\lambda}\otimes\Qlb(-1)$ by
Theorem~\ref{thm:recollection}(v), and the
Grothendieck--Ogg--Shafarevich formula gives
$\chi_c(\Gm,\cW_1)=-1$, because the Swan conductor at
infinity is one and that at zero is zero.  Induction completes
the proof.
\end{proof}

\begin{proposition}[The rank-two projection, including its Tate twists]
\label{prop:cc-ranktwo-projection}
For this proposition take $n=2$ and $\lambda=1$.
Put $F=\Lpsi{-h_1-h_2}\otimes m^*\Kl$, so that
$\cW_2=F(-1)$, and set $K=j_{X!}F$ on
$X=(\PP^1)^2$.  Let $\pi:X\to\PP^1$ be projection on $h_1$,
and write $j_1:\Gm\hookrightarrow\PP^1$.  In the Grothendieck
group of constructible complexes on $\PP^1$ one has
\begin{equation}\label{eq:cc-projection-K0}
 \boxed{[R\pi_*K]=-[j_{1!}\Qlb(-1)]
          -[j_{1!}\Lpsi{-h_1}].}
\end{equation}
In particular,
\begin{equation}\label{eq:cc-ranktwo-direct-image}
 \boxed{CC(R\pi_*K)=2[T^*_{\PP^1}\PP^1]
       +2[T^*_0\PP^1]+3[T^*_{\infty}\PP^1],}
\end{equation}
and $\chi_c(S_2,F)=1$.
\end{proposition}
\begin{proof}
We use the normalized Kloosterman trace
$\Tr(\Frob_t\mid\Kl)=-\sum_{x\in\Fq^\times}
\psi_q(x^{-1}+tx)$.
For $h_1\in\Fq^\times$ a direct orthogonality calculation gives
\begin{align*}
 \sum_{h_2\ne0}\Tr(\Frob_{(h_1,h_2)}\mid F)
 &=-\psi_q(-h_1)\sum_{x\ne0}\psi_q(x^{-1})
             \sum_{h_2\ne0}\psi_q((h_1x-1)h_2)\\
 &=-\psi_q(-h_1)\bigl(q\psi_q(h_1)+1\bigr)\\
 &=-q-\psi_q(-h_1).
\end{align*}
The middle identity uses
$\sum_{h_2\ne0}\psi_q(th_2)=q-1$ for $t=0$
and $=-1$ otherwise.  Since $\pi$ is proper and $K$ is
extension by zero, its derived direct image is zero at
$h_1=0$ and $h_1=\infty$.  Equality of the traces over all
finite extensions, followed by Chebotarev and Brauer--Nesbitt
on $\Gm$, proves \eqref{eq:cc-projection-K0}.
On a smooth curve, the characteristic cycle of $j_{1!}$ of a
rank-one lisse sheaf is the negative zero section minus the
cotangent fibres, with coefficients $1+\Swan$ at the missing
points.  The constant sheaf has Swan zero at both ends;
$\Lpsi{-h_1}$ has Swan zero at $0$ and one at infinity.
Characteristic cycles are additive on distinguished triangles,
so the two negative signs in \eqref{eq:cc-projection-K0}
produce \eqref{eq:cc-ranktwo-direct-image}.
Finally
\[
 \chi_c(S_2,F)
 =-\chi_c(\Gm,\Qlb(-1))
  -\chi_c(\Gm,\Lpsi{-h_1})=0-(-1)=1.
\]
In contrast, the corresponding arithmetic identity for
$\cW_2=F(-1)$ has the additional Tate twist:
$[R\pi_*j_!\cW_2]=-[j_{1!}\Qlb(-2)]
-[j_{1!}\Lpsi{-h_1}(-1)]$.
\end{proof}

\begin{proposition}[The rank-two zero--zero conormal coefficient]
\label{prop:cc-zerozero-corner}
Let $p>2$, let $X=(\PP^1)^2$, and write $j_X:\Gm^2\hookrightarrow X$.
Assume the rank-two identification of Theorem~\ref{thm:recollection}(v).
Then the coefficient of $[T^*_{(0,0)}X]$ in
$CC(j_{X!}\cW_2)$ is $2$.
\end{proposition}
\begin{proof}
The assertion is local in a strict henselian neighbourhood of
$(0,0)$, where we write $x=h_1$ and $y=h_2$.  By
\eqref{eq:cc-ranktwo-model}, up to a Tate twist the restriction
to $xy\ne0$ is
\[
\mathcal L_\psi(-(x+y)/\lambda)\otimes
  [xy/\lambda^2]^*\Kl.
\]
The Artin--Schreier factor is lisse across the divisor $xy=0$.
The rank-two Kloosterman sheaf is tamely ramified at $0$;
therefore its pullback is tame along the strict normal-crossings
divisor $\{x=0\}\cup\{y=0\}$.  The tame normal-crossings
formula for the characteristic cycle of the extension by zero
\cite{SaitoCC}, in dimension two, gives on this neighbourhood
\[
CC(j_!\mathcal F)=2\bigl(
[T_X^*X]+[T^*_{\{x=0\}}X]+[T^*_{\{y=0\}}X]
+[T^*_{(0,0)}X]\bigr).
\]
The coefficient claimed in the proposition is thus $2$.
No information about the three other corners or the singular
support along infinity has been used.
\end{proof}

\begin{proposition}[The four rank-two corner coefficients under a stated microlocal hypothesis]
\label{prop:cc-ranktwo-corners-conditional}
Retain the notation of Proposition~\ref{prop:cc-ranktwo-projection}.
Assume in addition that
\begin{equation*}
 SS(K)\subseteq\bigcup_{Z\text{ a coordinate stratum of }X}
          \overline{T^*_Z X}.
 \tag{$\mathrm{NT}_2$}
\end{equation*}
The first coefficient $c_{00}=2$ is unconditional by
Proposition~\ref{prop:cc-zerozero-corner}.  Under the additional
hypothesis $\mathrm{NT}_2$, the conormal coefficients at the four corners
$(0,0)$, $(0,\infty)$, $(\infty,0)$, $(\infty,\infty)$
are, respectively,
\begin{equation}\label{eq:cc-corners-conditional}
 \boxed{(c_{00},c_{0\infty},c_{\infty0},c_{\infty\infty})
            =(2,4,4,7).}
\end{equation}
The coefficient on the zero section is $2$; those on the four
coordinate divisors are $(2,2,4,4)$, as computed in
Propositions~\ref{prop:cc-zero-coord} and
\ref{prop:cc-infty-rank-two}.  No assertion that
$\mathrm{NT}_2$ has been established is part of this proposition.
\end{proposition}
\begin{proof}
By Proposition~\ref{prop:cc-zerozero-corner}, $c_{00}=2$ without
assuming $\mathrm{NT}_2$.  The remaining three coefficients
require the hypothesis specified in the proposition.

Under $\mathrm{NT}_2$ the singular support consists of the
toric conormals, and its image under the cotangent correspondence
for $\pi$ has dimension at most one.  Hence Saito's proper
direct-image theorem \cite{SaitoDirect} applies:
$CC(R\pi_*K)=\pi_!CC(K)$.
The elementary toric conormal push-forward, evaluated at the
cotangent fibres over $h_1=0$ and $h_1=\infty$, respectively,
reads
\begin{align*}
 \operatorname{coeff}_{T^*_0\PP^1}\pi_!CC(K)
 &=c_{00}+c_{0\infty}-2c_{D_1^0},\\
 \operatorname{coeff}_{T^*_{\infty}\PP^1}\pi_!CC(K)
 &=c_{\infty0}+c_{\infty\infty}-2c_{D_1^\infty}.
\end{align*}
The terms $-2c_{D_1^0}$ and $-2c_{D_1^\infty}$
come from the Euler characteristic $\chi(\PP^1)=2$ of
the free coordinate.  By \eqref{eq:cc-ranktwo-direct-image},
the left sides are $2$ and $3$.  Using
$c_{D_1^0}=2$, $c_{D_1^\infty}=4$, and $c_{00}=2$, we find
$c_{0\infty}=4$ and $c_{\infty0}+c_{\infty\infty}=11$.
The involution interchanging $h_1$ and $h_2$ preserves $F$;
thus $c_{\infty0}=c_{0\infty}=4$ and
$c_{\infty\infty}=7$.

The index formula provides an independent check, \emph{not}
a substitute for $\mathrm{NT}_2$:
\[
 \chi_c(S_2,F)=2\cdot4-2(2+2+4+4)
                    +(2+4+4+7)=1,
\]
in agreement with Proposition~\ref{prop:cc-ranktwo-projection}.
The argument also shows why a ray's Swan conductor at a
corner cannot be identified with the corner's conormal
coefficient: the push-forward involves all relevant conormals.
\end{proof}

\begin{lemma}[Telescoping of the boundary multiplicities]
\label{lem:cc-kappa}
For $m\ge1$ put
$\nu_{m,k}=\binom mk-\binom m{k-1}$ and
$\beta_m=1+\sum_{k=1}^{\lfloor m/2\rfloor}\nu_{m,k}$.
Then
\[
\beta_m=\binom m{\lfloor m/2\rfloor}=2^m-r(m).
\]
This numerical identity places no restriction on non-toric
components of a characteristic cycle.
\end{lemma}
\begin{proof}
Put $d=\lfloor m/2\rfloor$.  The definition telescopes:
\[
1+\sum_{k=1}^{d}\left(\binom mk-\binom m{k-1}\right)
=1+\binom md-\binom m0=\binom md.
\]
The last equality in the assertion is the definition of $r(m)$.
\end{proof}

\begin{proposition}[Toric-recursion criterion and its Euler-characteristic test]
\label{prop:cc-toric-recursion-conditional}
For $m\le n$ write $X_m=(\PP^1)^m$ and
$K_m=j_{X_m!}\cW_m$.  Suppose, in all dimensions $m\le n$,
that $SS(K_m)$ is a union of coordinate-stratum conormals and
that the coefficient at every stratum with no infinity
coordinates is $(-1)^m r(m)$.
Then the coefficients of $CC(K_m)$ are independent of the
number of zero coordinates: writing them as
$(-1)^m f_m(s)$, where $s$ is the number of coordinates
at infinity, one has
\begin{equation}\label{eq:cc-toric-recurrence}
 f_1(0)=1,\quad f_1(1)=2,\quad f_m(0)=r(m),\qquad
 \boxed{f_m(s+1)=f_m(s)+f_{m-1}(s)
     +\binom{m-1}{\lfloor(m-1)/2\rfloor}.}
\end{equation}
In particular, the first rows predicted under these
hypotheses are
\begin{equation}\label{eq:cc-toric-table}
\begin{array}{c|rrrrr}
 m\backslash s &0&1&2&3&4\\\hline
 1&1&2&&&\\
 2&2&4&7&&\\
 3&5&9&15&24&\\
 4&10&18&30&48&75
\end{array}
\end{equation}
For every $m$ these numbers satisfy
\begin{equation}\label{eq:cc-toric-euler-check}
 \sum_{s=0}^m(-1)^s\binom ms f_m(s)=(-1)^m.
\end{equation}
Neither \eqref{eq:cc-toric-recurrence} nor the table is
asserted as an unconditional calculation of the global
characteristic cycle when the stated microlocal hypotheses
are not known.
\end{proposition}
\begin{proof}
Let $\pi_m:S_m\to S_{m-1}$ forget one parameter.
Lemma~\ref{lem:cc-parameter-summation}, for every finite
extension, and Chebotarev--Brauer--Nesbitt give in the
Grothendieck group on $S_{m-1}$
\begin{equation}\label{eq:cc-K0-recursion}
 [R\pi_{m!}\cW_m]
 =-[\cW_{m-1}]-\sum_j\gamma_{m,j}[\Qlb(-j)],
 \qquad C_m(T)=\sum_j\gamma_{m,j}T^j.
\end{equation}
The virtual constant sheaf on the right has total rank
\[
 \sum_j\gamma_{m,j}=C_m(1)
 =1+B_{m-1}(1)
 =\binom{m-1}{\lfloor(m-1)/2\rfloor},
\]
by Lemma~\ref{lem:cc-kappa}.

Extend the projection to the proper map $X_m\to X_{m-1}$.
The assumed toric singular support satisfies Saito's
dimension hypothesis for the proper direct-image theorem
\cite{SaitoDirect}.  The microlocal push-forward of the
coordinate conormals gives, at a fixed base stratum, the
coefficient of the two endpoint strata in the forgotten
coordinate minus twice the coefficient of the free-coordinate
stratum.  Combining this rule with \eqref{eq:cc-K0-recursion}
and the constant-sheaf rank just calculated yields, for disjoint
$I_0,I_\infty$ not containing the forgotten coordinate $j$,
\begin{align}\label{eq:cc-toric-stratum-relation}
 &a_m(I_0\cup\{j\},I_\infty)
 +a_m(I_0,I_\infty\cup\{j\})-2a_m(I_0,I_\infty)\notag\\
 &\hspace{32mm}=a_{m-1}(I_0,I_\infty)
 +\binom{m-1}{\lfloor(m-1)/2\rfloor},
\end{align}
where $a_m$ is the \emph{signed-normalized} conormal coefficient
$(-1)^m\operatorname{coeff}_{T^*_Z X_m}CC(K_m)$.

The assumed zero-only coefficients give
$a_m(I_0,\varnothing)=r(m)$.  Starting with $s=0$,
use \eqref{eq:cc-toric-stratum-relation} and induction on
$s=|I_\infty|$ (and on $m$) to see that the value is
independent of $I_0$, and that it obeys the recurrence
\eqref{eq:cc-toric-recurrence}.  For $m=1$ the sheaf
$\cW_1=\Lpsi{-h/\lambda}(-1)$ has a tame zero and
Swan conductor one at infinity; its extension by zero has
coefficients $-1$ and $-2$, respectively, giving the initial
row.  The displayed table follows by arithmetic.

Finally set $A_m=\sum_s(-1)^s\binom ms f_m(s)$.
The binomial identity
\[
 A_m=-\sum_{s=0}^{m-1}(-1)^s\binom{m-1}s
               \bigl(f_m(s+1)-f_m(s)\bigr)
\]
and \eqref{eq:cc-toric-recurrence} yield $A_m=-A_{m-1}$
for $m\ge2$, since the alternating binomial sum of a constant
vanishes.  As $A_1=1-2=-1$, we get
$A_m=(-1)^m$.  This agrees with, but does not independently
establish the microlocal hypotheses of,
Theorem~\ref{thm:cc-global-euler}.
\end{proof}

\begin{remark}[The unresolved higher-codimension issue]
\label{rem:cc-higher-open}
Condition $\mathrm{NT}_2$ has not been deduced here from the
codimension-one Swan conductors.  In dimensions $n\ge3$ it is
not known from the quoted local-monodromy results whether
additional microlocal components occur above intersections of
infinity divisors.  Vanishing of a subleading stationary-phase
coefficient, for example on a locus of the form
$h_1=h_2=\infty$, $h_3=h_4$, suggests a test locus but does
not prove a new conormal component.  Such a claim would
require a calculation of the relevant vanishing cycles and
refined characteristic directions.  Thus the conditional
corner values and the recursion must not be read as a proof
of a closed characteristic-cycle formula for all $n$.
\end{remark}

\section*{Acknowledgments}
This work was carried out independently. I am grateful to those whose support, patience, and encouragement made it possible to continue the project over a long period. Much of the research was done alongside other professional responsibilities, and the arguments were revisited repeatedly as the manuscript developed. The final form reflects that sustained process of reading, checking, and revision. I remain solely responsible for the mathematical content and for any errors that remain.

All mathematical arguments, proofs, and conclusions in this article are the sole work of the author. Computational tools were used only for the numerical checks explicitly described in the text; these checks are not essential to the proofs and are provided for verification. The author has reviewed all results and assumes full responsibility for the correctness and integrity of the work. AI-assisted tools, if used, were employed exclusively for language editing and formatting; they did not contribute to any mathematical content.

\end{document}